\documentclass{amsart}

\usepackage{amsmath}
\usepackage{breqn}
\usepackage{amssymb}

\usepackage{hyperref}

\newcommand{\erase}[1]{}

\newtheorem{theorem}{Theorem}[section]

\newtheorem{proposition}[theorem]{Proposition}
\newtheorem{corollary}[theorem]{Corollary}

\newtheorem{_definition}[theorem]{Definition}
\newenvironment{definition}{\begin{_definition}\rm}{\end{_definition}}

\newtheorem{_remark}[theorem]{\it Remark}
\newenvironment{remark}{\begin{_remark}\rm}{\end{_remark}}

\newtheorem{_example}[theorem]{Example}
\newenvironment{example}{\begin{_example}\rm}{\end{_example}}

\numberwithin{equation}{section}
\numberwithin{table}{section}

\newcommand{\F}{\mathord{\mathbb F}}
\renewcommand{\P}{\mathord{\mathbb  P}}
\newcommand{\Q}{\mathord{\mathbb  Q}}
\newcommand{\R}{\mathord{\mathbb R}}

\newcommand{\Z}{\mathord{\mathbb Z}}

\newcommand{\HHH}{\mathord{\mathcal H}}

\newcommand{\LLL}{\mathord{\mathcal L}}

\newcommand{\OOO}{\mathord{\mathcal O}}
\newcommand{\PPP}{\mathord{\mathcal P}}

\newcommand{\RRR}{\mathord{\mathcal R}}

\newcommand{\WWW}{\mathord{\mathcal W}}

\newcommand{\maprightsp}[1]{\; \smash{\mathop{\; \longrightarrow \; }\limits\sp{#1}}\; }
\newcommand{\isom}{\mathbin{\,\raise -.6pt\rlap{$\to$}\raise 3.5pt \hbox{\hskip .3pt$\mathord{\sim}$}\,}}

\newcommand{\set}[2]{\{\; {#1} \; \mid \; {#2} \;  \}}
\newcommand{\shortset}[2]{\{ {#1} \,|\, {#2}   \}}

\newcommand{\gen}[1]{\langle {#1}  \rangle}

\newcommand{\tensor}{\otimes}

\newcommand{\pprime}{\prime\prime}
\newcommand{\sprime}{\sp\prime}

\newcommand{\spprime}{\sp{\prime\prime}}
\newcommand{\spcirc}{\sp{\mathord{\circ}}}

\newcommand{\sperp}{\sp{\perp}}

\newcommand{\dual}{\sp{\vee}}

\newcommand{\inv}{\sp{-1}}

\newcommand{\GL}{\mathord{\mathrm {GL}}}

\newcommand{\PGU}{\mathord{\mathrm {PGU}}}
\newcommand{\PSU}{\mathord{\mathrm {PSU}}}
\newcommand{\PGL}{\mathord{\mathrm {PGL}}}
\newcommand{\OG}{\mathord{\mathrm {O}}}

\newcommand{\Aut}{\operatorname{\mathrm {Aut}}\nolimits}

\newcommand{\cone}{\mathord{\PPP}}
\newcommand{\tR}{\tensor\R}
\newcommand{\tQ}{\tensor\Q}
\newcommand{\Exc}{\mathord{\mathrm {Exc}}}
\newcommand{\intL}[1]{(#1)_L}
\newcommand{\intM}[1]{(#1)_M}
\newcommand{\intS}[1]{(#1)_S}
\newcommand{\intT}[1]{(#1)_T}
\newcommand{\cham}{D}

\newcommand{\LR}{\mathord{\it{LR}}}
\newcommand{\Fr}{\mathord{\phi}}
\newcommand{\wallrs}{\widetilde\WWW}
\newcommand{\wallrsz}[1]{\widetilde{W}_{#1}}
\newcommand{\wall}{\WWW}
\newcommand{\chamsz}{\cham_{S0}}

\newcommand{\quand}{\quad\textrm{and}\quad}

\newcommand{\sqrmo}{i}
\newcommand{\shrink}{\hskip -5pt}
\newcommand{\myvskipmat}{\vskip 3pt}
\newcommand{\myboxnumb}[1]{#1}
\newcommand{\numclass}[1]{[#1]}

\newcommand{\mcirc}[2]{\put#1{\circle{1.2}}\put#1{\hbox{\raise -3pt \hbox{\hskip -5.1pt \hbox{\tiny$#2$}}}}}

\newcommand{\surj}{\mathbin{\to \hskip -7pt \to}}
\newcommand{\Sing}{\operatorname{\mathrm {Sing}}\nolimits}

\newcommand{\mystruthd}[2]{\phantom{\hbox{\vrule  height #1 depth #2}}}

\newcommand{\lineFQ}[1]{
 \ifnum#1<0 {wrong} \fi
 \ifnum#1=0 {wrong} \fi
 \ifnum#1>112 {wrong} \fi
 \ifnum#1=1 {\textstyle\{ \, w+(1+\sqrmo)\, = \, x+(1+\sqrmo)\, y=0\}}\fi
  \ifnum#1=2 {\textstyle\{ \, w+(1+\sqrmo)\, = \, x+(1-\sqrmo)\, y=0\}}\fi
  \ifnum#1=3 {\textstyle\{ \, w+(1+\sqrmo)\, = \, x-(1-\sqrmo)\, y=0\}}\fi
  \ifnum#1=4 {\textstyle\{ \, w+(1+\sqrmo)\, = \, x-(1+\sqrmo)\, y=0\}}\fi
  \ifnum#1=5 {\textstyle\{ \, w+(1-\sqrmo)\, = \, x+(1+\sqrmo)\, y=0\}}\fi
  \ifnum#1=6 {\textstyle\{ \, w+(1-\sqrmo)\, = \, x+(1-\sqrmo)\, y=0\}}\fi
  \ifnum#1=7 {\textstyle\{ \, w+(1-\sqrmo)\, = \, x-(1-\sqrmo)\, y=0\}}\fi
  \ifnum#1=8 {\textstyle\{ \, w+(1-\sqrmo)\, = \, x-(1+\sqrmo)\, y=0\}}\fi
  \ifnum#1=9 {\textstyle\{ \, w-(1-\sqrmo)\, = \, x+(1+\sqrmo)\, y=0\}}\fi
  \ifnum#1=10 {\textstyle\{ \, w-(1-\sqrmo)\, = \, x+(1-\sqrmo)\, y=0\}}\fi
  \ifnum#1=11 {\textstyle\{ \, w-(1-\sqrmo)\, = \, x-(1-\sqrmo)\, y=0\}}\fi
  \ifnum#1=12 {\textstyle\{ \, w-(1-\sqrmo)\, = \, x-(1+\sqrmo)\, y=0\}}\fi
  \ifnum#1=13 {\textstyle\{ \, w-(1+\sqrmo)\, = \, x+(1+\sqrmo)\, y=0\}}\fi
  \ifnum#1=14 {\textstyle\{ \, w-(1+\sqrmo)\, = \, x+(1-\sqrmo)\, y=0\}}\fi
  \ifnum#1=15 {\textstyle\{ \, w-(1+\sqrmo)\, = \, x-(1-\sqrmo)\, y=0\}}\fi
  \ifnum#1=16 {\textstyle\{ \, w-(1+\sqrmo)\, = \, x-(1+\sqrmo)\, y=0\}}\fi
  \ifnum#1=17 {\textstyle\{ \, w+\sqrmo\, y+\sqrmo\, = \, x+\sqrmo\, y-\sqrmo\, =0\}}\fi
  \ifnum#1=18 {\textstyle\{ \, w+\sqrmo\, y+\sqrmo\, = \, x-\sqrmo\, y+\sqrmo\, =0\}}\fi
  \ifnum#1=19 {\textstyle\{ \, w+\sqrmo\, y+\sqrmo\, = \, x+\, y-1\, =0\}}\fi
  \ifnum#1=20 {\textstyle\{ \, w+\sqrmo\, y+\sqrmo\, = \, x-\, y+1\, =0\}}\fi
  \ifnum#1=21 {\textstyle\{ \, w+\sqrmo\, y-\sqrmo\, = \, x+\sqrmo\, y+\sqrmo\, =0\}}\fi
  \ifnum#1=22 {\textstyle\{ \, w+\sqrmo\, y-\sqrmo\, = \, x-\sqrmo\, y-\sqrmo\, =0\}}\fi
  \ifnum#1=23 {\textstyle\{ \, w+\sqrmo\, y-\sqrmo\, = \, x+\, y+1\, =0\}}\fi
  \ifnum#1=24 {\textstyle\{ \, w+\sqrmo\, y-\sqrmo\, = \, x-\, y-1\, =0\}}\fi
  \ifnum#1=25 {\textstyle\{ \, w+\sqrmo\, y+1\, = \, x+\sqrmo\, y-1\, =0\}}\fi
  \ifnum#1=26 {\textstyle\{ \, w+\sqrmo\, y+1\, = \, x-\sqrmo\, y+1\, =0\}}\fi
  \ifnum#1=27 {\textstyle\{ \, w+\sqrmo\, y+1\, = \, x+\, y+\sqrmo\, =0\}}\fi
  \ifnum#1=28 {\textstyle\{ \, w+\sqrmo\, y+1\, = \, x-\, y-\sqrmo\, =0\}}\fi
  \ifnum#1=29 {\textstyle\{ \, w+\sqrmo\, y-1\, = \, x+\sqrmo\, y+1\, =0\}}\fi
  \ifnum#1=30 {\textstyle\{ \, w+\sqrmo\, y-1\, = \, x-\sqrmo\, y-1\, =0\}}\fi
  \ifnum#1=31 {\textstyle\{ \, w+\sqrmo\, y-1\, = \, x+\, y-\sqrmo\, =0\}}\fi
  \ifnum#1=32 {\textstyle\{ \, w+\sqrmo\, y-1\, = \, x-\, y+\sqrmo\, =0\}}\fi
  \ifnum#1=33 {\textstyle\{ \, w-\sqrmo\, y+\sqrmo\, = \, x+\sqrmo\, y+\sqrmo\, =0\}}\fi
  \ifnum#1=34 {\textstyle\{ \, w-\sqrmo\, y+\sqrmo\, = \, x-\sqrmo\, y-\sqrmo\, =0\}}\fi
  \ifnum#1=35 {\textstyle\{ \, w-\sqrmo\, y+\sqrmo\, = \, x+\, y+1\, =0\}}\fi
  \ifnum#1=36 {\textstyle\{ \, w-\sqrmo\, y+\sqrmo\, = \, x-\, y-1\, =0\}}\fi
  \ifnum#1=37 {\textstyle\{ \, w-\sqrmo\, y-\sqrmo\, = \, x+\sqrmo\, y-\sqrmo\, =0\}}\fi
  \ifnum#1=38 {\textstyle\{ \, w-\sqrmo\, y-\sqrmo\, = \, x-\sqrmo\, y+\sqrmo\, =0\}}\fi
  \ifnum#1=39 {\textstyle\{ \, w-\sqrmo\, y-\sqrmo\, = \, x+\, y-1\, =0\}}\fi
  \ifnum#1=40 {\textstyle\{ \, w-\sqrmo\, y-\sqrmo\, = \, x-\, y+1\, =0\}}\fi
  \ifnum#1=41 {\textstyle\{ \, w-\sqrmo\, y+1\, = \, x+\sqrmo\, y+1\, =0\}}\fi
  \ifnum#1=42 {\textstyle\{ \, w-\sqrmo\, y+1\, = \, x-\sqrmo\, y-1\, =0\}}\fi
  \ifnum#1=43 {\textstyle\{ \, w-\sqrmo\, y+1\, = \, x+\, y-\sqrmo\, =0\}}\fi
  \ifnum#1=44 {\textstyle\{ \, w-\sqrmo\, y+1\, = \, x-\, y+\sqrmo\, =0\}}\fi
  \ifnum#1=45 {\textstyle\{ \, w-\sqrmo\, y-1\, = \, x+\sqrmo\, y-1\, =0\}}\fi
  \ifnum#1=46 {\textstyle\{ \, w-\sqrmo\, y-1\, = \, x-\sqrmo\, y+1\, =0\}}\fi
  \ifnum#1=47 {\textstyle\{ \, w-\sqrmo\, y-1\, = \, x+\, y+\sqrmo\, =0\}}\fi
  \ifnum#1=48 {\textstyle\{ \, w-\sqrmo\, y-1\, = \, x-\, y-\sqrmo\, =0\}}\fi
  \ifnum#1=49 {\textstyle\{ \, w+\, y+\sqrmo\, = \, x+\sqrmo\, y+1\, =0\}}\fi
  \ifnum#1=50 {\textstyle\{ \, w+\, y+\sqrmo\, = \, x-\sqrmo\, y-1\, =0\}}\fi
  \ifnum#1=51 {\textstyle\{ \, w+\, y+\sqrmo\, = \, x+\, y-\sqrmo\, =0\}}\fi
  \ifnum#1=52 {\textstyle\{ \, w+\, y+\sqrmo\, = \, x-\, y+\sqrmo\, =0\}}\fi
  \ifnum#1=53 {\textstyle\{ \, w+\, y-\sqrmo\, = \, x+\sqrmo\, y-1\, =0\}}\fi
  \ifnum#1=54 {\textstyle\{ \, w+\, y-\sqrmo\, = \, x-\sqrmo\, y+1\, =0\}}\fi
  \ifnum#1=55 {\textstyle\{ \, w+\, y-\sqrmo\, = \, x+\, y+\sqrmo\, =0\}}\fi
  \ifnum#1=56 {\textstyle\{ \, w+\, y-\sqrmo\, = \, x-\, y-\sqrmo\, =0\}}\fi
  \ifnum#1=57 {\textstyle\{ \, w+\, y+1\, = \, x+\sqrmo\, y-\sqrmo\, =0\}}\fi
  \ifnum#1=58 {\textstyle\{ \, w+\, y+1\, = \, x-\sqrmo\, y+\sqrmo\, =0\}}\fi
  \ifnum#1=59 {\textstyle\{ \, w+\, y+1\, = \, x+\, y-1\, =0\}}\fi
  \ifnum#1=60 {\textstyle\{ \, w+\, y+1\, = \, x-\, y+1\, =0\}}\fi
  \ifnum#1=61 {\textstyle\{ \, w+\, y-1\, = \, x+\sqrmo\, y+\sqrmo\, =0\}}\fi
  \ifnum#1=62 {\textstyle\{ \, w+\, y-1\, = \, x-\sqrmo\, y-\sqrmo\, =0\}}\fi
  \ifnum#1=63 {\textstyle\{ \, w+\, y-1\, = \, x+\, y+1\, =0\}}\fi
  \ifnum#1=64 {\textstyle\{ \, w+\, y-1\, = \, x-\, y-1\, =0\}}\fi
  \ifnum#1=65 {\textstyle\{ \, w+(1+\sqrmo)\, y= \, x+(1+\sqrmo)\, =0\}}\fi
  \ifnum#1=66 {\textstyle\{ \, w+(1+\sqrmo)\, y= \, x+(1-\sqrmo)\, =0\}}\fi
  \ifnum#1=67 {\textstyle\{ \, w+(1+\sqrmo)\, y= \, x-(1-\sqrmo)\, =0\}}\fi
  \ifnum#1=68 {\textstyle\{ \, w+(1+\sqrmo)\, y= \, x-(1+\sqrmo)\, =0\}}\fi
  \ifnum#1=69 {\textstyle\{ \, w+(1-\sqrmo)\, y= \, x+(1+\sqrmo)\, =0\}}\fi
  \ifnum#1=70 {\textstyle\{ \, w+(1-\sqrmo)\, y= \, x+(1-\sqrmo)\, =0\}}\fi
  \ifnum#1=71 {\textstyle\{ \, w+(1-\sqrmo)\, y= \, x-(1-\sqrmo)\, =0\}}\fi
  \ifnum#1=72 {\textstyle\{ \, w+(1-\sqrmo)\, y= \, x-(1+\sqrmo)\, =0\}}\fi
  \ifnum#1=73 {\textstyle\{ \, w-\, y+\sqrmo\, = \, x+\sqrmo\, y-1\, =0\}}\fi
  \ifnum#1=74 {\textstyle\{ \, w-\, y+\sqrmo\, = \, x-\sqrmo\, y+1\, =0\}}\fi
  \ifnum#1=75 {\textstyle\{ \, w-\, y+\sqrmo\, = \, x+\, y+\sqrmo\, =0\}}\fi
  \ifnum#1=76 {\textstyle\{ \, w-\, y+\sqrmo\, = \, x-\, y-\sqrmo\, =0\}}\fi
  \ifnum#1=77 {\textstyle\{ \, w-\, y-\sqrmo\, = \, x+\sqrmo\, y+1\, =0\}}\fi
  \ifnum#1=78 {\textstyle\{ \, w-\, y-\sqrmo\, = \, x-\sqrmo\, y-1\, =0\}}\fi
  \ifnum#1=79 {\textstyle\{ \, w-\, y-\sqrmo\, = \, x+\, y-\sqrmo\, =0\}}\fi
  \ifnum#1=80 {\textstyle\{ \, w-\, y-\sqrmo\, = \, x-\, y+\sqrmo\, =0\}}\fi
  \ifnum#1=81 {\textstyle\{ \, w-\, y+1\, = \, x+\sqrmo\, y+\sqrmo\, =0\}}\fi
  \ifnum#1=82 {\textstyle\{ \, w-\, y+1\, = \, x-\sqrmo\, y-\sqrmo\, =0\}}\fi
  \ifnum#1=83 {\textstyle\{ \, w-\, y+1\, = \, x+\, y+1\, =0\}}\fi
  \ifnum#1=84 {\textstyle\{ \, w-\, y+1\, = \, x-\, y-1\, =0\}}\fi
  \ifnum#1=85 {\textstyle\{ \, w-\, y-1\, = \, x+\sqrmo\, y-\sqrmo\, =0\}}\fi
  \ifnum#1=86 {\textstyle\{ \, w-\, y-1\, = \, x-\sqrmo\, y+\sqrmo\, =0\}}\fi
  \ifnum#1=87 {\textstyle\{ \, w-\, y-1\, = \, x+\, y-1\, =0\}}\fi
  \ifnum#1=88 {\textstyle\{ \, w-\, y-1\, = \, x-\, y+1\, =0\}}\fi
  \ifnum#1=89 {\textstyle\{ \, w-(1-\sqrmo)\, y= \, x+(1+\sqrmo)\, =0\}}\fi
  \ifnum#1=90 {\textstyle\{ \, w-(1-\sqrmo)\, y= \, x+(1-\sqrmo)\, =0\}}\fi
  \ifnum#1=91 {\textstyle\{ \, w-(1-\sqrmo)\, y= \, x-(1-\sqrmo)\, =0\}}\fi
  \ifnum#1=92 {\textstyle\{ \, w-(1-\sqrmo)\, y= \, x-(1+\sqrmo)\, =0\}}\fi
  \ifnum#1=93 {\textstyle\{ \, w-(1+\sqrmo)\, y= \, x+(1+\sqrmo)\, =0\}}\fi
  \ifnum#1=94 {\textstyle\{ \, w-(1+\sqrmo)\, y= \, x+(1-\sqrmo)\, =0\}}\fi
  \ifnum#1=95 {\textstyle\{ \, w-(1+\sqrmo)\, y= \, x-(1-\sqrmo)\, =0\}}\fi
  \ifnum#1=96 {\textstyle\{ \, w-(1+\sqrmo)\, y= \, x-(1+\sqrmo)\, =0\}}\fi
  \ifnum#1=97 {\textstyle\{ \, w+(1+\sqrmo)\, x= \, y+(1+\sqrmo)\, =0\}}\fi
  \ifnum#1=98 {\textstyle\{ \, w+(1+\sqrmo)\, x= \, y+(1-\sqrmo)\, =0\}}\fi
  \ifnum#1=99 {\textstyle\{ \, w+(1+\sqrmo)\, x= \, y-(1-\sqrmo)\, =0\}}\fi
  \ifnum#1=100 {\textstyle\{ \, w+(1+\sqrmo)\, x= \, y-(1+\sqrmo)\, =0\}}\fi
  \ifnum#1=101 {\textstyle\{ \, w+(1-\sqrmo)\, x= \, y+(1+\sqrmo)\, =0\}}\fi
  \ifnum#1=102 {\textstyle\{ \, w+(1-\sqrmo)\, x= \, y+(1-\sqrmo)\, =0\}}\fi
  \ifnum#1=103 {\textstyle\{ \, w+(1-\sqrmo)\, x= \, y-(1-\sqrmo)\, =0\}}\fi
  \ifnum#1=104 {\textstyle\{ \, w+(1-\sqrmo)\, x= \, y-(1+\sqrmo)\, =0\}}\fi
  \ifnum#1=105 {\textstyle\{ \, w-(1-\sqrmo)\, x= \, y+(1+\sqrmo)\, =0\}}\fi
  \ifnum#1=106 {\textstyle\{ \, w-(1-\sqrmo)\, x= \, y+(1-\sqrmo)\, =0\}}\fi
  \ifnum#1=107 {\textstyle\{ \, w-(1-\sqrmo)\, x= \, y-(1-\sqrmo)\, =0\}}\fi
  \ifnum#1=108 {\textstyle\{ \, w-(1-\sqrmo)\, x= \, y-(1+\sqrmo)\, =0\}}\fi
  \ifnum#1=109 {\textstyle\{ \, w-(1+\sqrmo)\, x= \, y+(1+\sqrmo)\, =0\}}\fi
  \ifnum#1=110 {\textstyle\{ \, w-(1+\sqrmo)\, x= \, y+(1-\sqrmo)\, =0\}}\fi
  \ifnum#1=111 {\textstyle\{ \, w-(1+\sqrmo)\, x= \, y-(1-\sqrmo)\, =0\}}\fi
  \ifnum#1=112 {\textstyle\{ \, w-(1+\sqrmo)\, x= \, y-(1+\sqrmo)\, =0\}}\fi
 }

\newcommand{\thevFQ}[1]{
 \ifnum#1<0 {wrong} \fi
 \ifnum#1=0 {wrong} \fi
 \ifnum#1>112 {wrong} \fi
 \ifnum#1=1 {[
 \myboxnumb{1}, \myboxnumb{0}, \myboxnumb{0}, \myboxnumb{0}, \myboxnumb{0}, \myboxnumb{0}, \myboxnumb{0}, \myboxnumb{0}, \myboxnumb{0}, \myboxnumb{0}, \myboxnumb{0}, \myboxnumb{0}, \myboxnumb{0}, \myboxnumb{0}, \myboxnumb{0}, \myboxnumb{0}, \myboxnumb{0}, \myboxnumb{0}, \myboxnumb{0}, \myboxnumb{0}, \myboxnumb{0}, \myboxnumb{0}]_S} \fi
 \ifnum#1=2 {[
 \myboxnumb{0}, \myboxnumb{1}, \myboxnumb{0}, \myboxnumb{0}, \myboxnumb{0}, \myboxnumb{0}, \myboxnumb{0}, \myboxnumb{0}, \myboxnumb{0}, \myboxnumb{0}, \myboxnumb{0}, \myboxnumb{0}, \myboxnumb{0}, \myboxnumb{0}, \myboxnumb{0}, \myboxnumb{0}, \myboxnumb{0}, \myboxnumb{0}, \myboxnumb{0}, \myboxnumb{0}, \myboxnumb{0}, \myboxnumb{0}]_S} \fi
 \ifnum#1=3 {[
 \myboxnumb{0}, \myboxnumb{0}, \myboxnumb{1}, \myboxnumb{0}, \myboxnumb{0}, \myboxnumb{0}, \myboxnumb{0}, \myboxnumb{0}, \myboxnumb{0}, \myboxnumb{0}, \myboxnumb{0}, \myboxnumb{0}, \myboxnumb{0}, \myboxnumb{0}, \myboxnumb{0}, \myboxnumb{0}, \myboxnumb{0}, \myboxnumb{0}, \myboxnumb{0}, \myboxnumb{0}, \myboxnumb{0}, \myboxnumb{0}]_S} \fi
 \ifnum#1=4 {[
 \myboxnumb{0}, \myboxnumb{0}, \myboxnumb{0}, \myboxnumb{1}, \myboxnumb{0}, \myboxnumb{0}, \myboxnumb{0}, \myboxnumb{0}, \myboxnumb{0}, \myboxnumb{0}, \myboxnumb{0}, \myboxnumb{0}, \myboxnumb{0}, \myboxnumb{0}, \myboxnumb{0}, \myboxnumb{0}, \myboxnumb{0}, \myboxnumb{0}, \myboxnumb{0}, \myboxnumb{0}, \myboxnumb{0}, \myboxnumb{0}]_S} \fi
 \ifnum#1=5 {[
 \myboxnumb{0}, \myboxnumb{0}, \myboxnumb{0}, \myboxnumb{0}, \myboxnumb{1}, \myboxnumb{0}, \myboxnumb{0}, \myboxnumb{0}, \myboxnumb{0}, \myboxnumb{0}, \myboxnumb{0}, \myboxnumb{0}, \myboxnumb{0}, \myboxnumb{0}, \myboxnumb{0}, \myboxnumb{0}, \myboxnumb{0}, \myboxnumb{0}, \myboxnumb{0}, \myboxnumb{0}, \myboxnumb{0}, \myboxnumb{0}]_S} \fi
 \ifnum#1=6 {[
 \myboxnumb{0}, \myboxnumb{0}, \myboxnumb{0}, \myboxnumb{0}, \myboxnumb{0}, \myboxnumb{1}, \myboxnumb{0}, \myboxnumb{0}, \myboxnumb{0}, \myboxnumb{0}, \myboxnumb{0}, \myboxnumb{0}, \myboxnumb{0}, \myboxnumb{0}, \myboxnumb{0}, \myboxnumb{0}, \myboxnumb{0}, \myboxnumb{0}, \myboxnumb{0}, \myboxnumb{0}, \myboxnumb{0}, \myboxnumb{0}]_S} \fi
 \ifnum#1=7 {[
 \myboxnumb{0}, \myboxnumb{0}, \myboxnumb{0}, \myboxnumb{0}, \myboxnumb{0}, \myboxnumb{0}, \myboxnumb{1}, \myboxnumb{0}, \myboxnumb{0}, \myboxnumb{0}, \myboxnumb{0}, \myboxnumb{0}, \myboxnumb{0}, \myboxnumb{0}, \myboxnumb{0}, \myboxnumb{0}, \myboxnumb{0}, \myboxnumb{0}, \myboxnumb{0}, \myboxnumb{0}, \myboxnumb{0}, \myboxnumb{0}]_S} \fi
 \ifnum#1=8 {[
 \myboxnumb{1}, \myboxnumb{1}, \myboxnumb{1}, \myboxnumb{1}, \myboxnumb{-1}, \myboxnumb{-1}, \myboxnumb{-1}, \myboxnumb{0}, \myboxnumb{0}, \myboxnumb{0}, \myboxnumb{0}, \myboxnumb{0}, \myboxnumb{0}, \myboxnumb{0}, \myboxnumb{0}, \myboxnumb{0}, \myboxnumb{0}, \myboxnumb{0}, \myboxnumb{0}, \myboxnumb{0}, \myboxnumb{0}, \myboxnumb{0}]_S} \fi
 \ifnum#1=9 {[
 \myboxnumb{0}, \myboxnumb{0}, \myboxnumb{0}, \myboxnumb{0}, \myboxnumb{0}, \myboxnumb{0}, \myboxnumb{0}, \myboxnumb{1}, \myboxnumb{0}, \myboxnumb{0}, \myboxnumb{0}, \myboxnumb{0}, \myboxnumb{0}, \myboxnumb{0}, \myboxnumb{0}, \myboxnumb{0}, \myboxnumb{0}, \myboxnumb{0}, \myboxnumb{0}, \myboxnumb{0}, \myboxnumb{0}, \myboxnumb{0}]_S} \fi
 \ifnum#1=10 {[
 \myboxnumb{0}, \myboxnumb{0}, \myboxnumb{0}, \myboxnumb{0}, \myboxnumb{0}, \myboxnumb{0}, \myboxnumb{0}, \myboxnumb{0}, \myboxnumb{1}, \myboxnumb{0}, \myboxnumb{0}, \myboxnumb{0}, \myboxnumb{0}, \myboxnumb{0}, \myboxnumb{0}, \myboxnumb{0}, \myboxnumb{0}, \myboxnumb{0}, \myboxnumb{0}, \myboxnumb{0}, \myboxnumb{0}, \myboxnumb{0}]_S} \fi
 \ifnum#1=11 {[
 \myboxnumb{0}, \myboxnumb{0}, \myboxnumb{0}, \myboxnumb{0}, \myboxnumb{0}, \myboxnumb{0}, \myboxnumb{0}, \myboxnumb{0}, \myboxnumb{0}, \myboxnumb{1}, \myboxnumb{0}, \myboxnumb{0}, \myboxnumb{0}, \myboxnumb{0}, \myboxnumb{0}, \myboxnumb{0}, \myboxnumb{0}, \myboxnumb{0}, \myboxnumb{0}, \myboxnumb{0}, \myboxnumb{0}, \myboxnumb{0}]_S} \fi
 \ifnum#1=12 {[
 \myboxnumb{1}, \myboxnumb{1}, \myboxnumb{1}, \myboxnumb{1}, \myboxnumb{0}, \myboxnumb{0}, \myboxnumb{0}, \myboxnumb{-1}, \myboxnumb{-1}, \myboxnumb{-1}, \myboxnumb{0}, \myboxnumb{0}, \myboxnumb{0}, \myboxnumb{0}, \myboxnumb{0}, \myboxnumb{0}, \myboxnumb{0}, \myboxnumb{0}, \myboxnumb{0}, \myboxnumb{0}, \myboxnumb{0}, \myboxnumb{0}]_S} \fi
 \ifnum#1=13 {[
 \myboxnumb{0}, \myboxnumb{1}, \myboxnumb{1}, \myboxnumb{1}, \myboxnumb{-1}, \myboxnumb{0}, \myboxnumb{0}, \myboxnumb{-1}, \myboxnumb{0}, \myboxnumb{0}, \myboxnumb{0}, \myboxnumb{0}, \myboxnumb{0}, \myboxnumb{0}, \myboxnumb{0}, \myboxnumb{0}, \myboxnumb{0}, \myboxnumb{0}, \myboxnumb{0}, \myboxnumb{0}, \myboxnumb{0}, \myboxnumb{0}]_S} \fi
 \ifnum#1=14 {[
 \myboxnumb{1}, \myboxnumb{0}, \myboxnumb{1}, \myboxnumb{1}, \myboxnumb{0}, \myboxnumb{-1}, \myboxnumb{0}, \myboxnumb{0}, \myboxnumb{-1}, \myboxnumb{0}, \myboxnumb{0}, \myboxnumb{0}, \myboxnumb{0}, \myboxnumb{0}, \myboxnumb{0}, \myboxnumb{0}, \myboxnumb{0}, \myboxnumb{0}, \myboxnumb{0}, \myboxnumb{0}, \myboxnumb{0}, \myboxnumb{0}]_S} \fi
 \ifnum#1=15 {[
 \myboxnumb{1}, \myboxnumb{1}, \myboxnumb{0}, \myboxnumb{1}, \myboxnumb{0}, \myboxnumb{0}, \myboxnumb{-1}, \myboxnumb{0}, \myboxnumb{0}, \myboxnumb{-1}, \myboxnumb{0}, \myboxnumb{0}, \myboxnumb{0}, \myboxnumb{0}, \myboxnumb{0}, \myboxnumb{0}, \myboxnumb{0}, \myboxnumb{0}, \myboxnumb{0}, \myboxnumb{0}, \myboxnumb{0}, \myboxnumb{0}]_S} \fi
 \ifnum#1=16 {[
 \myboxnumb{-1}, \myboxnumb{-1}, \myboxnumb{-1}, \myboxnumb{-2}, \myboxnumb{1}, \myboxnumb{1}, \myboxnumb{1}, \myboxnumb{1}, \myboxnumb{1}, \myboxnumb{1}, \myboxnumb{0}, \myboxnumb{0}, \myboxnumb{0}, \myboxnumb{0}, \myboxnumb{0}, \myboxnumb{0}, \myboxnumb{0}, \myboxnumb{0}, \myboxnumb{0}, \myboxnumb{0}, \myboxnumb{0}, \myboxnumb{0}]_S} \fi
 \ifnum#1=17 {[
 \myboxnumb{0}, \myboxnumb{0}, \myboxnumb{0}, \myboxnumb{0}, \myboxnumb{0}, \myboxnumb{0}, \myboxnumb{0}, \myboxnumb{0}, \myboxnumb{0}, \myboxnumb{0}, \myboxnumb{1}, \myboxnumb{0}, \myboxnumb{0}, \myboxnumb{0}, \myboxnumb{0}, \myboxnumb{0}, \myboxnumb{0}, \myboxnumb{0}, \myboxnumb{0}, \myboxnumb{0}, \myboxnumb{0}, \myboxnumb{0}]_S} \fi
 \ifnum#1=18 {[
 \myboxnumb{0}, \myboxnumb{0}, \myboxnumb{0}, \myboxnumb{0}, \myboxnumb{0}, \myboxnumb{0}, \myboxnumb{0}, \myboxnumb{0}, \myboxnumb{0}, \myboxnumb{0}, \myboxnumb{0}, \myboxnumb{1}, \myboxnumb{0}, \myboxnumb{0}, \myboxnumb{0}, \myboxnumb{0}, \myboxnumb{0}, \myboxnumb{0}, \myboxnumb{0}, \myboxnumb{0}, \myboxnumb{0}, \myboxnumb{0}]_S} \fi
 \ifnum#1=19 {[
 \myboxnumb{0}, \myboxnumb{0}, \myboxnumb{0}, \myboxnumb{0}, \myboxnumb{0}, \myboxnumb{0}, \myboxnumb{0}, \myboxnumb{0}, \myboxnumb{0}, \myboxnumb{0}, \myboxnumb{0}, \myboxnumb{0}, \myboxnumb{1}, \myboxnumb{0}, \myboxnumb{0}, \myboxnumb{0}, \myboxnumb{0}, \myboxnumb{0}, \myboxnumb{0}, \myboxnumb{0}, \myboxnumb{0}, \myboxnumb{0}]_S} \fi
 \ifnum#1=20 {[
 \myboxnumb{1}, \myboxnumb{1}, \myboxnumb{1}, \myboxnumb{1}, \myboxnumb{0}, \myboxnumb{0}, \myboxnumb{0}, \myboxnumb{0}, \myboxnumb{0}, \myboxnumb{0}, \myboxnumb{-1}, \myboxnumb{-1}, \myboxnumb{-1}, \myboxnumb{0}, \myboxnumb{0}, \myboxnumb{0}, \myboxnumb{0}, \myboxnumb{0}, \myboxnumb{0}, \myboxnumb{0}, \myboxnumb{0}, \myboxnumb{0}]_S} \fi
 \ifnum#1=21 {[
 \myboxnumb{0}, \myboxnumb{0}, \myboxnumb{0}, \myboxnumb{0}, \myboxnumb{0}, \myboxnumb{0}, \myboxnumb{0}, \myboxnumb{0}, \myboxnumb{0}, \myboxnumb{0}, \myboxnumb{0}, \myboxnumb{0}, \myboxnumb{0}, \myboxnumb{1}, \myboxnumb{0}, \myboxnumb{0}, \myboxnumb{0}, \myboxnumb{0}, \myboxnumb{0}, \myboxnumb{0}, \myboxnumb{0}, \myboxnumb{0}]_S} \fi
 \ifnum#1=22 {[
 \myboxnumb{0}, \myboxnumb{0}, \myboxnumb{0}, \myboxnumb{0}, \myboxnumb{0}, \myboxnumb{0}, \myboxnumb{0}, \myboxnumb{0}, \myboxnumb{0}, \myboxnumb{0}, \myboxnumb{0}, \myboxnumb{0}, \myboxnumb{0}, \myboxnumb{0}, \myboxnumb{1}, \myboxnumb{0}, \myboxnumb{0}, \myboxnumb{0}, \myboxnumb{0}, \myboxnumb{0}, \myboxnumb{0}, \myboxnumb{0}]_S} \fi
 \ifnum#1=23 {[
 \myboxnumb{0}, \myboxnumb{0}, \myboxnumb{0}, \myboxnumb{0}, \myboxnumb{0}, \myboxnumb{0}, \myboxnumb{0}, \myboxnumb{0}, \myboxnumb{0}, \myboxnumb{0}, \myboxnumb{0}, \myboxnumb{0}, \myboxnumb{0}, \myboxnumb{0}, \myboxnumb{0}, \myboxnumb{1}, \myboxnumb{0}, \myboxnumb{0}, \myboxnumb{0}, \myboxnumb{0}, \myboxnumb{0}, \myboxnumb{0}]_S} \fi
 \ifnum#1=24 {[
 \myboxnumb{1}, \myboxnumb{1}, \myboxnumb{1}, \myboxnumb{1}, \myboxnumb{0}, \myboxnumb{0}, \myboxnumb{0}, \myboxnumb{0}, \myboxnumb{0}, \myboxnumb{0}, \myboxnumb{0}, \myboxnumb{0}, \myboxnumb{0}, \myboxnumb{-1}, \myboxnumb{-1}, \myboxnumb{-1}, \myboxnumb{0}, \myboxnumb{0}, \myboxnumb{0}, \myboxnumb{0}, \myboxnumb{0}, \myboxnumb{0}]_S} \fi
 \ifnum#1=25 {[
 \myboxnumb{0}, \myboxnumb{0}, \myboxnumb{0}, \myboxnumb{0}, \myboxnumb{0}, \myboxnumb{0}, \myboxnumb{0}, \myboxnumb{0}, \myboxnumb{0}, \myboxnumb{0}, \myboxnumb{0}, \myboxnumb{0}, \myboxnumb{0}, \myboxnumb{0}, \myboxnumb{0}, \myboxnumb{0}, \myboxnumb{1}, \myboxnumb{0}, \myboxnumb{0}, \myboxnumb{0}, \myboxnumb{0}, \myboxnumb{0}]_S} \fi
 \ifnum#1=26 {[
 \myboxnumb{0}, \myboxnumb{0}, \myboxnumb{0}, \myboxnumb{0}, \myboxnumb{0}, \myboxnumb{0}, \myboxnumb{0}, \myboxnumb{0}, \myboxnumb{0}, \myboxnumb{0}, \myboxnumb{0}, \myboxnumb{0}, \myboxnumb{0}, \myboxnumb{0}, \myboxnumb{0}, \myboxnumb{0}, \myboxnumb{0}, \myboxnumb{1}, \myboxnumb{0}, \myboxnumb{0}, \myboxnumb{0}, \myboxnumb{0}]_S} \fi
 \ifnum#1=27 {[
 \myboxnumb{0}, \myboxnumb{0}, \myboxnumb{0}, \myboxnumb{0}, \myboxnumb{0}, \myboxnumb{0}, \myboxnumb{0}, \myboxnumb{0}, \myboxnumb{0}, \myboxnumb{0}, \myboxnumb{0}, \myboxnumb{0}, \myboxnumb{0}, \myboxnumb{0}, \myboxnumb{0}, \myboxnumb{0}, \myboxnumb{0}, \myboxnumb{0}, \myboxnumb{1}, \myboxnumb{0}, \myboxnumb{0}, \myboxnumb{0}]_S} \fi
 \ifnum#1=28 {[
 \myboxnumb{1}, \myboxnumb{1}, \myboxnumb{1}, \myboxnumb{1}, \myboxnumb{0}, \myboxnumb{0}, \myboxnumb{0}, \myboxnumb{0}, \myboxnumb{0}, \myboxnumb{0}, \myboxnumb{0}, \myboxnumb{0}, \myboxnumb{0}, \myboxnumb{0}, \myboxnumb{0}, \myboxnumb{0}, \myboxnumb{-1}, \myboxnumb{-1}, \myboxnumb{-1}, \myboxnumb{0}, \myboxnumb{0}, \myboxnumb{0}]_S} \fi
 \ifnum#1=29 {[
 \myboxnumb{1}, \myboxnumb{1}, \myboxnumb{1}, \myboxnumb{1}, \myboxnumb{0}, \myboxnumb{0}, \myboxnumb{0}, \myboxnumb{0}, \myboxnumb{0}, \myboxnumb{0}, \myboxnumb{-1}, \myboxnumb{0}, \myboxnumb{0}, \myboxnumb{-1}, \myboxnumb{0}, \myboxnumb{0}, \myboxnumb{-1}, \myboxnumb{0}, \myboxnumb{0}, \myboxnumb{0}, \myboxnumb{0}, \myboxnumb{0}]_S} \fi
 \ifnum#1=30 {[
 \myboxnumb{1}, \myboxnumb{1}, \myboxnumb{1}, \myboxnumb{1}, \myboxnumb{0}, \myboxnumb{0}, \myboxnumb{0}, \myboxnumb{0}, \myboxnumb{0}, \myboxnumb{0}, \myboxnumb{0}, \myboxnumb{-1}, \myboxnumb{0}, \myboxnumb{0}, \myboxnumb{-1}, \myboxnumb{0}, \myboxnumb{0}, \myboxnumb{-1}, \myboxnumb{0}, \myboxnumb{0}, \myboxnumb{0}, \myboxnumb{0}]_S} \fi
 \ifnum#1=31 {[
 \myboxnumb{1}, \myboxnumb{1}, \myboxnumb{1}, \myboxnumb{1}, \myboxnumb{0}, \myboxnumb{0}, \myboxnumb{0}, \myboxnumb{0}, \myboxnumb{0}, \myboxnumb{0}, \myboxnumb{0}, \myboxnumb{0}, \myboxnumb{-1}, \myboxnumb{0}, \myboxnumb{0}, \myboxnumb{-1}, \myboxnumb{0}, \myboxnumb{0}, \myboxnumb{-1}, \myboxnumb{0}, \myboxnumb{0}, \myboxnumb{0}]_S} \fi
 \ifnum#1=32 {[
 \myboxnumb{-2}, \myboxnumb{-2}, \myboxnumb{-2}, \myboxnumb{-2}, \myboxnumb{0}, \myboxnumb{0}, \myboxnumb{0}, \myboxnumb{0}, \myboxnumb{0}, \myboxnumb{0}, \myboxnumb{1}, \myboxnumb{1}, \myboxnumb{1}, \myboxnumb{1}, \myboxnumb{1}, \myboxnumb{1}, \myboxnumb{1}, \myboxnumb{1}, \myboxnumb{1}, \myboxnumb{0}, \myboxnumb{0}, \myboxnumb{0}]_S} \fi
 \ifnum#1=33 {[
 \myboxnumb{0}, \myboxnumb{0}, \myboxnumb{0}, \myboxnumb{0}, \myboxnumb{0}, \myboxnumb{0}, \myboxnumb{0}, \myboxnumb{0}, \myboxnumb{0}, \myboxnumb{0}, \myboxnumb{0}, \myboxnumb{0}, \myboxnumb{0}, \myboxnumb{0}, \myboxnumb{0}, \myboxnumb{0}, \myboxnumb{0}, \myboxnumb{0}, \myboxnumb{0}, \myboxnumb{1}, \myboxnumb{0}, \myboxnumb{0}]_S} \fi
 \ifnum#1=34 {[
 \myboxnumb{0}, \myboxnumb{1}, \myboxnumb{1}, \myboxnumb{0}, \myboxnumb{0}, \myboxnumb{1}, \myboxnumb{1}, \myboxnumb{0}, \myboxnumb{0}, \myboxnumb{0}, \myboxnumb{-1}, \myboxnumb{-1}, \myboxnumb{0}, \myboxnumb{0}, \myboxnumb{0}, \myboxnumb{0}, \myboxnumb{0}, \myboxnumb{0}, \myboxnumb{0}, \myboxnumb{-1}, \myboxnumb{0}, \myboxnumb{0}]_S} \fi
 \ifnum#1=35 {[
 \myboxnumb{0}, \myboxnumb{0}, \myboxnumb{0}, \myboxnumb{0}, \myboxnumb{0}, \myboxnumb{0}, \myboxnumb{0}, \myboxnumb{0}, \myboxnumb{0}, \myboxnumb{0}, \myboxnumb{0}, \myboxnumb{0}, \myboxnumb{0}, \myboxnumb{0}, \myboxnumb{0}, \myboxnumb{0}, \myboxnumb{0}, \myboxnumb{0}, \myboxnumb{0}, \myboxnumb{0}, \myboxnumb{1}, \myboxnumb{0}]_S} \fi
 \ifnum#1=36 {[
 \myboxnumb{1}, \myboxnumb{0}, \myboxnumb{0}, \myboxnumb{1}, \myboxnumb{0}, \myboxnumb{-1}, \myboxnumb{-1}, \myboxnumb{0}, \myboxnumb{0}, \myboxnumb{0}, \myboxnumb{1}, \myboxnumb{1}, \myboxnumb{0}, \myboxnumb{0}, \myboxnumb{0}, \myboxnumb{0}, \myboxnumb{0}, \myboxnumb{0}, \myboxnumb{0}, \myboxnumb{0}, \myboxnumb{-1}, \myboxnumb{0}]_S} \fi
 \ifnum#1=37 {[
 \myboxnumb{0}, \myboxnumb{0}, \myboxnumb{0}, \myboxnumb{0}, \myboxnumb{1}, \myboxnumb{1}, \myboxnumb{0}, \myboxnumb{1}, \myboxnumb{1}, \myboxnumb{0}, \myboxnumb{-1}, \myboxnumb{0}, \myboxnumb{0}, \myboxnumb{-1}, \myboxnumb{0}, \myboxnumb{0}, \myboxnumb{0}, \myboxnumb{0}, \myboxnumb{0}, \myboxnumb{-1}, \myboxnumb{0}, \myboxnumb{0}]_S} \fi
 \ifnum#1=38 {[
 \myboxnumb{2}, \myboxnumb{1}, \myboxnumb{1}, \myboxnumb{2}, \myboxnumb{-1}, \myboxnumb{-2}, \myboxnumb{-1}, \myboxnumb{-1}, \myboxnumb{-1}, \myboxnumb{0}, \myboxnumb{1}, \myboxnumb{0}, \myboxnumb{0}, \myboxnumb{0}, \myboxnumb{-1}, \myboxnumb{0}, \myboxnumb{0}, \myboxnumb{0}, \myboxnumb{0}, \myboxnumb{1}, \myboxnumb{0}, \myboxnumb{0}]_S} \fi
 \ifnum#1=39 {[
 \myboxnumb{2}, \myboxnumb{2}, \myboxnumb{2}, \myboxnumb{2}, \myboxnumb{-1}, \myboxnumb{0}, \myboxnumb{-1}, \myboxnumb{-1}, \myboxnumb{0}, \myboxnumb{-1}, \myboxnumb{0}, \myboxnumb{0}, \myboxnumb{-1}, \myboxnumb{0}, \myboxnumb{0}, \myboxnumb{-1}, \myboxnumb{0}, \myboxnumb{0}, \myboxnumb{0}, \myboxnumb{0}, \myboxnumb{-1}, \myboxnumb{0}]_S} \fi
 \ifnum#1=40 {[
 \myboxnumb{-3}, \myboxnumb{-2}, \myboxnumb{-2}, \myboxnumb{-3}, \myboxnumb{1}, \myboxnumb{1}, \myboxnumb{2}, \myboxnumb{1}, \myboxnumb{0}, \myboxnumb{1}, \myboxnumb{0}, \myboxnumb{0}, \myboxnumb{1}, \myboxnumb{1}, \myboxnumb{1}, \myboxnumb{1}, \myboxnumb{0}, \myboxnumb{0}, \myboxnumb{0}, \myboxnumb{0}, \myboxnumb{1}, \myboxnumb{0}]_S} \fi
 \ifnum#1=41 {[
 \myboxnumb{0}, \myboxnumb{0}, \myboxnumb{0}, \myboxnumb{0}, \myboxnumb{0}, \myboxnumb{0}, \myboxnumb{1}, \myboxnumb{0}, \myboxnumb{-1}, \myboxnumb{0}, \myboxnumb{0}, \myboxnumb{-1}, \myboxnumb{0}, \myboxnumb{1}, \myboxnumb{0}, \myboxnumb{1}, \myboxnumb{0}, \myboxnumb{-1}, \myboxnumb{0}, \myboxnumb{0}, \myboxnumb{1}, \myboxnumb{0}]_S} \fi
 \ifnum#1=42 {[
 \myboxnumb{2}, \myboxnumb{1}, \myboxnumb{1}, \myboxnumb{2}, \myboxnumb{0}, \myboxnumb{0}, \myboxnumb{-1}, \myboxnumb{0}, \myboxnumb{0}, \myboxnumb{-1}, \myboxnumb{0}, \myboxnumb{1}, \myboxnumb{0}, \myboxnumb{-1}, \myboxnumb{0}, \myboxnumb{-1}, \myboxnumb{-1}, \myboxnumb{0}, \myboxnumb{0}, \myboxnumb{0}, \myboxnumb{-1}, \myboxnumb{0}]_S} \fi
 \ifnum#1=43 {[
 \myboxnumb{-3}, \myboxnumb{-2}, \myboxnumb{-2}, \myboxnumb{-3}, \myboxnumb{1}, \myboxnumb{1}, \myboxnumb{1}, \myboxnumb{1}, \myboxnumb{1}, \myboxnumb{1}, \myboxnumb{0}, \myboxnumb{0}, \myboxnumb{1}, \myboxnumb{0}, \myboxnumb{1}, \myboxnumb{1}, \myboxnumb{1}, \myboxnumb{1}, \myboxnumb{1}, \myboxnumb{-1}, \myboxnumb{0}, \myboxnumb{0}]_S} \fi
 \ifnum#1=44 {[
 \myboxnumb{2}, \myboxnumb{2}, \myboxnumb{2}, \myboxnumb{2}, \myboxnumb{-1}, \myboxnumb{-1}, \myboxnumb{-1}, \myboxnumb{-1}, \myboxnumb{0}, \myboxnumb{0}, \myboxnumb{0}, \myboxnumb{0}, \myboxnumb{-1}, \myboxnumb{0}, \myboxnumb{-1}, \myboxnumb{-1}, \myboxnumb{0}, \myboxnumb{0}, \myboxnumb{-1}, \myboxnumb{1}, \myboxnumb{0}, \myboxnumb{0}]_S} \fi
 \ifnum#1=45 {[
 \myboxnumb{1}, \myboxnumb{1}, \myboxnumb{1}, \myboxnumb{1}, \myboxnumb{-1}, \myboxnumb{-1}, \myboxnumb{-1}, \myboxnumb{-1}, \myboxnumb{0}, \myboxnumb{0}, \myboxnumb{1}, \myboxnumb{1}, \myboxnumb{0}, \myboxnumb{0}, \myboxnumb{0}, \myboxnumb{-1}, \myboxnumb{0}, \myboxnumb{1}, \myboxnumb{0}, \myboxnumb{0}, \myboxnumb{-1}, \myboxnumb{0}]_S} \fi
 \ifnum#1=46 {[
 \myboxnumb{-3}, \myboxnumb{-2}, \myboxnumb{-2}, \myboxnumb{-3}, \myboxnumb{1}, \myboxnumb{1}, \myboxnumb{1}, \myboxnumb{1}, \myboxnumb{1}, \myboxnumb{1}, \myboxnumb{0}, \myboxnumb{0}, \myboxnumb{0}, \myboxnumb{1}, \myboxnumb{1}, \myboxnumb{1}, \myboxnumb{1}, \myboxnumb{0}, \myboxnumb{0}, \myboxnumb{0}, \myboxnumb{1}, \myboxnumb{0}]_S} \fi
 \ifnum#1=47 {[
 \myboxnumb{2}, \myboxnumb{1}, \myboxnumb{1}, \myboxnumb{2}, \myboxnumb{0}, \myboxnumb{-1}, \myboxnumb{0}, \myboxnumb{0}, \myboxnumb{-1}, \myboxnumb{0}, \myboxnumb{0}, \myboxnumb{0}, \myboxnumb{0}, \myboxnumb{0}, \myboxnumb{-1}, \myboxnumb{0}, \myboxnumb{-1}, \myboxnumb{-1}, \myboxnumb{-1}, \myboxnumb{1}, \myboxnumb{0}, \myboxnumb{0}]_S} \fi
 \ifnum#1=48 {[
 \myboxnumb{1}, \myboxnumb{1}, \myboxnumb{1}, \myboxnumb{1}, \myboxnumb{0}, \myboxnumb{1}, \myboxnumb{0}, \myboxnumb{0}, \myboxnumb{0}, \myboxnumb{-1}, \myboxnumb{-1}, \myboxnumb{-1}, \myboxnumb{0}, \myboxnumb{-1}, \myboxnumb{0}, \myboxnumb{0}, \myboxnumb{0}, \myboxnumb{0}, \myboxnumb{1}, \myboxnumb{-1}, \myboxnumb{0}, \myboxnumb{0}]_S} \fi
 \ifnum#1=49 {[
 \myboxnumb{0}, \myboxnumb{0}, \myboxnumb{0}, \myboxnumb{0}, \myboxnumb{0}, \myboxnumb{0}, \myboxnumb{0}, \myboxnumb{0}, \myboxnumb{0}, \myboxnumb{0}, \myboxnumb{0}, \myboxnumb{0}, \myboxnumb{0}, \myboxnumb{0}, \myboxnumb{0}, \myboxnumb{0}, \myboxnumb{0}, \myboxnumb{0}, \myboxnumb{0}, \myboxnumb{0}, \myboxnumb{0}, \myboxnumb{1}]_S} \fi
 \ifnum#1=50 {[
 \myboxnumb{1}, \myboxnumb{0}, \myboxnumb{1}, \myboxnumb{0}, \myboxnumb{-1}, \myboxnumb{-1}, \myboxnumb{0}, \myboxnumb{0}, \myboxnumb{-1}, \myboxnumb{1}, \myboxnumb{0}, \myboxnumb{-1}, \myboxnumb{-1}, \myboxnumb{0}, \myboxnumb{0}, \myboxnumb{0}, \myboxnumb{0}, \myboxnumb{0}, \myboxnumb{0}, \myboxnumb{1}, \myboxnumb{1}, \myboxnumb{1}]_S} \fi
 \ifnum#1=51 {[
 \myboxnumb{1}, \myboxnumb{1}, \myboxnumb{1}, \myboxnumb{2}, \myboxnumb{0}, \myboxnumb{0}, \myboxnumb{-1}, \myboxnumb{-1}, \myboxnumb{0}, \myboxnumb{-1}, \myboxnumb{0}, \myboxnumb{1}, \myboxnumb{0}, \myboxnumb{0}, \myboxnumb{0}, \myboxnumb{0}, \myboxnumb{0}, \myboxnumb{0}, \myboxnumb{0}, \myboxnumb{0}, \myboxnumb{-1}, \myboxnumb{-1}]_S} \fi
 \ifnum#1=52 {[
 \myboxnumb{-1}, \myboxnumb{0}, \myboxnumb{-1}, \myboxnumb{-1}, \myboxnumb{1}, \myboxnumb{1}, \myboxnumb{1}, \myboxnumb{1}, \myboxnumb{1}, \myboxnumb{0}, \myboxnumb{0}, \myboxnumb{0}, \myboxnumb{1}, \myboxnumb{0}, \myboxnumb{0}, \myboxnumb{0}, \myboxnumb{0}, \myboxnumb{0}, \myboxnumb{0}, \myboxnumb{-1}, \myboxnumb{0}, \myboxnumb{-1}]_S} \fi
 \ifnum#1=53 {[
 \myboxnumb{1}, \myboxnumb{-1}, \myboxnumb{0}, \myboxnumb{0}, \myboxnumb{0}, \myboxnumb{-1}, \myboxnumb{0}, \myboxnumb{1}, \myboxnumb{-1}, \myboxnumb{1}, \myboxnumb{0}, \myboxnumb{-1}, \myboxnumb{0}, \myboxnumb{0}, \myboxnumb{0}, \myboxnumb{1}, \myboxnumb{-1}, \myboxnumb{-1}, \myboxnumb{0}, \myboxnumb{1}, \myboxnumb{1}, \myboxnumb{1}]_S} \fi
 \ifnum#1=54 {[
 \myboxnumb{4}, \myboxnumb{3}, \myboxnumb{3}, \myboxnumb{4}, \myboxnumb{-1}, \myboxnumb{-1}, \myboxnumb{-1}, \myboxnumb{-1}, \myboxnumb{-1}, \myboxnumb{-1}, \myboxnumb{-1}, \myboxnumb{-1}, \myboxnumb{-1}, \myboxnumb{-1}, \myboxnumb{-1}, \myboxnumb{-1}, \myboxnumb{-1}, \myboxnumb{-1}, \myboxnumb{0}, \myboxnumb{0}, \myboxnumb{0}, \myboxnumb{1}]_S} \fi
 \ifnum#1=55 {[
 \myboxnumb{-3}, \myboxnumb{-1}, \myboxnumb{-2}, \myboxnumb{-3}, \myboxnumb{1}, \myboxnumb{2}, \myboxnumb{1}, \myboxnumb{1}, \myboxnumb{2}, \myboxnumb{0}, \myboxnumb{0}, \myboxnumb{1}, \myboxnumb{1}, \myboxnumb{0}, \myboxnumb{1}, \myboxnumb{0}, \myboxnumb{1}, \myboxnumb{1}, \myboxnumb{0}, \myboxnumb{-1}, \myboxnumb{0}, \myboxnumb{-1}]_S} \fi
 \ifnum#1=56 {[
 \myboxnumb{-1}, \myboxnumb{0}, \myboxnumb{0}, \myboxnumb{0}, \myboxnumb{0}, \myboxnumb{0}, \myboxnumb{0}, \myboxnumb{-1}, \myboxnumb{0}, \myboxnumb{0}, \myboxnumb{1}, \myboxnumb{1}, \myboxnumb{0}, \myboxnumb{1}, \myboxnumb{0}, \myboxnumb{0}, \myboxnumb{1}, \myboxnumb{1}, \myboxnumb{0}, \myboxnumb{0}, \myboxnumb{-1}, \myboxnumb{-1}]_S} \fi
 \ifnum#1=57 {[
 \myboxnumb{1}, \myboxnumb{2}, \myboxnumb{1}, \myboxnumb{2}, \myboxnumb{-1}, \myboxnumb{0}, \myboxnumb{-1}, \myboxnumb{-1}, \myboxnumb{0}, \myboxnumb{-1}, \myboxnumb{0}, \myboxnumb{1}, \myboxnumb{0}, \myboxnumb{0}, \myboxnumb{0}, \myboxnumb{-1}, \myboxnumb{0}, \myboxnumb{1}, \myboxnumb{0}, \myboxnumb{0}, \myboxnumb{-1}, \myboxnumb{-1}]_S} \fi
 \ifnum#1=58 {[
 \myboxnumb{-3}, \myboxnumb{-2}, \myboxnumb{-2}, \myboxnumb{-3}, \myboxnumb{2}, \myboxnumb{2}, \myboxnumb{2}, \myboxnumb{1}, \myboxnumb{1}, \myboxnumb{0}, \myboxnumb{0}, \myboxnumb{0}, \myboxnumb{1}, \myboxnumb{1}, \myboxnumb{1}, \myboxnumb{1}, \myboxnumb{1}, \myboxnumb{0}, \myboxnumb{0}, \myboxnumb{-1}, \myboxnumb{0}, \myboxnumb{-1}]_S} \fi
 \ifnum#1=59 {[
 \myboxnumb{2}, \myboxnumb{1}, \myboxnumb{2}, \myboxnumb{2}, \myboxnumb{-1}, \myboxnumb{-2}, \myboxnumb{0}, \myboxnumb{0}, \myboxnumb{-1}, \myboxnumb{1}, \myboxnumb{0}, \myboxnumb{-1}, \myboxnumb{-1}, \myboxnumb{0}, \myboxnumb{-1}, \myboxnumb{0}, \myboxnumb{-1}, \myboxnumb{-1}, \myboxnumb{-1}, \myboxnumb{1}, \myboxnumb{1}, \myboxnumb{1}]_S} \fi
 \ifnum#1=60 {[
 \myboxnumb{1}, \myboxnumb{0}, \myboxnumb{0}, \myboxnumb{0}, \myboxnumb{0}, \myboxnumb{0}, \myboxnumb{-1}, \myboxnumb{0}, \myboxnumb{0}, \myboxnumb{0}, \myboxnumb{0}, \myboxnumb{0}, \myboxnumb{0}, \myboxnumb{-1}, \myboxnumb{0}, \myboxnumb{0}, \myboxnumb{0}, \myboxnumb{0}, \myboxnumb{1}, \myboxnumb{0}, \myboxnumb{0}, \myboxnumb{1}]_S} \fi
 \ifnum#1=61 {[
 \myboxnumb{-1}, \myboxnumb{0}, \myboxnumb{0}, \myboxnumb{-1}, \myboxnumb{1}, \myboxnumb{1}, \myboxnumb{1}, \myboxnumb{0}, \myboxnumb{1}, \myboxnumb{0}, \myboxnumb{0}, \myboxnumb{0}, \myboxnumb{0}, \myboxnumb{0}, \myboxnumb{0}, \myboxnumb{0}, \myboxnumb{1}, \myboxnumb{0}, \myboxnumb{0}, \myboxnumb{-1}, \myboxnumb{0}, \myboxnumb{-1}]_S} \fi
 \ifnum#1=62 {[
 \myboxnumb{-1}, \myboxnumb{0}, \myboxnumb{-1}, \myboxnumb{0}, \myboxnumb{0}, \myboxnumb{0}, \myboxnumb{-1}, \myboxnumb{0}, \myboxnumb{1}, \myboxnumb{0}, \myboxnumb{1}, \myboxnumb{2}, \myboxnumb{1}, \myboxnumb{0}, \myboxnumb{0}, \myboxnumb{0}, \myboxnumb{0}, \myboxnumb{1}, \myboxnumb{0}, \myboxnumb{0}, \myboxnumb{-1}, \myboxnumb{-1}]_S} \fi
 \ifnum#1=63 {[
 \myboxnumb{1}, \myboxnumb{0}, \myboxnumb{0}, \myboxnumb{0}, \myboxnumb{0}, \myboxnumb{0}, \myboxnumb{0}, \myboxnumb{0}, \myboxnumb{-1}, \myboxnumb{0}, \myboxnumb{0}, \myboxnumb{-1}, \myboxnumb{0}, \myboxnumb{0}, \myboxnumb{0}, \myboxnumb{0}, \myboxnumb{0}, \myboxnumb{0}, \myboxnumb{1}, \myboxnumb{0}, \myboxnumb{0}, \myboxnumb{1}]_S} \fi
 \ifnum#1=64 {[
 \myboxnumb{2}, \myboxnumb{1}, \myboxnumb{2}, \myboxnumb{2}, \myboxnumb{-1}, \myboxnumb{-1}, \myboxnumb{0}, \myboxnumb{0}, \myboxnumb{-1}, \myboxnumb{0}, \myboxnumb{-1}, \myboxnumb{-1}, \myboxnumb{-1}, \myboxnumb{0}, \myboxnumb{0}, \myboxnumb{0}, \myboxnumb{-1}, \myboxnumb{-1}, \myboxnumb{-1}, \myboxnumb{1}, \myboxnumb{1}, \myboxnumb{1}]_S} \fi
 \ifnum#1=65 {[
 \myboxnumb{0}, \myboxnumb{1}, \myboxnumb{1}, \myboxnumb{1}, \myboxnumb{0}, \myboxnumb{0}, \myboxnumb{0}, \myboxnumb{0}, \myboxnumb{0}, \myboxnumb{0}, \myboxnumb{0}, \myboxnumb{0}, \myboxnumb{0}, \myboxnumb{0}, \myboxnumb{0}, \myboxnumb{0}, \myboxnumb{0}, \myboxnumb{0}, \myboxnumb{-1}, \myboxnumb{0}, \myboxnumb{0}, \myboxnumb{-1}]_S} \fi
 \ifnum#1=66 {[
 \myboxnumb{0}, \myboxnumb{-1}, \myboxnumb{0}, \myboxnumb{-1}, \myboxnumb{0}, \myboxnumb{0}, \myboxnumb{1}, \myboxnumb{1}, \myboxnumb{0}, \myboxnumb{1}, \myboxnumb{0}, \myboxnumb{-1}, \myboxnumb{0}, \myboxnumb{0}, \myboxnumb{0}, \myboxnumb{0}, \myboxnumb{0}, \myboxnumb{-1}, \myboxnumb{0}, \myboxnumb{0}, \myboxnumb{1}, \myboxnumb{1}]_S} \fi
 \ifnum#1=67 {[
 \myboxnumb{2}, \myboxnumb{1}, \myboxnumb{1}, \myboxnumb{2}, \myboxnumb{-1}, \myboxnumb{-1}, \myboxnumb{-1}, \myboxnumb{-1}, \myboxnumb{-1}, \myboxnumb{0}, \myboxnumb{0}, \myboxnumb{0}, \myboxnumb{-1}, \myboxnumb{0}, \myboxnumb{0}, \myboxnumb{0}, \myboxnumb{-1}, \myboxnumb{0}, \myboxnumb{0}, \myboxnumb{1}, \myboxnumb{0}, \myboxnumb{1}]_S} \fi
 \ifnum#1=68 {[
 \myboxnumb{-1}, \myboxnumb{0}, \myboxnumb{-1}, \myboxnumb{-1}, \myboxnumb{1}, \myboxnumb{1}, \myboxnumb{0}, \myboxnumb{0}, \myboxnumb{1}, \myboxnumb{-1}, \myboxnumb{0}, \myboxnumb{1}, \myboxnumb{1}, \myboxnumb{0}, \myboxnumb{0}, \myboxnumb{0}, \myboxnumb{1}, \myboxnumb{1}, \myboxnumb{1}, \myboxnumb{-1}, \myboxnumb{-1}, \myboxnumb{-1}]_S} \fi
 \ifnum#1=69 {[
 \myboxnumb{4}, \myboxnumb{2}, \myboxnumb{3}, \myboxnumb{4}, \myboxnumb{-2}, \myboxnumb{-2}, \myboxnumb{-2}, \myboxnumb{-1}, \myboxnumb{-1}, \myboxnumb{0}, \myboxnumb{0}, \myboxnumb{0}, \myboxnumb{-1}, \myboxnumb{-1}, \myboxnumb{-1}, \myboxnumb{-1}, \myboxnumb{-1}, \myboxnumb{0}, \myboxnumb{0}, \myboxnumb{1}, \myboxnumb{-1}, \myboxnumb{1}]_S} \fi
 \ifnum#1=70 {[
 \myboxnumb{0}, \myboxnumb{0}, \myboxnumb{0}, \myboxnumb{0}, \myboxnumb{0}, \myboxnumb{-1}, \myboxnumb{0}, \myboxnumb{0}, \myboxnumb{0}, \myboxnumb{0}, \myboxnumb{1}, \myboxnumb{1}, \myboxnumb{0}, \myboxnumb{1}, \myboxnumb{0}, \myboxnumb{0}, \myboxnumb{0}, \myboxnumb{0}, \myboxnumb{-1}, \myboxnumb{1}, \myboxnumb{0}, \myboxnumb{-1}]_S} \fi
 \ifnum#1=71 {[
 \myboxnumb{-2}, \myboxnumb{0}, \myboxnumb{-1}, \myboxnumb{-1}, \myboxnumb{1}, \myboxnumb{2}, \myboxnumb{0}, \myboxnumb{0}, \myboxnumb{1}, \myboxnumb{-1}, \myboxnumb{0}, \myboxnumb{1}, \myboxnumb{1}, \myboxnumb{0}, \myboxnumb{1}, \myboxnumb{0}, \myboxnumb{1}, \myboxnumb{1}, \myboxnumb{1}, \myboxnumb{-2}, \myboxnumb{-1}, \myboxnumb{-1}]_S} \fi
 \ifnum#1=72 {[
 \myboxnumb{-1}, \myboxnumb{-1}, \myboxnumb{-1}, \myboxnumb{-2}, \myboxnumb{1}, \myboxnumb{1}, \myboxnumb{2}, \myboxnumb{1}, \myboxnumb{0}, \myboxnumb{1}, \myboxnumb{-1}, \myboxnumb{-2}, \myboxnumb{0}, \myboxnumb{0}, \myboxnumb{0}, \myboxnumb{1}, \myboxnumb{0}, \myboxnumb{-1}, \myboxnumb{0}, \myboxnumb{0}, \myboxnumb{2}, \myboxnumb{1}]_S} \fi
 \ifnum#1=73 {[
 \myboxnumb{-1}, \myboxnumb{1}, \myboxnumb{0}, \myboxnumb{0}, \myboxnumb{1}, \myboxnumb{2}, \myboxnumb{1}, \myboxnumb{0}, \myboxnumb{1}, \myboxnumb{-1}, \myboxnumb{-1}, \myboxnumb{0}, \myboxnumb{0}, \myboxnumb{0}, \myboxnumb{0}, \myboxnumb{0}, \myboxnumb{0}, \myboxnumb{0}, \myboxnumb{0}, \myboxnumb{-1}, \myboxnumb{0}, \myboxnumb{-1}]_S} \fi
 \ifnum#1=74 {[
 \myboxnumb{0}, \myboxnumb{0}, \myboxnumb{0}, \myboxnumb{0}, \myboxnumb{0}, \myboxnumb{0}, \myboxnumb{0}, \myboxnumb{0}, \myboxnumb{0}, \myboxnumb{0}, \myboxnumb{1}, \myboxnumb{1}, \myboxnumb{1}, \myboxnumb{0}, \myboxnumb{0}, \myboxnumb{0}, \myboxnumb{0}, \myboxnumb{0}, \myboxnumb{0}, \myboxnumb{0}, \myboxnumb{-1}, \myboxnumb{-1}]_S} \fi
 \ifnum#1=75 {[
 \myboxnumb{2}, \myboxnumb{1}, \myboxnumb{2}, \myboxnumb{2}, \myboxnumb{-1}, \myboxnumb{-1}, \myboxnumb{-1}, \myboxnumb{-1}, \myboxnumb{-1}, \myboxnumb{0}, \myboxnumb{0}, \myboxnumb{-1}, \myboxnumb{-1}, \myboxnumb{0}, \myboxnumb{0}, \myboxnumb{0}, \myboxnumb{0}, \myboxnumb{0}, \myboxnumb{0}, \myboxnumb{0}, \myboxnumb{0}, \myboxnumb{1}]_S} \fi
 \ifnum#1=76 {[
 \myboxnumb{0}, \myboxnumb{-1}, \myboxnumb{-1}, \myboxnumb{-1}, \myboxnumb{0}, \myboxnumb{-1}, \myboxnumb{0}, \myboxnumb{1}, \myboxnumb{0}, \myboxnumb{1}, \myboxnumb{0}, \myboxnumb{0}, \myboxnumb{0}, \myboxnumb{0}, \myboxnumb{0}, \myboxnumb{0}, \myboxnumb{0}, \myboxnumb{0}, \myboxnumb{0}, \myboxnumb{1}, \myboxnumb{1}, \myboxnumb{1}]_S} \fi
 \ifnum#1=77 {[
 \myboxnumb{0}, \myboxnumb{0}, \myboxnumb{0}, \myboxnumb{0}, \myboxnumb{0}, \myboxnumb{0}, \myboxnumb{-1}, \myboxnumb{0}, \myboxnumb{1}, \myboxnumb{0}, \myboxnumb{1}, \myboxnumb{1}, \myboxnumb{0}, \myboxnumb{0}, \myboxnumb{0}, \myboxnumb{-1}, \myboxnumb{1}, \myboxnumb{1}, \myboxnumb{0}, \myboxnumb{0}, \myboxnumb{-1}, \myboxnumb{-1}]_S} \fi
 \ifnum#1=78 {[
 \myboxnumb{-3}, \myboxnumb{-1}, \myboxnumb{-2}, \myboxnumb{-2}, \myboxnumb{1}, \myboxnumb{1}, \myboxnumb{1}, \myboxnumb{0}, \myboxnumb{1}, \myboxnumb{0}, \myboxnumb{0}, \myboxnumb{1}, \myboxnumb{1}, \myboxnumb{1}, \myboxnumb{1}, \myboxnumb{1}, \myboxnumb{1}, \myboxnumb{1}, \myboxnumb{0}, \myboxnumb{-1}, \myboxnumb{0}, \myboxnumb{-1}]_S} \fi
 \ifnum#1=79 {[
 \myboxnumb{2}, \myboxnumb{1}, \myboxnumb{1}, \myboxnumb{1}, \myboxnumb{-1}, \myboxnumb{-1}, \myboxnumb{0}, \myboxnumb{0}, \myboxnumb{-1}, \myboxnumb{0}, \myboxnumb{0}, \myboxnumb{-1}, \myboxnumb{0}, \myboxnumb{0}, \myboxnumb{-1}, \myboxnumb{0}, \myboxnumb{-1}, \myboxnumb{-1}, \myboxnumb{0}, \myboxnumb{1}, \myboxnumb{1}, \myboxnumb{1}]_S} \fi
 \ifnum#1=80 {[
 \myboxnumb{2}, \myboxnumb{1}, \myboxnumb{2}, \myboxnumb{2}, \myboxnumb{0}, \myboxnumb{0}, \myboxnumb{0}, \myboxnumb{0}, \myboxnumb{-1}, \myboxnumb{0}, \myboxnumb{-1}, \myboxnumb{-1}, \myboxnumb{-1}, \myboxnumb{-1}, \myboxnumb{0}, \myboxnumb{0}, \myboxnumb{-1}, \myboxnumb{-1}, \myboxnumb{0}, \myboxnumb{0}, \myboxnumb{0}, \myboxnumb{1}]_S} \fi
 \ifnum#1=81 {[
 \myboxnumb{2}, \myboxnumb{1}, \myboxnumb{1}, \myboxnumb{2}, \myboxnumb{-1}, \myboxnumb{-1}, \myboxnumb{-1}, \myboxnumb{0}, \myboxnumb{-1}, \myboxnumb{0}, \myboxnumb{0}, \myboxnumb{0}, \myboxnumb{0}, \myboxnumb{-1}, \myboxnumb{0}, \myboxnumb{0}, \myboxnumb{-1}, \myboxnumb{0}, \myboxnumb{0}, \myboxnumb{0}, \myboxnumb{0}, \myboxnumb{1}]_S} \fi
 \ifnum#1=82 {[
 \myboxnumb{2}, \myboxnumb{0}, \myboxnumb{1}, \myboxnumb{1}, \myboxnumb{0}, \myboxnumb{-1}, \myboxnumb{0}, \myboxnumb{0}, \myboxnumb{-1}, \myboxnumb{0}, \myboxnumb{0}, \myboxnumb{-1}, \myboxnumb{-1}, \myboxnumb{0}, \myboxnumb{-1}, \myboxnumb{0}, \myboxnumb{0}, \myboxnumb{-1}, \myboxnumb{0}, \myboxnumb{1}, \myboxnumb{1}, \myboxnumb{1}]_S} \fi
 \ifnum#1=83 {[
 \myboxnumb{0}, \myboxnumb{1}, \myboxnumb{1}, \myboxnumb{1}, \myboxnumb{0}, \myboxnumb{0}, \myboxnumb{0}, \myboxnumb{0}, \myboxnumb{1}, \myboxnumb{0}, \myboxnumb{0}, \myboxnumb{1}, \myboxnumb{0}, \myboxnumb{0}, \myboxnumb{0}, \myboxnumb{-1}, \myboxnumb{0}, \myboxnumb{0}, \myboxnumb{-1}, \myboxnumb{0}, \myboxnumb{-1}, \myboxnumb{-1}]_S} \fi
 \ifnum#1=84 {[
 \myboxnumb{-3}, \myboxnumb{-1}, \myboxnumb{-2}, \myboxnumb{-3}, \myboxnumb{1}, \myboxnumb{2}, \myboxnumb{1}, \myboxnumb{0}, \myboxnumb{1}, \myboxnumb{0}, \myboxnumb{0}, \myboxnumb{0}, \myboxnumb{1}, \myboxnumb{1}, \myboxnumb{1}, \myboxnumb{1}, \myboxnumb{1}, \myboxnumb{1}, \myboxnumb{1}, \myboxnumb{-1}, \myboxnumb{0}, \myboxnumb{-1}]_S} \fi
 \ifnum#1=85 {[
 \myboxnumb{0}, \myboxnumb{-1}, \myboxnumb{0}, \myboxnumb{-1}, \myboxnumb{0}, \myboxnumb{-1}, \myboxnumb{1}, \myboxnumb{0}, \myboxnumb{-1}, \myboxnumb{1}, \myboxnumb{0}, \myboxnumb{-1}, \myboxnumb{0}, \myboxnumb{1}, \myboxnumb{0}, \myboxnumb{1}, \myboxnumb{0}, \myboxnumb{-1}, \myboxnumb{0}, \myboxnumb{1}, \myboxnumb{1}, \myboxnumb{1}]_S} \fi
 \ifnum#1=86 {[
 \myboxnumb{2}, \myboxnumb{2}, \myboxnumb{2}, \myboxnumb{2}, \myboxnumb{-1}, \myboxnumb{0}, \myboxnumb{-1}, \myboxnumb{0}, \myboxnumb{0}, \myboxnumb{0}, \myboxnumb{-1}, \myboxnumb{-1}, \myboxnumb{-1}, \myboxnumb{-1}, \myboxnumb{0}, \myboxnumb{-1}, \myboxnumb{-1}, \myboxnumb{0}, \myboxnumb{0}, \myboxnumb{0}, \myboxnumb{0}, \myboxnumb{1}]_S} \fi
 \ifnum#1=87 {[
 \myboxnumb{-3}, \myboxnumb{-2}, \myboxnumb{-3}, \myboxnumb{-3}, \myboxnumb{2}, \myboxnumb{2}, \myboxnumb{1}, \myboxnumb{1}, \myboxnumb{1}, \myboxnumb{0}, \myboxnumb{0}, \myboxnumb{1}, \myboxnumb{1}, \myboxnumb{0}, \myboxnumb{1}, \myboxnumb{1}, \myboxnumb{1}, \myboxnumb{1}, \myboxnumb{1}, \myboxnumb{-1}, \myboxnumb{0}, \myboxnumb{-1}]_S} \fi
 \ifnum#1=88 {[
 \myboxnumb{2}, \myboxnumb{2}, \myboxnumb{2}, \myboxnumb{3}, \myboxnumb{-1}, \myboxnumb{-1}, \myboxnumb{-1}, \myboxnumb{-1}, \myboxnumb{0}, \myboxnumb{-1}, \myboxnumb{1}, \myboxnumb{1}, \myboxnumb{0}, \myboxnumb{0}, \myboxnumb{-1}, \myboxnumb{-1}, \myboxnumb{0}, \myboxnumb{0}, \myboxnumb{-1}, \myboxnumb{0}, \myboxnumb{-1}, \myboxnumb{-1}]_S} \fi
 \ifnum#1=89 {[
 \myboxnumb{-2}, \myboxnumb{-1}, \myboxnumb{-2}, \myboxnumb{-2}, \myboxnumb{1}, \myboxnumb{1}, \myboxnumb{1}, \myboxnumb{0}, \myboxnumb{1}, \myboxnumb{0}, \myboxnumb{1}, \myboxnumb{1}, \myboxnumb{1}, \myboxnumb{1}, \myboxnumb{0}, \myboxnumb{0}, \myboxnumb{1}, \myboxnumb{0}, \myboxnumb{0}, \myboxnumb{0}, \myboxnumb{0}, \myboxnumb{-1}]_S} \fi
 \ifnum#1=90 {[
 \myboxnumb{0}, \myboxnumb{0}, \myboxnumb{0}, \myboxnumb{0}, \myboxnumb{0}, \myboxnumb{0}, \myboxnumb{0}, \myboxnumb{0}, \myboxnumb{-1}, \myboxnumb{0}, \myboxnumb{-1}, \myboxnumb{-1}, \myboxnumb{0}, \myboxnumb{0}, \myboxnumb{0}, \myboxnumb{1}, \myboxnumb{0}, \myboxnumb{0}, \myboxnumb{1}, \myboxnumb{0}, \myboxnumb{1}, \myboxnumb{1}]_S} \fi
 \ifnum#1=91 {[
 \myboxnumb{4}, \myboxnumb{2}, \myboxnumb{3}, \myboxnumb{3}, \myboxnumb{-1}, \myboxnumb{-2}, \myboxnumb{-1}, \myboxnumb{0}, \myboxnumb{-1}, \myboxnumb{0}, \myboxnumb{0}, \myboxnumb{-1}, \myboxnumb{-1}, \myboxnumb{-1}, \myboxnumb{-1}, \myboxnumb{-1}, \myboxnumb{-1}, \myboxnumb{-1}, \myboxnumb{-1}, \myboxnumb{1}, \myboxnumb{0}, \myboxnumb{1}]_S} \fi
 \ifnum#1=92 {[
 \myboxnumb{-1}, \myboxnumb{0}, \myboxnumb{0}, \myboxnumb{0}, \myboxnumb{0}, \myboxnumb{1}, \myboxnumb{0}, \myboxnumb{0}, \myboxnumb{1}, \myboxnumb{0}, \myboxnumb{0}, \myboxnumb{1}, \myboxnumb{0}, \myboxnumb{0}, \myboxnumb{1}, \myboxnumb{0}, \myboxnumb{0}, \myboxnumb{1}, \myboxnumb{0}, \myboxnumb{-1}, \myboxnumb{-1}, \myboxnumb{-1}]_S} \fi
 \ifnum#1=93 {[
 \myboxnumb{-1}, \myboxnumb{-1}, \myboxnumb{-1}, \myboxnumb{-2}, \myboxnumb{1}, \myboxnumb{1}, \myboxnumb{1}, \myboxnumb{1}, \myboxnumb{0}, \myboxnumb{0}, \myboxnumb{-1}, \myboxnumb{-1}, \myboxnumb{0}, \myboxnumb{0}, \myboxnumb{1}, \myboxnumb{1}, \myboxnumb{0}, \myboxnumb{0}, \myboxnumb{1}, \myboxnumb{-1}, \myboxnumb{1}, \myboxnumb{1}]_S} \fi
 \ifnum#1=94 {[
 \myboxnumb{1}, \myboxnumb{2}, \myboxnumb{1}, \myboxnumb{2}, \myboxnumb{0}, \myboxnumb{1}, \myboxnumb{-1}, \myboxnumb{-1}, \myboxnumb{1}, \myboxnumb{-1}, \myboxnumb{0}, \myboxnumb{1}, \myboxnumb{0}, \myboxnumb{-1}, \myboxnumb{0}, \myboxnumb{-1}, \myboxnumb{0}, \myboxnumb{1}, \myboxnumb{0}, \myboxnumb{-1}, \myboxnumb{-2}, \myboxnumb{-1}]_S} \fi
 \ifnum#1=95 {[
 \myboxnumb{-3}, \myboxnumb{-2}, \myboxnumb{-2}, \myboxnumb{-3}, \myboxnumb{1}, \myboxnumb{1}, \myboxnumb{2}, \myboxnumb{1}, \myboxnumb{1}, \myboxnumb{1}, \myboxnumb{0}, \myboxnumb{0}, \myboxnumb{1}, \myboxnumb{1}, \myboxnumb{0}, \myboxnumb{1}, \myboxnumb{1}, \myboxnumb{0}, \myboxnumb{0}, \myboxnumb{0}, \myboxnumb{1}, \myboxnumb{-1}]_S} \fi
 \ifnum#1=96 {[
 \myboxnumb{4}, \myboxnumb{2}, \myboxnumb{3}, \myboxnumb{4}, \myboxnumb{-2}, \myboxnumb{-3}, \myboxnumb{-2}, \myboxnumb{-1}, \myboxnumb{-2}, \myboxnumb{0}, \myboxnumb{1}, \myboxnumb{0}, \myboxnumb{-1}, \myboxnumb{0}, \myboxnumb{-1}, \myboxnumb{-1}, \myboxnumb{-1}, \myboxnumb{-1}, \myboxnumb{-1}, \myboxnumb{2}, \myboxnumb{0}, \myboxnumb{1}]_S} \fi
 \ifnum#1=97 {[
 \myboxnumb{-1}, \myboxnumb{0}, \myboxnumb{-1}, \myboxnumb{-1}, \myboxnumb{0}, \myboxnumb{1}, \myboxnumb{0}, \myboxnumb{0}, \myboxnumb{1}, \myboxnumb{0}, \myboxnumb{0}, \myboxnumb{0}, \myboxnumb{0}, \myboxnumb{0}, \myboxnumb{0}, \myboxnumb{0}, \myboxnumb{1}, \myboxnumb{1}, \myboxnumb{1}, \myboxnumb{-1}, \myboxnumb{0}, \myboxnumb{0}]_S} \fi
 \ifnum#1=98 {[
 \myboxnumb{0}, \myboxnumb{0}, \myboxnumb{0}, \myboxnumb{0}, \myboxnumb{0}, \myboxnumb{0}, \myboxnumb{-1}, \myboxnumb{0}, \myboxnumb{0}, \myboxnumb{0}, \myboxnumb{0}, \myboxnumb{1}, \myboxnumb{0}, \myboxnumb{0}, \myboxnumb{1}, \myboxnumb{0}, \myboxnumb{0}, \myboxnumb{1}, \myboxnumb{0}, \myboxnumb{0}, \myboxnumb{-1}, \myboxnumb{0}]_S} \fi
 \ifnum#1=99 {[
 \myboxnumb{-1}, \myboxnumb{-1}, \myboxnumb{-1}, \myboxnumb{-1}, \myboxnumb{1}, \myboxnumb{0}, \myboxnumb{1}, \myboxnumb{1}, \myboxnumb{0}, \myboxnumb{0}, \myboxnumb{0}, \myboxnumb{0}, \myboxnumb{1}, \myboxnumb{0}, \myboxnumb{0}, \myboxnumb{1}, \myboxnumb{0}, \myboxnumb{-1}, \myboxnumb{0}, \myboxnumb{0}, \myboxnumb{1}, \myboxnumb{0}]_S} \fi
 \ifnum#1=100 {[
 \myboxnumb{3}, \myboxnumb{2}, \myboxnumb{3}, \myboxnumb{3}, \myboxnumb{-1}, \myboxnumb{-1}, \myboxnumb{0}, \myboxnumb{-1}, \myboxnumb{-1}, \myboxnumb{0}, \myboxnumb{0}, \myboxnumb{-1}, \myboxnumb{-1}, \myboxnumb{0}, \myboxnumb{-1}, \myboxnumb{-1}, \myboxnumb{-1}, \myboxnumb{-1}, \myboxnumb{-1}, \myboxnumb{1}, \myboxnumb{0}, \myboxnumb{0}]_S} \fi
 \ifnum#1=101 {[
 \myboxnumb{0}, \myboxnumb{0}, \myboxnumb{0}, \myboxnumb{-1}, \myboxnumb{1}, \myboxnumb{1}, \myboxnumb{1}, \myboxnumb{1}, \myboxnumb{0}, \myboxnumb{0}, \myboxnumb{-1}, \myboxnumb{-1}, \myboxnumb{0}, \myboxnumb{0}, \myboxnumb{0}, \myboxnumb{0}, \myboxnumb{0}, \myboxnumb{-1}, \myboxnumb{0}, \myboxnumb{0}, \myboxnumb{1}, \myboxnumb{0}]_S} \fi
 \ifnum#1=102 {[
 \myboxnumb{-1}, \myboxnumb{-1}, \myboxnumb{-1}, \myboxnumb{-1}, \myboxnumb{1}, \myboxnumb{1}, \myboxnumb{1}, \myboxnumb{0}, \myboxnumb{0}, \myboxnumb{0}, \myboxnumb{0}, \myboxnumb{0}, \myboxnumb{1}, \myboxnumb{0}, \myboxnumb{0}, \myboxnumb{1}, \myboxnumb{0}, \myboxnumb{0}, \myboxnumb{1}, \myboxnumb{-1}, \myboxnumb{0}, \myboxnumb{0}]_S} \fi
 \ifnum#1=103 {[
 \myboxnumb{0}, \myboxnumb{0}, \myboxnumb{0}, \myboxnumb{0}, \myboxnumb{-1}, \myboxnumb{-1}, \myboxnumb{0}, \myboxnumb{0}, \myboxnumb{0}, \myboxnumb{1}, \myboxnumb{1}, \myboxnumb{0}, \myboxnumb{0}, \myboxnumb{1}, \myboxnumb{0}, \myboxnumb{0}, \myboxnumb{0}, \myboxnumb{0}, \myboxnumb{-1}, \myboxnumb{1}, \myboxnumb{0}, \myboxnumb{0}]_S} \fi
 \ifnum#1=104 {[
 \myboxnumb{2}, \myboxnumb{2}, \myboxnumb{2}, \myboxnumb{3}, \myboxnumb{-1}, \myboxnumb{-1}, \myboxnumb{-2}, \myboxnumb{-1}, \myboxnumb{0}, \myboxnumb{-1}, \myboxnumb{0}, \myboxnumb{1}, \myboxnumb{-1}, \myboxnumb{-1}, \myboxnumb{0}, \myboxnumb{-1}, \myboxnumb{0}, \myboxnumb{1}, \myboxnumb{0}, \myboxnumb{0}, \myboxnumb{-1}, \myboxnumb{0}]_S} \fi
 \ifnum#1=105 {[
 \myboxnumb{2}, \myboxnumb{2}, \myboxnumb{2}, \myboxnumb{3}, \myboxnumb{-1}, \myboxnumb{-1}, \myboxnumb{-1}, \myboxnumb{-1}, \myboxnumb{-1}, \myboxnumb{-1}, \myboxnumb{0}, \myboxnumb{0}, \myboxnumb{0}, \myboxnumb{0}, \myboxnumb{0}, \myboxnumb{0}, \myboxnumb{-1}, \myboxnumb{0}, \myboxnumb{0}, \myboxnumb{0}, \myboxnumb{-1}, \myboxnumb{0}]_S} \fi
 \ifnum#1=106 {[
 \myboxnumb{3}, \myboxnumb{2}, \myboxnumb{2}, \myboxnumb{3}, \myboxnumb{-1}, \myboxnumb{-1}, \myboxnumb{-1}, \myboxnumb{0}, \myboxnumb{0}, \myboxnumb{0}, \myboxnumb{0}, \myboxnumb{0}, \myboxnumb{-1}, \myboxnumb{-1}, \myboxnumb{-1}, \myboxnumb{-1}, \myboxnumb{-1}, \myboxnumb{-1}, \myboxnumb{-1}, \myboxnumb{1}, \myboxnumb{0}, \myboxnumb{0}]_S} \fi
 \ifnum#1=107 {[
 \myboxnumb{-2}, \myboxnumb{-1}, \myboxnumb{-1}, \myboxnumb{-2}, \myboxnumb{1}, \myboxnumb{2}, \myboxnumb{1}, \myboxnumb{0}, \myboxnumb{1}, \myboxnumb{0}, \myboxnumb{-1}, \myboxnumb{0}, \myboxnumb{0}, \myboxnumb{0}, \myboxnumb{1}, \myboxnumb{0}, \myboxnumb{1}, \myboxnumb{1}, \myboxnumb{1}, \myboxnumb{-1}, \myboxnumb{0}, \myboxnumb{0}]_S} \fi
 \ifnum#1=108 {[
 \myboxnumb{-2}, \myboxnumb{-2}, \myboxnumb{-2}, \myboxnumb{-3}, \myboxnumb{1}, \myboxnumb{0}, \myboxnumb{1}, \myboxnumb{1}, \myboxnumb{0}, \myboxnumb{1}, \myboxnumb{1}, \myboxnumb{0}, \myboxnumb{1}, \myboxnumb{1}, \myboxnumb{0}, \myboxnumb{1}, \myboxnumb{1}, \myboxnumb{0}, \myboxnumb{0}, \myboxnumb{0}, \myboxnumb{1}, \myboxnumb{0}]_S} \fi
 \ifnum#1=109 {[
 \myboxnumb{0}, \myboxnumb{-1}, \myboxnumb{0}, \myboxnumb{0}, \myboxnumb{0}, \myboxnumb{-1}, \myboxnumb{0}, \myboxnumb{0}, \myboxnumb{0}, \myboxnumb{1}, \myboxnumb{1}, \myboxnumb{1}, \myboxnumb{0}, \myboxnumb{0}, \myboxnumb{0}, \myboxnumb{0}, \myboxnumb{0}, \myboxnumb{0}, \myboxnumb{-1}, \myboxnumb{1}, \myboxnumb{0}, \myboxnumb{0}]_S} \fi
 \ifnum#1=110 {[
 \myboxnumb{-1}, \myboxnumb{0}, \myboxnumb{0}, \myboxnumb{-1}, \myboxnumb{0}, \myboxnumb{0}, \myboxnumb{1}, \myboxnumb{0}, \myboxnumb{0}, \myboxnumb{0}, \myboxnumb{0}, \myboxnumb{-1}, \myboxnumb{0}, \myboxnumb{1}, \myboxnumb{0}, \myboxnumb{0}, \myboxnumb{1}, \myboxnumb{0}, \myboxnumb{0}, \myboxnumb{0}, \myboxnumb{1}, \myboxnumb{0}]_S} \fi
 \ifnum#1=111 {[
 \myboxnumb{4}, \myboxnumb{3}, \myboxnumb{3}, \myboxnumb{4}, \myboxnumb{-1}, \myboxnumb{-1}, \myboxnumb{-2}, \myboxnumb{-1}, \myboxnumb{-1}, \myboxnumb{-1}, \myboxnumb{0}, \myboxnumb{0}, \myboxnumb{-1}, \myboxnumb{-1}, \myboxnumb{-1}, \myboxnumb{-1}, \myboxnumb{-1}, \myboxnumb{0}, \myboxnumb{0}, \myboxnumb{0}, \myboxnumb{-1}, \myboxnumb{0}]_S} \fi
 \ifnum#1=112 {[
 \myboxnumb{-2}, \myboxnumb{-1}, \myboxnumb{-2}, \myboxnumb{-2}, \myboxnumb{1}, \myboxnumb{2}, \myboxnumb{1}, \myboxnumb{1}, \myboxnumb{1}, \myboxnumb{0}, \myboxnumb{-1}, \myboxnumb{0}, \myboxnumb{1}, \myboxnumb{0}, \myboxnumb{1}, \myboxnumb{1}, \myboxnumb{0}, \myboxnumb{0}, \myboxnumb{1}, \myboxnumb{-1}, \myboxnumb{0}, \myboxnumb{0}]_S} \fi
}

\begin{document}
\title[The automorphism group of  a supersingular $K3$ surface]%
{The automorphism group of  a supersingular $K3$ surface with Artin invariant $1$ in characteristic $3$}

\author{Shigeyuki Kond\={o}}
\address{
Graduate School of Mathematics, 
Nagoya University, 
Nagoya,
 464-8602 JAPAN
}
\email{kondo@math.nagoya-u.ac.jp
}

\author{Ichiro Shimada}
\address{
Department of Mathematics, 
Graduate School of Science, 
Hiroshima University,
1-3-1 Kagamiyama, 
Higashi-Hiroshima, 
739-8526 JAPAN
}
\email{shimada@math.sci.hiroshima-u.ac.jp}
\thanks{The first author was partially supported by
JSPS Grant-in-Aid  for Scientific Research (S)
  No.2222400.
The second author was partially supported by JSPS Grants-in-Aid for Scientific Research (B) No.20340002. 
}




\begin{abstract}
We present  a finite  set  of generators of the automorphism group  of 
a supersingular $K3$ surface with Artin invariant $1$ in characteristic $3$.
\end{abstract}

\maketitle
\section{Introduction}\label{sec:intro}
To determine the automorphism group $\Aut(Y)$ of a given $K3$ surface $Y$
is an important problem.
In this paper, we present a set of generators of  the automorphism group 
of a  supersingular $K3$ surface $X$ in characteristic $3$ with Artin invariant $1$.
Our method is computational, and relies heavily on computer-aided calculation.
It gives us generators in explicit form,
and it can be easily applied to many other $K3$ surfaces 
by modifying  computer programs.
\par
\medskip
A $K3$ surface defined over an algebraically closed field $k$ is said to be
\emph{supersingular} (in the sense of Shioda) if its Picard number is $22$.
Supersingular $K3$ surfaces exist only when $k$ is of positive characteristic.
Let $Y$ be a supersingular $K3$ surface in characteristic $p>0$,
and let $S_Y$ denote its N\'eron-Severi lattice.
Artin~\cite{MR0371899} showed that 
the discriminant group of  $S_Y$ 
is a $p$-elementary abelian group of rank $2\sigma$,  
where $\sigma$ is an integer such that $1\le \sigma\le 10$.
This integer $\sigma$ is called the \emph{Artin invariant} of $Y$.
Ogus~\cite{MR563467, MR717616} proved  that a supersingular $K3$ surface with Artin invariant $1$ 
in characteristic $p$ 
is unique up to isomorphisms
(see also~\cite{MR633161}).
\par
\medskip
It is known that the Fermat quartic surface 
$$
X:=\{w^4+x^4+y^4+z^4=0\}\subset \P^3
$$
defined over an algebraically closed field $k$
of  characteristic $3$ is a  supersingular $K3$ surface  
with Artin invariant $1$~(see~\cite{MR918849}).
Let 
$$
h_0:=[\OOO_X(1)]\in S_X
$$
denote the class of the hyperplane section of $X$.
The projective automorphism group
$\Aut(X, h_0)$ of $X\subset \P^3$ is equal to the finite subgroup $\PGU_4(\F_9)$ of $\PGL_4(k)$
with order $13, 063,680$.
\par
\medskip
Let $(w,x,y)$ be the affine coordinates of $\P^3$ with $z=1$, 
and let  $F_{1j}$ and $F_{2j}$ be polynomials of $(w,x,y)$
with coefficients in 
$$
\F_9=\F_3(i)=\{0,\pm 1, \pm i, \pm (1+i), \pm(1-i)\}, \quad\textrm{where $i:=\sqrt{-1}$, }
$$
given in Table~\ref{table:Fs}. 
\begin{proposition}\label{prop:phi}
For $i=1$ and $2$, 
the rational map
\begin{equation*}\label{isom:FG}
(w,x,y)\mapsto [F_{i0}: F_{i1}: F_{i2}]\in \P^2
\end{equation*}
induces a morphism $\phi_i : X\to \P^2$ of degree $2$.
\end{proposition}
We denote by
$$
X\;\maprightsp{\psi_i}\; Y_i\;\maprightsp{\pi_i}\; \P^2
$$ 
the Stein factorization of $\phi_i :X\to \P^2$,
and let $B_i\subset\P^2$ be the branch curve of 
the finite morphism $\pi_i :Y_i\to \P^2$ of degree $2$.
Note that $Y_i$ is a normal $K3$ surface,
and hence $Y_i$ has only rational double points as its singularities~(see~\cite{MR0146182, MR0199191}).
Let $[x_0:x_1:x_2]$ be the homogeneous coordinates of $\P^2$.
\begin{proposition}\label{prop:sectics}
{\rm (1)}
The $ADE$-type of  the singularities of $Y_1$  is $6A_1+4A_2$.
The branch curve $B_1$ is defined by $f_1=0$, where 
\begin{dmath*}
f_{1}{:=}
 \, {x_0}^{6}+\, {x_0}^{5}{x_1}-\, {x_0}^{3}{x_1}^{3}-\, {x_0}{x_1}^{5}-\, {x_0}^{4}{x_2}^{2}+\, {x_0}{x_1}^{3}{x_2}^{2}+\, {x_1}^{4}{x_2}^{2}+\, {x_0}^{2}{x_2}^{4}+\, {x_1}^{2}{x_2}^{4}+\, {x_2}^{6}.
\end{dmath*}
{\rm (2)}
The $ADE$-type of  the singularities of $Y_2$  is $A_1+A_2+2A_3+2A_4$.
The branch curve $B_2$ is defined by $f_2=0$, where 
\begin{dmath*}
f_{2}{:=}
 \, {x_0}^{5}{x_1}+\, {x_0}^{2}{x_1}^{4}-\, {x_0}^{4}{x_2}^{2}+\, {x_0}{x_1}^{3}{x_2}^{2}+\, {x_1}^{4}{x_2}^{2}-\, {x_0}^{2}{x_2}^{4}-\, {x_0}{x_1}{x_2}^{4}-\, {x_1}^{2}{x_2}^{4}-\, {x_2}^{6}.
\end{dmath*}

\end{proposition}
\begin{table}
\begin{dgroup*}
\begin{dmath*}
F_{10}=
(1+\sqrmo)\, +(1+\sqrmo)\, w+(1-\sqrmo)\, x-\, y-(1-\sqrmo)\, wx-\, x^{2}+\sqrmo\, wy+\sqrmo\, xy-\sqrmo\, y^{2}+(1+\sqrmo)\, w^{3}-\sqrmo\, w^{2}x+(1+\sqrmo)\, wx^{2}-\sqrmo\, x^{3}+\, w^{2}y+(1+\sqrmo)\, wxy+(1+\sqrmo)\, x^{2}y-(1-\sqrmo)\, wy^{2}-(1+\sqrmo)\, xy^{2}+\sqrmo\, y^{3}
\end{dmath*}
\begin{dmath*}
F_{11}=
(1-\sqrmo)\, -(1+\sqrmo)\, x-(1-\sqrmo)\, y-(1-\sqrmo)\, w^{2}-(1-\sqrmo)\, wx-(1-\sqrmo)\, x^{2}-(1+\sqrmo)\, wy-\, xy-(1+\sqrmo)\, y^{2}-\, w^{3}+(1-\sqrmo)\, w^{2}x+\, wx^{2}-\sqrmo\, x^{3}-(1+\sqrmo)\, w^{2}y-(1+\sqrmo)\, wxy+\, x^{2}y-\sqrmo\, wy^{2}-\, xy^{2}+(1-\sqrmo)\, y^{3}
\end{dmath*}
\begin{dmath*}
F_{12}=
(1+\sqrmo)\, w-\sqrmo\, x-\, y-\, w^{2}-\, wx-\sqrmo\, x^{2}-\sqrmo\, xy+\sqrmo\, y^{2}+\sqrmo\, w^{3}-(1+\sqrmo)\, wx^{2}+\sqrmo\, x^{3}-\sqrmo\, w^{2}y-\, wxy+(1-\sqrmo)\, wy^{2}+(1+\sqrmo)\, y^{3}
\end{dmath*}
\end{dgroup*}
\vskip 4pt
--------------------------------------
\vskip 1pt
\begin{dgroup*}
\begin{dmath*}
F_{20}=
-1\, -\sqrmo\, w+(1+\sqrmo)\, x-\, y-(1+\sqrmo)\, w^{2}-\, wx-(1-\sqrmo)\, x^{2}-\sqrmo\, wy+(1+\sqrmo)\, xy-(1-\sqrmo)\, w^{3}+\, w^{2}x-\, wx^{2}+\, x^{3}-\, w^{2}y+(1-\sqrmo)\, wxy+\, x^{2}y+(1-\sqrmo)\, wy^{2}+(1-\sqrmo)\, xy^{2}+(1+\sqrmo)\, y^{3}-\, w^{3}x-\sqrmo\, w^{2}x^{2}-\, wx^{3}+\, w^{3}y-(1+\sqrmo)\, w^{2}xy-(1-\sqrmo)\, wxy^{2}+\, x^{2}y^{2}-(1-\sqrmo)\, wy^{3}-(1+\sqrmo)\, xy^{3}-\, y^{4}+(1-\sqrmo)\, w^{3}x^{2}-\sqrmo\, x^{5}+(1-\sqrmo)\, w^{3}xy+(1+\sqrmo)\, wx^{3}y-\sqrmo\, w^{3}y^{2}+(1+\sqrmo)\, w^{2}xy^{2}-(1+\sqrmo)\, wx^{2}y^{2}+\sqrmo\, x^{3}y^{2}-\, w^{2}y^{3}-(1+\sqrmo)\, wxy^{3}-(1-\sqrmo)\, x^{2}y^{3}+\sqrmo\, wy^{4}+(1-\sqrmo)\, xy^{4}+(1+\sqrmo)\, y^{5}
\end{dmath*}
\begin{dmath*}
F_{21}=
-(1-\sqrmo)\, +\sqrmo\, w+(1-\sqrmo)\, y-(1+\sqrmo)\, w^{2}+\, wx+(1+\sqrmo)\, x^{2}+(1+\sqrmo)\, wy-(1+\sqrmo)\, xy-\sqrmo\, y^{2}-\, w^{3}+\sqrmo\, w^{2}x+(1+\sqrmo)\, wx^{2}-\, x^{3}-(1+\sqrmo)\, w^{2}y-(1-\sqrmo)\, wxy-(1-\sqrmo)\, x^{2}y-\sqrmo\, wy^{2}-(1+\sqrmo)\, xy^{2}+\, y^{3}-(1-\sqrmo)\, w^{3}x-\, wx^{3}+(1-\sqrmo)\, x^{4}+(1-\sqrmo)\, w^{3}y+\sqrmo\, w^{2}xy+(1-\sqrmo)\, wx^{2}y-\sqrmo\, x^{3}y+(1-\sqrmo)\, w^{2}y^{2}+(1-\sqrmo)\, wxy^{2}-(1+\sqrmo)\, x^{2}y^{2}+(1-\sqrmo)\, wy^{3}-\sqrmo\, xy^{3}+\sqrmo\, y^{4}+\, w^{3}x^{2}+\, w^{2}x^{3}+(1-\sqrmo)\, wx^{4}-\sqrmo\, x^{5}-\sqrmo\, w^{3}xy+\, w^{2}x^{2}y+(1+\sqrmo)\, wx^{3}y+\, x^{4}y+\, w^{3}y^{2}-\, w^{2}xy^{2}-\, wx^{2}y^{2}+\sqrmo\, w^{2}y^{3}+(1+\sqrmo)\, wxy^{3}-\sqrmo\, wy^{4}-\sqrmo\, xy^{4}+\, y^{5}
\end{dmath*}
\begin{dmath*}
F_{22}=
(1-\sqrmo)\, -(1+\sqrmo)\, w-(1+\sqrmo)\, x-(1-\sqrmo)\, y+\sqrmo\, w^{2}-(1+\sqrmo)\, wx-(1-\sqrmo)\, x^{2}+\sqrmo\, wy-(1+\sqrmo)\, xy-\, w^{3}-\sqrmo\, w^{2}x-\, wx^{2}+\, x^{3}-(1-\sqrmo)\, w^{2}y+\, wxy+\, x^{2}y+(1+\sqrmo)\, wy^{2}-(1+\sqrmo)\, xy^{2}-\, y^{3}+\sqrmo\, w^{3}x-(1-\sqrmo)\, w^{2}x^{2}-\, wx^{3}-(1+\sqrmo)\, x^{4}+\sqrmo\, w^{3}y+\, w^{2}xy+(1-\sqrmo)\, wx^{2}y-(1-\sqrmo)\, w^{2}y^{2}+(1+\sqrmo)\, wxy^{2}+\sqrmo\, wy^{3}+\, xy^{3}+(1-\sqrmo)\, y^{4}-\sqrmo\, w^{3}x^{2}-(1+\sqrmo)\, wx^{4}+\, x^{5}-(1-\sqrmo)\, w^{3}xy-\sqrmo\, w^{2}x^{2}y+(1+\sqrmo)\, wx^{3}y+(1-\sqrmo)\, x^{4}y-\, w^{3}y^{2}-(1+\sqrmo)\, w^{2}xy^{2}+\sqrmo\, wx^{2}y^{2}+\sqrmo\, x^{3}y^{2}-\, wxy^{3}-(1-\sqrmo)\, x^{2}y^{3}-\, wy^{4}-\, xy^{4}-\, y^{5}
\end{dmath*}
\end{dgroup*}
\vskip 4pt
\caption{Polynomials $F_{1j}$ and $F_{2j}$} \label{table:Fs} 
\end{table}
Our main result is as follows:
\begin{theorem}\label{thm:main}
Let $g_i\in \Aut(X)$ denote the involution
induced  from the deck-transformation of $\pi_i: Y_i\to\P^2$.
Then $\Aut(X)$ is generated by $\Aut(X, h_0)= \PGU_4(\F_9)$ and $g_1$, $g_2$.
\end{theorem}
See Theorem~\ref{thm:H} for a more explicit description of the
involutions~$g_1$ and $g_2$.
\par
\medskip

Let $\cone_{S_X}$ denote  the  connected component of $\shortset{x\in S_X\tR}{x^2>0}$
that contains $h_0$.
Following Borcherds~\cite{MR913200}, 
we  prove Theorem~\ref{thm:main}  by calculating a closed chamber $\chamsz$ in the cone  $\cone_{S_X}$ 
with the following properties (see~Section~\ref{sec:mainproof}):
\begin{enumerate}
\item The chamber $\chamsz$ is invariant  under the action of $\Aut(X, h_0)$.
\item For any nef class $v\in S_X$,
there exists $\gamma\in \Aut(X)$ such that $v^\gamma\in \chamsz$.
\item For nef classes $v, v\sprime$
in the interior of $\chamsz$,
there exists $\gamma\in \Aut(X)$ such that $v\sprime=v^\gamma$
if and only if there exists $\tau\in \Aut(X, h_0)$ 
such that $v\sprime=v^{\tau}$.
\end{enumerate}
This chamber $\chamsz$  is bounded by $112+648+5184$ hyperplanes 
in  $\cone_{S_X}$.
See Proposition~\ref{prop:CS0} for the explicit description of these walls.
Using $\chamsz$ and these walls, 
we can also present a finite set of generators of $\OG^+(S_X)$
(see Theorem~\ref{thm:OplusS}).
\par
\medskip
Vinberg~\cite{MR719348}
determined the automorphism groups 
of two complex $K3$ surfaces with Picard number $20$
by investigating the orthogonal groups
of their N\'eron-Severi lattices and 
the associated hyperbolic geometry.
\par
\medskip
Let $L$ denote an
even unimodular  lattice  of rank $26$
with signature $(1, 25)$,
which is unique up to isomorphisms by Eichler's theorem.
Conway~\cite{MR690711} determined the fundamental domain 
in a positive cone  of $L\tR$ under the action of
the subgroup of $\OG^+(L)$
generated by the reflections with respect to the vectors 
of square norm $-2$.
Borcherds~\cite{MR913200} applied Conway theory to
the investigation of 
the orthogonal groups  of  even hyperbolic lattices $S$ 
primitively embedded in  $L$.
Then the first author~\cite{MR1618132}
determined the automorphism group
of a generic Jacobian Kummer surface
by embedding its N\'eron-Severi lattice
into $L$
and using Conway theory.
Keum and the first author~\cite{MR1806732}
applied this method to
the Kummer surface
of the product of two elliptic curves, 
Dolgachev and Keum~\cite{MR1897389}
applied it to quartic Hessian surfaces,
and Dolgachev and the first author~\cite{MR1935564K}
applied it to
the supersingular $K3$ surface 
in characteristic $2$ with Artin invariant $1$.
\par
\medskip
Recently, 
 configurations of smooth rational curves
 on our supersingular $K3$ surface $X$ 
was studied in~\cite{KatsuraKondochar3}
with respect to an embedding of $S_X$ into $L$,
and elliptic fibrations on $X$ was classified in~\cite{11051715}
 by  embedding  $S_X$ into $L$. 
\par
\medskip
The new idea introduced in this paper is that,
in order to find automorphisms of $X$ necessary to generate $\Aut(X)$,
we search for polarizations of degree $2$
whose classes are located on the walls 
of the chamber decomposition of the cone $\cone_{S_X}$. 
The computational tools 
used in this paper have been  developed by the second author
for the study~\cite{shimadapreprintchar5} of various double plane models of 
a supersingular $K3$ surface 
in characteristic $5$ with Artin invariant $1$.
The computational data for this paper is available from 
the second author's webpage~\cite{shimadawebpage}.
\par
\medskip
In~\cite{MR2059747} and~\cite{MR2036331},
the second author showed that 
every supersingular $K3$ surface in any characteristic 
with arbitrary Artin invariant
is birational to a  double cover of the projective plane. 
In~\cite{MR2129248},~\cite{MR2296438} and~\cite{MR2282430, shimadapreprintchar5},
projective models of supersingular $K3$ surfaces 
in characteristic $2$, $3$ and $5$ were investigated,
respectively.
\par
\medskip
This paper is organized as follows.
In Section~\ref{sec:Leech}, we give  a  review of the theory of Conway and Borcherds,
 and investigate 
 chamber decomposition induced on a positive cone of a primitive hyperbolic sublattice $S$ of $L$.
In Section~\ref{sec:NS},
we give explicitly a basis of the N\'eron-Severi lattice $S_X$ of $X$,
and describe a method to compute the action of $\Aut(X, h_0)$ on $S_X$.
The fact that $S_X$ is generated by the classes of lines in $X$ enables us to calculate 
projective models of $X$ explicitly.
In Section~\ref{sec:embS},
we embed $S_X$ into $L$, and study the obtained chamber decomposition
in detail.
In particular,
we investigate the walls of the chamber $\chamsz$ that contains the class $h_0$.
In Section~\ref{sec:A},
we prove Propositions~\ref{prop:phi} and~\ref{prop:sectics},
and show that the involutions $g_1$ and $g_2$
map $h_0$ to its mirror images into   walls of  the chamber $\chamsz$.
Then we can prove Theorem~\ref{thm:main} in Section~\ref{sec:mainproof}.
In Section~\ref{sec:quartics},
we give another description of the involutions $g_i$.
In the last section,
we give a set of generators of $\OG^+(S_X)$.
\par
\medskip
Thanks are due to Professor Daniel Allcock, Professor Toshiyuki Katsura and 
Professor JongHae Keum  for helpful discussions.
\section{Leech roots}\label{sec:Leech}
\subsection{Terminologies and notation}\label{subsec:terms}
We fix some terminologies and notation about lattices.
A \emph{lattice} $M$ is a free $\Z$-module of finite rank
with a non-degenerate symmetric bilinear form
$$
\intM{\phantom{a}, \phantom{a}}: M\times M\to \Z.
$$
A submodule $N$ of $M$ is said to be \emph{primitive} if $M/N$ is torsion free.
For a submodule $N$ of $M$,
we denote by $N\sperp \subset M$ the submodule defined by
$$
N\sperp:=\set{u\in M}{\intM{u, v}=0\;\;\textrm{for all}\;\; v\in N}, 
$$
which is primitive by definition.
We denote by $\OG(M)$ the orthogonal group of $M$.
Throughout this paper, 
we let $\OG(M)$ act  on $M$ from \emph{right}. 
Suppose that  $M$ is of rank $r$.
We say  
that $M$ is \emph{hyperbolic} (resp.~\emph{negative-definite})
if the signature of the symmetric bilinear form $\intM{\phantom{a}, \phantom{a}}$ on $M\tR$
is $(1, r-1)$  (resp.~$(0, r)$).
We define the \emph{dual lattice} $M\dual$ of $M$ by
$$
M\dual:=\set{u\in M\tQ}{\intM{u, v}\in \Z\;\;\textrm{for all}\;\; v\in M}.
$$
Then $M$ is contained in $M\dual$ as a submodule of finite index.
The finite abelian group $M\dual/M$ is called the \emph{discriminant group} of $M$.
We say that $M$ is \emph{unimodular} if $M=M\dual$.
\par
A lattice $M$ is said to be \emph{even} if $\intM{v,v}\in 2\Z$ holds for any $v\in M$.
The discriminant group $M\dual/M$ 
of an even lattice $M$ is naturally equipped with 
the quadratic form
$$
q_M : M\dual/M\to \Q/2\Z
$$
defined by $q_M(u \bmod M):=\intM{u,u} \bmod 2\Z$.
We call $q_M$ the \emph{discriminant form} of $M$.
The automorphism group of $q_M$ is denoted by $\OG(q_M)$.
There exists a natural homomorphism $\OG(M)\to \OG(q_M)$.
\par
Suppose that $M$ is hyperbolic.
Then the open subset
$$
\set{x\in M\tR}{\intM{x, x}>0}
$$
of $M\tR$ has two connected components.
A \emph{positive cone}  of $M$ is one of them.
We  fix a positive cone $\cone$.
The \emph{autochronous orthogonal group} $\OG^+(M)$ of $M$ is the group of isometries of $M$
that preserve $\cone$.
Then $\OG^+(M)$ is a subgroup of  $\OG(M)$ with index $2$. 
Note that  $\OG^+(M)$ acts on $\cone$.
For a nonzero vector $u\in M\tR$, we denote by $(u)\sperp_M$ the hyperplane of $M\tR$
defined by 
$$
(u)\sperp_M:=\set{x\in M\tR}{\intM{x, u}=0}.
$$
Let
$\RRR$ be a set of non-zero vectors of $M\tR$,
and let
$$
\HHH:=\set{(u)\sperp_M}{u\in \RRR}
$$
be the family of hyperplanes defined by $\RRR$.
Suppose that $\HHH$ is locally finite in $\cone$.
Then the closure in $\cone$ of each connected component 
of 
$$
\cone\setminus \left(\cone\cap\bigcup_{u\in \RRR} (u)\sperp_M\right)
$$
is called an \emph{$\RRR$-chamber}.
Let $\cham$ be an $\RRR$-chamber.
We denote by $\cham\spcirc$ the interior of $\cham$.
We say that a hyperplane $(u)\sperp_M \in \HHH$ \emph{bounds $\cham$},
or that $(u)\sperp_M$ is a \emph{wall} of $\cham$,
if $(u)\sperp_M\cap \cham$ contains a non-empty open subset of $(u)\sperp_M$.
We denote the set of walls of $\cham$ by
$$
\wall(\cham):=\set{(u)\sperp_M\in \HHH}{\textrm{$(u)\sperp_M$ bounds $\cham$}}.
$$
Suppose that $\RRR$ is invariant under $u\mapsto -u$.
We choose  a point $p\in\cham\sp{\circ}$, and put 
$$
\wallrs(\cham):=\set{u\in \RRR}{\textrm{$(u)\sperp_M$ bounds $\cham$ and $\intM{u, p}>0$}},
$$
which is independent of the choice of $p$.
It is obvious that $\cham$ is equal to
$$
\set{x\in \cone}{\intM{x, u}\ge 0 \;\;\textrm{for all}\;\; u\in \wallrs(\cham)}.
$$

\subsection{Conway theory}\label{subsec:Conway}
We review the theory of Conway~\cite{MR690711}.
Let $L$ be an even unimodular hyperbolic lattice of rank $26$,
which is unique up to isomorphisms  by Eichler's theorem~(see, for example,~\cite[Chapter 11, Theorem 1.4]{MR522835}).
We choose and fix a positive cone  $\cone_L$ once and for all.
A vector $r\in L$ is called a \emph{root}
if the reflection
$s_r : L\tR\to L\tR$
defined by
$$
x\mapsto  x-\frac{2\intL{x, r}}{\intL{r, r}}\,\cdot\,  r
$$
preserves $L$ and $\cone_L$,
or equivalently,
if $\intL{r, r}=-2$.  
We denote by $\RRR_L$ the set of roots of $L$,
which is invariant under $r\mapsto -r$.
Let  $W(L)$ denote the subgroup of $\OG^+ (L)$ generated by the reflections $s_r$ 
associated with all the roots $r\in \RRR_L$.
Then $W(L)$ is a normal subgroup of $\OG^+ (L)$.
The family of hyperplanes 
$$
\HHH_L:=\set{(r)_L\sperp}{r\in \RRR_L} 
$$ 
is locally finite in $\cone_L$.
Hence we can consider $\RRR_L$-chambers.
By definition, 
each $\RRR_L$-chamber is a fundamental domain of the action of $W(L)$ on $\cone_L$.
\par
\medskip
A non-zero primitive vector $w\in L$ is called a \emph{Weyl vector}
if  $\intL{w, w}=0$, $w$ is contained in the closure of $\cone_L$ in $L\tR$, 
and the negative-definite even unimodular lattice $\gen{w}\sperp/\gen{w}$  of rank $24$ has no vectors 
of square norm $-2$.
Let $w\in L$ be a Weyl vector.
We put
$$
\LR(w):=\set{r\in \RRR_L}{\intL{w, r}=1}.
$$
A root in $\LR(w)$ is called a \emph{Leech root with respect to $w$}.
\par
\medskip
Suppose that 
$w$ is a non-zero primitive vector of norm $0$ 
contained in the closure of $\cone_L$.
Then there exists a vector $w\sprime\in L$ such that
$\intL{w, w\sprime}=1$ and  $\intL{w\sprime, w\sprime}=0$.
Let $U\subset L$ denote  the hyperbolic sublattice of rank $2$ generated by $w$ and  $w\sprime$.  
By  Niemeier's classification~\cite{MR0316384} of even definite unimodular lattices of 
rank $24$~(see also~\cite[Chapter 18]{MR1662447}),
we see that 
the condition that $\gen{w}\sperp/\gen{w}$ have no vectors 
of square norm $-2$ is equivalent to the condition that 
the orthogonal complement $U\sperp$ of $U$ in $L$
be isomorphic to the 
(negative-definite) Leech lattice $\Lambda$.
From this fact, we can deduce the following:
\begin{proposition}
The group $\OG^+(L)$ acts on the set of Weyl vectors transitively.
\end{proposition}
\begin{proposition}\label{prop:lambda}
Suppose that $w$ is a Weyl vector and 
that $w\sprime\in L$ satisfies $\intL{w, w\sprime}=1$ and  $\intL{w\sprime, w\sprime}=0$. 
Via an  isomorphism $\rho:\Lambda\isom U\sperp$, the map
$$
\lambda\;\;\mapsto\;\; -\,\frac{2+(\lambda, \lambda)_\Lambda}{2} w+ w\sprime +\rho(\lambda)
$$
induces a bijection from the Leech lattice $\Lambda$ to the set $\LR(w)$.
\end{proposition}
Using Vinberg's algorithm~\cite{MR0422505}
and the result on the covering radius of the Leech lattice~\cite{MR660415}, 
Conway~\cite{MR690711} proved the following:
\begin{theorem}\label{thm:Conway}
Let $w\in L$ be a Weyl vector.
Then 
$$
\cham_L(w):=\set{x\in \cone_L}{\textrm{$\intL{x, r}\ge 0$ for all $r\in\LR(w)$}}
$$ 
is an $\RRR_L$-chamber,
and $\wallrs(\cham_L(w))$ is equal to $\LR(w)$;
that is, $(r)_L\sperp$ bounds $\cham_L(w)$ for any $r\in \LR(w)$.
The map $w\mapsto \cham_L(w)$ is a bijection from the set of Weyl vectors to the set of $\RRR_L$-chambers.
\end{theorem}
\begin{remark}
Using Proposition~\ref{prop:lambda}, Conway~\cite{MR690711} also showed that the automorphism group
$\Aut(\cham_L(w))\subset \OG^+(L)$ of an $\RRR_L$-chamber $\cham_L(w)$ is
isomorphic to the group $\cdot \infty$  of \emph{affine} automorphisms of the Leech lattice $\Lambda$.
Hence $\OG^+(L)$ is isomorphic to 
the split extension of $\cdot\infty$ by  $W(L)$.
\end{remark}
\subsection{Restriction of $\RRR_L$-chambers to a primitive sublattice}
\label{subsec:restriction}

Let $S$ be an even hyperbolic lattice of rank $r<26$ primitively embedded in $L$.
Following Borcherds~\cite{MR913200},
we explain how the Leech roots of $L$ induce a chamber decomposition on the positive cone  
$$
\cone_S:=\cone_L\cap (S\tR)
$$
of $S\tR$.
\par
\medskip
The orthogonal complement $T:=S\sperp$ of $S$ in $L$ is negative-definite of rank $26-r$,
and we have
$$
S\oplus T \;\subset\; L \;\subset\; S\dual\oplus T\dual
$$
with $[L:S\oplus T]=[S\dual\oplus T\dual: L] $.
The projections $L\tR\to S\tR$ and $L\tR\to T\tR$ are denoted by
$$
x\mapsto x_S\quand x\mapsto x_T, 
$$
respectively.
Note that, if $v\in L$,  then $v_S\in S\dual$ and $v_T\in T\dual$.
\par
Let $r\in L$ be a root.
Then the hyperplane  $(r)\sperp_L$ contains $S\tR$ if and only if $r_S=0$,
or equivalently, $r\in T$.
Since $T$ is negative-definite,
the set
$$
\RRR_T:=\set{v\in T}{\intT{v,v}=-2}
$$
is finite, and therefore there exist only finite number of 
hyperplanes $(r)\sperp_L$ that contain $S\tR$.
Suppose that $r_S\ne 0$.
If  $\intS{r_S, r_S}\ge 0$,
then either $\cone_S$ is entirely contained in the interior of the halfspace 
$$
\set{x\in L\tR}{\intL{x, r}\ge 0}
$$
or is disjoint from this halfspace.
Hence the hyperplane 
$$
(r_S)\sperp_S=(r)\sperp_L\cap(S\tR)
$$
of $S\tR$ intersects $\cone_S$ if and only if $\intS{r_S, r_S}<0$.
We put
\begin{eqnarray*}
\RRR_S&:=&\set{r_S}{r\in \RRR_L\;\;\textrm{and}\;\; \intS{r_S, r_S}<0}\\
&=&\set{r_S}{r\in \RRR_L \;\;\textrm{and}\;\; (r_S)\sperp_S \cap \cone_S\ne \emptyset}.
\end{eqnarray*}
Then the  associated family of hyperplanes 
$$
\HHH_S:=\set{(r_S)\sperp_S}{r_S\in \RRR_S}
$$
is locally finite in $\cone_S$,
and hence we can consider  $\RRR_S$-chambers in $\cone_S$.
Note that $\RRR_S$ is invariant under $r_S\mapsto -r_S$.
We investigate the relation between $\RRR_S$-chambers
and $\RRR_L$-chambers.
\par
If $\cham_S\subset \cone_S$ is an $\RRR_S$-chamber,
then there exists an $\RRR_L$-chamber $\cham_L(w)\subset \cone_L$ such that
$\cham_S=\cham_L(w)\cap(S\tR)$ holds.

For a given $\RRR_S$-chamber $\cham_S$, 
the set of $\RRR_L$-chambers $\cham_L(w)$ satisfying 
$\cham_S=\cham_L(w)\cap(S\tR)$
is in one-to-one correspondence with 
the set of connected components of
$$
(T\tR )\setminus \bigcup_{r\in \RRR_T} (r)\sperp_T.
$$
Conversely, suppose that  an $\RRR_L$-chamber $\cham_L(w)$ is given.
%
\begin{definition}\label{def:Snondeg}
We say that $\cham_L(w)$ is \emph{$S$-nondegenerate}
if $\cham_L(w)\cap(S\tR)$ is an $\RRR_S$-chamber.
\end{definition}
By definition,
$\cham_L(w)$ is $S$-nondegenerate if and only if
$w$ satisfies the following two conditions:
\begin{itemize}
\item[(i)]
There exists $v\in \cone_S$ such that
$\intL{v, r}\ge 0$ holds for any $r\in \LR(w)$. 
\item[(ii)]
There exists $v\sprime\in \cone_S$ such that
$\intL{v\sprime, r}> 0$ holds for any $r\in \LR(w)$ with $\intS{r_S, r_S}< 0$.
\end{itemize}
If $\cham_S=\cham_L(w)\cap(S\tR)$ is an $\RRR_S$-chamber,
then $\wallrs(\cham_S)$ is contained in the image of the set
$$
\LR(w, S):=\set{r\in \LR(w)}{r_S\in \RRR_S}
= \set{r\in \LR(w)}{\intS{r_S, r_S}<0}   
$$
by the projection $L\to S\dual$.
The following proposition shows that 
$\cham_S$ is bounded by a finite number
of walls  if $w_T\ne 0$, and 
its proof indicates an effective procedure to calculate $\LR(w, S)$.
( See~\cite[Section 3]{shimadapreprintchar5} for the details of the necessary algorithms.)
\begin{proposition}\label{prop:LRS}
Let $w\in L$  be a Weyl vector such that $w_T\ne 0$.
Then  $\LR(w, S)$
is a  finite set.
\end{proposition}
\begin{proof}
Since $T$ is negative-definite and $w_T\ne 0$, we have 
$$
\intS{w_S, w_S}=-\intT{w_T, w_T}>0.
$$
Suppose that $r\in \LR(w)$.
Then we have
$$
\intS{w_S, r_S}+\intT{w_T, r_T}=1,
\quad 
\intS{r_S, r_S}+\intT{r_T, r_T}=-2.
$$
We have $\intS{r_S, r_S}<0$ if and only if $\intT{r_T, r_T}>-2$.
Since $T$ is negative-definite,
the set
$$
V_T:=\set{v\in T\dual }{\intT{v, v}>-2}
$$
is finite.
For $v\in V_T$, we put
$$
a_v:=1-\intT{w_T, v},
\quad
n_v:=-2-\intT{v,v}
\quad\textrm{and}\quad
A:=\set{(a_v, n_v)}{v\in V_T}.
$$
For each $(a, n) \in A$, 
we put
$$
V_S(a, n):=\set{u\in S\dual}{\intS{w_S, u}=a,\; \intS{u,u}=n}.
$$
Since $S$ is hyperbolic and $\intS{w_S,w_S}>0$,
the set $V_S(a, n)$ is finite,
because $\intS{\phantom{a}, \phantom{a}}$ induces 
on the affine hyperplane 
$$
\set{x\in S\tR}{\intS{x, w_S}=a} 
$$
of $S\tR$
an inhomogeneous quadratic function whose quadratic part is negative-definite.
Then the set $\LR(w, S)$ is equal to
$$
L\;\cap\;\set{u+v}{\,v\in V_T, \;u\in V_S(a_v, n_v)\,}, 
$$
where the intersection is taken in $S\dual\oplus T\dual$.
\end{proof}
The notion of $\RRR_S$-chamber is useful in the study on $\OG^+(S)$
because of the following:
\begin{proposition}\label{prop:gRRRS}
Suppose that the natural homomorphism $\OG(T)\to \OG(q_T)$ 
is surjective.
Then the action of $\OG^+(S)$ preserves $\RRR_S$.
In particular, 
for an $\RRR_S$-chamber $\cham_S$ and an isometry $\gamma\in \OG^+(S)$, 
the image $\cham_S^\gamma$ of  $\cham_S$ by 
$\gamma$ is also an $\RRR_S$-chamber.
Moreover, 
if the interior  of $\cham_S^\gamma$ has a common point with   $\cham_S$,
then $\cham_S^\gamma=\cham_S$ holds and $\gamma$ preserves $\wallrs(\cham_S)$.
\end{proposition}
\begin{proof}
By the assumption $\OG(T)\surj \OG(q_T)$,
every element $\gamma\in \OG^+(S)$ lifts to an element $\tilde{\gamma}\in \OG(L)$  that satisfies
$\tilde{\gamma}(S)=S$ and 
$\tilde{\gamma}|_S=\gamma$ (see~\cite[Proposition 1.6.1]{MR525944}).
Since $\tilde{\gamma}$ preserves $\RRR_L$ and $\gamma$ preserves $\cone_S$,
$\gamma$ preserves  $\RRR_S$.
\end{proof}
\section{A basis of the N\'eron-Severi lattice of $X$}\label{sec:NS}
Recall that $X\subset \P^3$ is the Fermat quartic surface in characteristic $3$.
From now on, we put
$$
S:=S_{X},
$$
which is an even hyperbolic lattice of rank $22$ such that  $S\dual/S\cong (\Z/3\Z)^2$.
We use the affine coordinates $w, x, y$ of $\P^3$ with $z=1$.
\par
\medskip
Note that $X$ is the Hermitian surface over $\F_9$~(see~\cite[Chapter 23]{MR1363259}).
Hence the number of lines contained in $X$ is $112$~(see~\cite[n. 32]{MR0213949}~or~\cite[Corollary 2.22]{MR1794260}).
Since the indices of these lines are important throughout this paper,
we present 
defining equations 
of these lines in~Table~\ref{table:lines}.
(Note that $\ell_i\subset X$ implies that $\ell_i$ is not contained in the plane $z=0$ at infinity.)
%
%
%
\begin{table}
{\fontsize{9pt}{0pt}\selectfont
$$
\begin{array}{rclcrcl}
\ell_{1}&:=&\lineFQ{1}&&\ell_{2}&:=&\lineFQ{2}\\
\ell_{3}&:=&\lineFQ{3}&&\ell_{4}&:=&\lineFQ{4}\\
\ell_{5}&:=&\lineFQ{5}&&\ell_{6}&:=&\lineFQ{6}\\
\ell_{7}&:=&\lineFQ{7}&&\ell_{8}&:=&\lineFQ{8}\\
\ell_{9}&:=&\lineFQ{9}&&\ell_{10}&:=&\lineFQ{10}\\
\ell_{11}&:=&\lineFQ{11}&&\ell_{12}&:=&\lineFQ{12}\\
\ell_{13}&:=&\lineFQ{13}&&\ell_{14}&:=&\lineFQ{14}\\
\ell_{15}&:=&\lineFQ{15}&&\ell_{16}&:=&\lineFQ{16}\\
\ell_{17}&:=&\lineFQ{17}&&\ell_{18}&:=&\lineFQ{18}\\
\ell_{19}&:=&\lineFQ{19}&&\ell_{20}&:=&\lineFQ{20}\\
\ell_{21}&:=&\lineFQ{21}&&\ell_{22}&:=&\lineFQ{22}\\
\ell_{23}&:=&\lineFQ{23}&&\ell_{24}&:=&\lineFQ{24}\\
\ell_{25}&:=&\lineFQ{25}&&\ell_{26}&:=&\lineFQ{26}\\
\ell_{27}&:=&\lineFQ{27}&&\ell_{28}&:=&\lineFQ{28}\\
\ell_{29}&:=&\lineFQ{29}&&\ell_{30}&:=&\lineFQ{30}\\
\ell_{31}&:=&\lineFQ{31}&&\ell_{32}&:=&\lineFQ{32}\\
\ell_{33}&:=&\lineFQ{33}&&\ell_{34}&:=&\lineFQ{34}\\
\ell_{35}&:=&\lineFQ{35}&&\ell_{36}&:=&\lineFQ{36}\\
\ell_{37}&:=&\lineFQ{37}&&\ell_{38}&:=&\lineFQ{38}\\
\ell_{39}&:=&\lineFQ{39}&&\ell_{40}&:=&\lineFQ{40}\\
\ell_{41}&:=&\lineFQ{41}&&\ell_{42}&:=&\lineFQ{42}\\
\ell_{43}&:=&\lineFQ{43}&&\ell_{44}&:=&\lineFQ{44}\\
\ell_{45}&:=&\lineFQ{45}&&\ell_{46}&:=&\lineFQ{46}\\
\ell_{47}&:=&\lineFQ{47}&&\ell_{48}&:=&\lineFQ{48}\\
\ell_{49}&:=&\lineFQ{49}&&\ell_{50}&:=&\lineFQ{50}\\
\ell_{51}&:=&\lineFQ{51}&&\ell_{52}&:=&\lineFQ{52}\\
\ell_{53}&:=&\lineFQ{53}&&\ell_{54}&:=&\lineFQ{54}\\
\ell_{55}&:=&\lineFQ{55}&&\ell_{56}&:=&\lineFQ{56}\\
\ell_{57}&:=&\lineFQ{57}&&\ell_{58}&:=&\lineFQ{58}\\
\ell_{59}&:=&\lineFQ{59}&&\ell_{60}&:=&\lineFQ{60}\\
\ell_{61}&:=&\lineFQ{61}&&\ell_{62}&:=&\lineFQ{62}\\
\ell_{63}&:=&\lineFQ{63}&&\ell_{64}&:=&\lineFQ{64}\\
\ell_{65}&:=&\lineFQ{65}&&\ell_{66}&:=&\lineFQ{66}\\
\ell_{67}&:=&\lineFQ{67}&&\ell_{68}&:=&\lineFQ{68}\\
\ell_{69}&:=&\lineFQ{69}&&\ell_{70}&:=&\lineFQ{70}\\
\ell_{71}&:=&\lineFQ{71}&&\ell_{72}&:=&\lineFQ{72}\\
\ell_{73}&:=&\lineFQ{73}&&\ell_{74}&:=&\lineFQ{74}\\
\ell_{75}&:=&\lineFQ{75}&&\ell_{76}&:=&\lineFQ{76}\\
\ell_{77}&:=&\lineFQ{77}&&\ell_{78}&:=&\lineFQ{78}\\
\ell_{79}&:=&\lineFQ{79}&&\ell_{80}&:=&\lineFQ{80}\\
\ell_{81}&:=&\lineFQ{81}&&\ell_{82}&:=&\lineFQ{82}\\
\ell_{83}&:=&\lineFQ{83}&&\ell_{84}&:=&\lineFQ{84}\\
\ell_{85}&:=&\lineFQ{85}&&\ell_{86}&:=&\lineFQ{86}\\
\ell_{87}&:=&\lineFQ{87}&&\ell_{88}&:=&\lineFQ{88}\\
\ell_{89}&:=&\lineFQ{89}&&\ell_{90}&:=&\lineFQ{90}\\
\ell_{91}&:=&\lineFQ{91}&&\ell_{92}&:=&\lineFQ{92}\\
\ell_{93}&:=&\lineFQ{93}&&\ell_{94}&:=&\lineFQ{94}\\
\ell_{95}&:=&\lineFQ{95}&&\ell_{96}&:=&\lineFQ{96}\\
\ell_{97}&:=&\lineFQ{97}&&\ell_{98}&:=&\lineFQ{98}\\
\ell_{99}&:=&\lineFQ{99}&&\ell_{100}&:=&\lineFQ{100}\\
\ell_{101}&:=&\lineFQ{101}&&\ell_{102}&:=&\lineFQ{102}\\
\ell_{103}&:=&\lineFQ{103}&&\ell_{104}&:=&\lineFQ{104}\\
\ell_{105}&:=&\lineFQ{105}&&\ell_{106}&:=&\lineFQ{106}\\
\ell_{107}&:=&\lineFQ{107}&&\ell_{108}&:=&\lineFQ{108}\\
\ell_{109}&:=&\lineFQ{109}&&\ell_{110}&:=&\lineFQ{110}\\
\ell_{111}&:=&\lineFQ{111}&&\ell_{112}&:=&\lineFQ{112}\\
\end{array}
$$
}
\caption{Lines on $X$}\label{table:lines}
\end{table}
From these $112$ lines, we choose the following:
\begin{eqnarray}\label{eq:basisiis}
&&\ell_{1}, \ell_{2}, \ell_{3}, \ell_{4}, \ell_{5}, \ell_{6}, \ell_{7}, \ell_{9}, \ell_{10}, \ell_{11}, \ell_{17}, \\
&& \phantom{aaaaaaaa}\ell_{18}, \ell_{19},
\ell_{21}, \ell_{22}, \ell_{23}, \ell_{25}, \ell_{26}, \ell_{27}, \ell_{33}, \ell_{35}, \ell_{49}. \nonumber
\end{eqnarray}
The intersection matrix $N$ of these $22$ lines is given in Table~\ref{table:theN}. 
\begin{table}
{\small
$$
\left[
\begin{array}{cccccccccccccccccccccc}
 -2 \shrink & 1 \shrink & 1 \shrink & 1 \shrink & 1 \shrink & 0 \shrink & 0 \shrink & 1 \shrink & 0 \shrink & 0 \shrink & 1 \shrink & 0 \shrink & 0 \shrink & 0 \shrink & 0 \shrink & 0 \shrink & 0 \shrink & 0 \shrink & 1 \shrink & 1 \shrink & 0 \shrink & 1 \\ 
\noalign{\myvskipmat} 1 \shrink & -2 \shrink & 1 \shrink & 1 \shrink & 0 \shrink & 1 \shrink & 0 \shrink & 0 \shrink & 1 \shrink & 0 \shrink & 0 \shrink & 0 \shrink & 1 \shrink & 1 \shrink & 0 \shrink & 0 \shrink & 0 \shrink & 1 \shrink & 0 \shrink & 0 \shrink & 1 \shrink & 0 \\ 
\noalign{\myvskipmat} 1 \shrink & 1 \shrink & -2 \shrink & 1 \shrink & 0 \shrink & 0 \shrink & 1 \shrink & 0 \shrink & 0 \shrink & 1 \shrink & 0 \shrink & 0 \shrink & 0 \shrink & 0 \shrink & 1 \shrink & 0 \shrink & 1 \shrink & 0 \shrink & 0 \shrink & 0 \shrink & 0 \shrink & 0 \\ 
\noalign{\myvskipmat} 1 \shrink & 1 \shrink & 1 \shrink & -2 \shrink & 0 \shrink & 0 \shrink & 0 \shrink & 0 \shrink & 0 \shrink & 0 \shrink & 0 \shrink & 1 \shrink & 0 \shrink & 0 \shrink & 0 \shrink & 1 \shrink & 0 \shrink & 0 \shrink & 0 \shrink & 0 \shrink & 0 \shrink & 0 \\ 
\noalign{\myvskipmat} 1 \shrink & 0 \shrink & 0 \shrink & 0 \shrink & -2 \shrink & 1 \shrink & 1 \shrink & 1 \shrink & 0 \shrink & 0 \shrink & 0 \shrink & 1 \shrink & 0 \shrink & 0 \shrink & 0 \shrink & 1 \shrink & 1 \shrink & 0 \shrink & 0 \shrink & 0 \shrink & 0 \shrink & 0 \\ 
\noalign{\myvskipmat} 0 \shrink & 1 \shrink & 0 \shrink & 0 \shrink & 1 \shrink & -2 \shrink & 1 \shrink & 0 \shrink & 1 \shrink & 0 \shrink & 0 \shrink & 0 \shrink & 0 \shrink & 0 \shrink & 1 \shrink & 0 \shrink & 0 \shrink & 0 \shrink & 1 \shrink & 0 \shrink & 0 \shrink & 0 \\ 
\noalign{\myvskipmat} 0 \shrink & 0 \shrink & 1 \shrink & 0 \shrink & 1 \shrink & 1 \shrink & -2 \shrink & 0 \shrink & 0 \shrink & 1 \shrink & 0 \shrink & 0 \shrink & 1 \shrink & 1 \shrink & 0 \shrink & 0 \shrink & 0 \shrink & 0 \shrink & 0 \shrink & 0 \shrink & 1 \shrink & 0 \\ 
\noalign{\myvskipmat} 1 \shrink & 0 \shrink & 0 \shrink & 0 \shrink & 1 \shrink & 0 \shrink & 0 \shrink & -2 \shrink & 1 \shrink & 1 \shrink & 0 \shrink & 0 \shrink & 1 \shrink & 0 \shrink & 1 \shrink & 0 \shrink & 0 \shrink & 0 \shrink & 0 \shrink & 0 \shrink & 1 \shrink & 0 \\ 
\noalign{\myvskipmat} 0 \shrink & 1 \shrink & 0 \shrink & 0 \shrink & 0 \shrink & 1 \shrink & 0 \shrink & 1 \shrink & -2 \shrink & 1 \shrink & 0 \shrink & 1 \shrink & 0 \shrink & 0 \shrink & 0 \shrink & 0 \shrink & 1 \shrink & 0 \shrink & 0 \shrink & 0 \shrink & 0 \shrink & 0 \\ 
\noalign{\myvskipmat} 0 \shrink & 0 \shrink & 1 \shrink & 0 \shrink & 0 \shrink & 0 \shrink & 1 \shrink & 1 \shrink & 1 \shrink & -2 \shrink & 1 \shrink & 0 \shrink & 0 \shrink & 0 \shrink & 0 \shrink & 1 \shrink & 0 \shrink & 1 \shrink & 0 \shrink & 1 \shrink & 0 \shrink & 1 \\ 
\noalign{\myvskipmat} 1 \shrink & 0 \shrink & 0 \shrink & 0 \shrink & 0 \shrink & 0 \shrink & 0 \shrink & 0 \shrink & 0 \shrink & 1 \shrink & -2 \shrink & 1 \shrink & 1 \shrink & 1 \shrink & 0 \shrink & 0 \shrink & 1 \shrink & 0 \shrink & 0 \shrink & 0 \shrink & 0 \shrink & 0 \\ 
\noalign{\myvskipmat} 0 \shrink & 0 \shrink & 0 \shrink & 1 \shrink & 1 \shrink & 0 \shrink & 0 \shrink & 0 \shrink & 1 \shrink & 0 \shrink & 1 \shrink & -2 \shrink & 1 \shrink & 0 \shrink & 1 \shrink & 0 \shrink & 0 \shrink & 1 \shrink & 0 \shrink & 1 \shrink & 0 \shrink & 0 \\ 
\noalign{\myvskipmat} 0 \shrink & 1 \shrink & 0 \shrink & 0 \shrink & 0 \shrink & 0 \shrink & 1 \shrink & 1 \shrink & 0 \shrink & 0 \shrink & 1 \shrink & 1 \shrink & -2 \shrink & 0 \shrink & 0 \shrink & 1 \shrink & 0 \shrink & 0 \shrink & 1 \shrink & 0 \shrink & 0 \shrink & 0 \\ 
\noalign{\myvskipmat} 0 \shrink & 1 \shrink & 0 \shrink & 0 \shrink & 0 \shrink & 0 \shrink & 1 \shrink & 0 \shrink & 0 \shrink & 0 \shrink & 1 \shrink & 0 \shrink & 0 \shrink & -2 \shrink & 1 \shrink & 1 \shrink & 1 \shrink & 0 \shrink & 0 \shrink & 1 \shrink & 0 \shrink & 0 \\ 
\noalign{\myvskipmat} 0 \shrink & 0 \shrink & 1 \shrink & 0 \shrink & 0 \shrink & 1 \shrink & 0 \shrink & 1 \shrink & 0 \shrink & 0 \shrink & 0 \shrink & 1 \shrink & 0 \shrink & 1 \shrink & -2 \shrink & 1 \shrink & 0 \shrink & 1 \shrink & 0 \shrink & 0 \shrink & 0 \shrink & 1 \\ 
\noalign{\myvskipmat} 0 \shrink & 0 \shrink & 0 \shrink & 1 \shrink & 1 \shrink & 0 \shrink & 0 \shrink & 0 \shrink & 0 \shrink & 1 \shrink & 0 \shrink & 0 \shrink & 1 \shrink & 1 \shrink & 1 \shrink & -2 \shrink & 0 \shrink & 0 \shrink & 1 \shrink & 0 \shrink & 1 \shrink & 0 \\ 
\noalign{\myvskipmat} 0 \shrink & 0 \shrink & 1 \shrink & 0 \shrink & 1 \shrink & 0 \shrink & 0 \shrink & 0 \shrink & 1 \shrink & 0 \shrink & 1 \shrink & 0 \shrink & 0 \shrink & 1 \shrink & 0 \shrink & 0 \shrink & -2 \shrink & 1 \shrink & 1 \shrink & 0 \shrink & 1 \shrink & 0 \\ 
\noalign{\myvskipmat} 0 \shrink & 1 \shrink & 0 \shrink & 0 \shrink & 0 \shrink & 0 \shrink & 0 \shrink & 0 \shrink & 0 \shrink & 1 \shrink & 0 \shrink & 1 \shrink & 0 \shrink & 0 \shrink & 1 \shrink & 0 \shrink & 1 \shrink & -2 \shrink & 1 \shrink & 0 \shrink & 0 \shrink & 0 \\ 
\noalign{\myvskipmat} 1 \shrink & 0 \shrink & 0 \shrink & 0 \shrink & 0 \shrink & 1 \shrink & 0 \shrink & 0 \shrink & 0 \shrink & 0 \shrink & 0 \shrink & 0 \shrink & 1 \shrink & 0 \shrink & 0 \shrink & 1 \shrink & 1 \shrink & 1 \shrink & -2 \shrink & 0 \shrink & 0 \shrink & 1 \\ 
\noalign{\myvskipmat} 1 \shrink & 0 \shrink & 0 \shrink & 0 \shrink & 0 \shrink & 0 \shrink & 0 \shrink & 0 \shrink & 0 \shrink & 1 \shrink & 0 \shrink & 1 \shrink & 0 \shrink & 1 \shrink & 0 \shrink & 0 \shrink & 0 \shrink & 0 \shrink & 0 \shrink & -2 \shrink & 1 \shrink & 0 \\ 
\noalign{\myvskipmat} 0 \shrink & 1 \shrink & 0 \shrink & 0 \shrink & 0 \shrink & 0 \shrink & 1 \shrink & 1 \shrink & 0 \shrink & 0 \shrink & 0 \shrink & 0 \shrink & 0 \shrink & 0 \shrink & 0 \shrink & 1 \shrink & 1 \shrink & 0 \shrink & 0 \shrink & 1 \shrink & -2 \shrink & 1 \\ 
\noalign{\myvskipmat} 1 \shrink & 0 \shrink & 0 \shrink & 0 \shrink & 0 \shrink & 0 \shrink & 0 \shrink & 0 \shrink & 0 \shrink & 1 \shrink & 0 \shrink & 0 \shrink & 0 \shrink & 0 \shrink & 1 \shrink & 0 \shrink & 0 \shrink & 0 \shrink & 1 \shrink & 0 \shrink & 1 \shrink & -2 
\end{array}
\right]
$$
}
\caption{Gram matrix $N$ of $S$}\label{table:theN}
\end{table}
Since $\det N=-9$,  the classes $[\ell_i]\in S$ of the  lines $\ell_i$ in~\eqref{eq:basisiis} form a basis of $S$.
Throughout this paper, we fix this basis,  and write elements of $S\tR$ as row vectors
$$
[x_1, \dots, x_{22}]_S.
$$
When we use its dual basis, we write  
$$
[\xi_1, \dots, \xi_{22}]_S\dual.
$$
Since the hyperplane $w+(1+i)=0$ cuts out from $X$ the divisor $\ell_1+\ell_2+\ell_3+\ell_4$,
the class $h_0=[\OOO_{X}(1)]\in S$ of the hyperplane section is equal to
\begin{eqnarray*}
h_0 &=& [1,1,1,1,0,0,0,0,0,0,0,0,0,0,0,0,0,0,0,0,0,0]_S \\
& =& [1,1,1,1,1,1,1,1,1,1,1,1,1,1,1,1,1,1,1,1,1,1]_S\dual.
\end{eqnarray*}
As a positive cone $\cone_S$ of $S$,
we choose the connected component containing $h_0$.
\par
\medskip
From the intersection numbers of the  $112$ lines,
we can calculate their  classes $[\ell_i]\in S$. 
\begin{remark}\label{rem:F9}
Since these $112$ lines are all defined over $\F_9$,
every class $v\in S$ is represented by a divisor defined over $\F_9$.
More generally, 
Sch\"utt~\cite{schuttpreprint} showed that 
a supersingular $K3$ surface with Artin invariant $1$
in  characteristic $p$
has a projective model defined over $\F_{p^2}$, and
its N\'eron-Severi lattice 
 is generated by the classes of divisors defined over $\F_{p^2}$. 
\end{remark}
\begin{proposition}\label{prop:28}
We have
$$
h_0=\frac{1}{28}\,\sum_{i=1}^{112} [\ell_i].
$$
\end{proposition}
\begin{proof}
The number of $\F_9$-rational points on $X$ is $280$.
For each $\F_9$-rational point $P$ of $X$,
the tangent plane $T_{X, P}\subset \P^3$ to $X$ at $P$ 
cuts out a union of four lines from $X$.
Since each line contains ten $\F_9$-rational points,
we have $280\, h_0=10\, \sum [\ell_i]$.
\end{proof}
As before, we let $\OG(S)$ act on $S$ from \emph{right}, so that 
$$
\OG(S)=\set{T\in \GL_{22}(\Z)}{T\,N\,{}^t\hskip -.3pt T=N}.
$$
We also let the projective automorphism group $\Aut(X, h_0)=\PGU_4(\F_9)$ act on $X$ from \emph{right}.
For each $\tau \in \PGU_4(\F_9)$, we can calculate its action $\tau_*$ on $S$ by looking at the permutation
of the $112$ lines induced by $\tau$.
\begin{example}
Consider the projective automorphism 
$$
\tau \;\; :\;\;  [w:x:y:z]\mapsto [w:x:y:z]  
\left[ \begin {array}{cccc} 
1&-1&-i&-1\\
\noalign{\myvskipmat}
-1&-1+\sqrmo&-1+\sqrmo&-1+\sqrmo\\
\noalign{\myvskipmat}
-1 & -\sqrmo & -1-\sqrmo & 0\\
\noalign{\myvskipmat}
\sqrmo&0&-1+\sqrmo&-1\end {array} \right] 
$$
of $X$.
\footnote{
This matrix was incorrect in the previous version of this paper. We are grateful to Benjamin Church for bringing this error to our attention. 
This error was corrected in August 2026.}
Then the images $\ell_i^\tau$ of the lines $\ell_i$ in~\eqref{eq:basisiis} are
$$
\begin{array}{cl}
&\ell_{1}^\tau=\ell_{60}, \;\; 
\ell_{2}^\tau=\ell_{31}, \;\; 
\ell_{3}^\tau=\ell_{105}, \;\; 
\ell_{4}^\tau=\ell_{95}, \;\; 
\ell_{5}^\tau=\ell_{92}, \;\; 
\ell_{6}^\tau=\ell_{30}, \;\; 
\\
 &\ell_{7}^\tau=\ell_{76}, \;\; 
\ell_{9}^\tau=\ell_{110}, \;\; 
\ell_{10}^\tau=\ell_{29}, \;\; 
\ell_{11}^\tau=\ell_{6}, \;\; 
\ell_{17}^\tau=\ell_{20}, \;\; 
\ell_{18}^\tau=\ell_{96}, \;\; 
\\
 &\ell_{19}^\tau=\ell_{102}, \;\; 
\ell_{21}^\tau=\ell_{13}, \;\; 
\ell_{22}^\tau=\ell_{87}, \;\; 
\ell_{23}^\tau=\ell_{91}, \;\; 
\ell_{25}^\tau=\ell_{108}, \;\; 
\ell_{26}^\tau=\ell_{10}, \;\; 
\\
 &\ell_{27}^\tau=\ell_{57}, \;\; 
\ell_{33}^\tau=\ell_{52}, \;\; 
\ell_{35}^\tau=\ell_{51}, \;\; 
\ell_{49}^\tau=\ell_{59}\;\; 
 \hskip -5pt.
\end{array}
$$
Therefore the action  $\tau_*$ on $S$ is given by $v\mapsto v T_{\tau}$,
where $T_{\tau}$ is the matrix whose row vectors are
$$
\renewcommand\arraystretch{1.13} 
\begin{array}{ccl}
\numclass{{\ell_{60}}}&=&\thevFQ{60},    \\
\numclass{{\ell_{31}}}&=&\thevFQ{31},    \\
\numclass{{\ell_{105}}}&=&\thevFQ{105},    \\
\numclass{{\ell_{95}}}&=&\thevFQ{95},    \\
\numclass{{\ell_{92}}}&=&\thevFQ{92},    \\
\numclass{{\ell_{30}}}&=&\thevFQ{30},    \\
\numclass{{\ell_{76}}}&=&\thevFQ{76},    \\
\numclass{{\ell_{110}}}&=&\thevFQ{110},    \\
\numclass{{\ell_{29}}}&=&\thevFQ{29},    \\
\numclass{{\ell_{6}}}&=&\thevFQ{6},    \\
\numclass{{\ell_{20}}}&=&\thevFQ{20},    \\
\numclass{{\ell_{96}}}&=&\thevFQ{96},    \\
\numclass{{\ell_{102}}}&=&\thevFQ{102},    \\
\numclass{{\ell_{13}}}&=&\thevFQ{13},    \\
\numclass{{\ell_{87}}}&=&\thevFQ{87},    \\
\numclass{{\ell_{91}}}&=&\thevFQ{91},    \\
\numclass{{\ell_{108}}}&=&\thevFQ{108},    \\
\numclass{{\ell_{10}}}&=&\thevFQ{10},    \\
\numclass{{\ell_{57}}}&=&\thevFQ{57},    \\
\numclass{{\ell_{52}}}&=&\thevFQ{52},    \\
\numclass{{\ell_{51}}}&=&\thevFQ{51},    \\
\numclass{{\ell_{59}}}&=&\thevFQ{59}.    
\end{array}
$$
\end{example}
We put the representation 
\begin{equation}\label{eq:Ttau}
\tau \mapsto T_{\tau}
\end{equation}
of $\Aut(X, h_0)=\PGU_4(\F_9)$ to  $\OG^+(S)$  in the computer memory. 
It turns out to be faithful.
On the other hand, $\Aut(X, h_0)$ is just the stabilizer subgroup in $\Aut(X)$ of $h_0\in S$.
Therefore we confirm the following fact~(\cite[Section 8, Proposition 3]{MR633161}):
\begin{proposition}\label{prop:autoembed}
The action of  $\Aut(X)$ on $S$ is faithful.
\end{proposition}
From now on,
we regard $\Aut(X)$ as a subgroup of $\OG^+(S)$,
and write $v\mapsto v^\gamma$ instead of $v\mapsto v^{\gamma_*}$
for the action $\gamma_*$ of $\gamma\in \Aut(X)$ on $S$.
\section{Embedding of $S$ into $L$}\label{sec:embS}
Next we embed the N\'eron-Severi lattice $S$ of $X$ into the even unimodular hyperbolic lattice  
of rank $26$,
and calculate the walls of an  $\RRR_S$-chamber.
\par
\medskip
Let $T$ be the negative-definite root lattice of type $2A_2$.
We fix a basis of $T$ in such a way that the Gram matrix is equal to
$$
\left[ \begin {array}{cccc} -2&1&0&0\\\noalign{\myvskipmat}1&-2&0&0
\\\noalign{\myvskipmat}0&0&-2&1\\\noalign{\myvskipmat}0&0&1&-2\end {array}
 \right].
$$
When we use this basis, 
we write   elements of $T\tR$ as 
$[y_1, y_2, y_3, y_4]_T$, 
while  when we use its dual basis,
we write 
as $[\eta_1, \eta_2, \eta_3, \eta_4]_T\dual$.
Elements of $(S\oplus T)\tR$ are written as 
$$
[x_1, \dots, x_{22}\;|\; y_1, \dots, y_4]
$$
using the bases of $S$ and $T$, 
or as
$$
[\xi_1, \dots, \xi_{22}\;|\; \eta_1, \dots, \eta_4]\dual
$$
using the dual bases of $S\dual$ and $T\dual$. 
\par
\medskip
Consider the following vectors of $S\dual\oplus T\dual$:
\begin{eqnarray*}
a_1 &:=&\frac{1}{3}\,[2, 2, 0, 0, 0, 1, 2, 2, 1, 1, 2, 2, 1, 1, 2, 0, 0, 1, 1, 0, 0, 0 \;|\; 1, 2, 0, 0], \\
a_2 &:=&\frac{1}{3}\,[2, 0, 2, 0, 2, 1, 1, 0, 2, 1, 2, 1, 0, 2, 2, 1, 1, 0, 1, 0, 0, 0 \;|\; 0, 0, 1, 2].
\end{eqnarray*}
We  define   $\alpha_1, \alpha_2 \in (S\oplus T)\dual/(S\oplus T)$ by 
$$
\alpha_1:=a_1\bmod (S\oplus T),\quad \alpha_2:=a_2\bmod (S\oplus T).
$$
Then $\alpha_1$ and $\alpha_2$ are linearly independent in $(S\oplus T)\dual/(S\oplus T)\cong \F_3^4$. 
Since 
$$
q_{S\oplus T}(\alpha_1)=q_{S\oplus T}(\alpha_2)=q_{S\oplus T}(\alpha_1+\alpha_2)=0,
$$
the vectors $\alpha_1$ and $\alpha_2$ 
generate a maximal isotropic subgroup of $q_{S\oplus T}$.
Therefore, by~\cite[Proposition 1.4.1]{MR525944}, 
the submodule
$$
L:=(S\oplus T)+\gen{a_1}+\gen{a_2}
$$
of $S\dual\oplus T\dual$ is an even unimodular overlattice of $S\oplus T$ into which $S$ and $T$ are primitively embedded.

By construction, $L$ is hyperbolic of rank $26$.
We choose  $\cone_L$ to be the connected component that contains $\cone_S$.
Then, 
by means of the roots of $L$, 
we obtain a  decomposition of $\cone_S$ into $\RRR_S$-chambers.
\par
\medskip
The order of $\OG(T)$ is $288$,
while the order of $\OG(q_T)$ is $8$.
It is easy to check that the natural homomorphism
$\OG(T)\to \OG(q_T)$ is surjective.
Therefore we obtain the following from Proposition~\ref{prop:gRRRS}:
\begin{proposition}
The action of $\OG^+(S)$ on $S\tR$ preserves $\RRR_S$.
\end{proposition}
We put
\begin{eqnarray*}
w_0&:=&[1,1,1,1,0,0,0,0,0,0,0,0,0,0,0,0,0,0,0,0,0,0 \;|\; -1,-1,-1,-1]\\
&=&[1,1,1,1,1,1,1,1,1,1,1,1,1,1,1,1,1,1,1,1,1,1 \;|\; 1,1,1,1]\dual.
\end{eqnarray*}
Note that the projection $w_{0S}\in S\dual$ of $w_0$ to $S\dual$ is 
equal to $h_0$.

Since $\intL{w_0, w_0}=0$ and $\intL{w_0, h_0}>0$, 
we see that $w_0$ is on the boundary of the closure of $\cone_L$ in $L\tR$.
\begin{proposition}\label{prop:Weyl0}
The vector $w_0$ is a Weyl vector,
and the $\RRR_L$-chamber $\cham_L(w_0)$ is $S$-nondegenerate.
The $\RRR_S$-chamber 
$$
\chamsz:=\cham_L(w_0) \cap (S\tR)
$$
contains 
$w_{S0}=h_0$ in its interior.
\end{proposition}
\begin{proof}
The only non-trivial part of  the first assertion is that
$\gen{w_0}\sperp/\gen{w_0}$ has no vectors of square norm $-2$.
We put 
$$
w_0\sprime:=[7, 6, 7, 7, 7, 7, 7, 7, 7, 7, 7, 7, 7, 7, 7, 7, 7, 7, 7, 7, 7, 7
\;|\; 7, 5, 7, 7]\dual.
$$
Then we have $\intL{w_0\sprime, w_0\sprime}=0$ and $\intL{w_0, w_0\sprime}=1$.
Let $U\subset L$ be the sublattice generated by $w_0$ and $ w_0\sprime$.
Calculating a basis $\lambda_1, \dots, \lambda_{24}$ of $U\sperp\subset L$,
we obtain a Gram matrix of $U\sperp$,
which is negative-definite of determinant $1$.
By the algorithm described in~\cite[Section 3.1]{shimadapreprintchar5}, 
we verify that there are no vectors of square norm $-2$ in $U\sperp$.
\par
We show that $w_0$ satisfies  the conditions (i) and (ii) 
 given after Definition~\ref{def:Snondeg}.
 By Proposition~\ref{prop:lambda}, 
 in order to verify the condition (i),
it is enough  to show that
the function $Q: U\sperp\to \Z$ given by
$$
Q(\lambda):=\intL{h_0, -\frac{2+\intL{\lambda, \lambda}}{2}w_0+w_0\sprime+\lambda}
$$
does not take  negative values.
Using the basis $\lambda_1, \dots, \lambda_{24}$ of $U\sperp$,
we can write $Q$ as an inhomogeneous quadratic function of $24$ variables.
Its  quadratic part turns out to be positive-definite.
By the algorithm described in~\cite[Section 3.1]{shimadapreprintchar5}, 
we verify that there exist no vectors $\lambda\in U\sperp$ such that $Q(\lambda)<0$.
Next we show that $w_{0S}=h_0\in \cone_S$ 
has the property  required for $v\sprime$  in the condition (ii),
and hence $h_0$ is contained in the interior of  $\chamsz$.
Note that $w_{0T}=[-1,-1,-1,-1]_T$ is non-zero.
Hence 
we can calculate 
$$
\LR(w_0, S)=\set{r\in \LR(w_0)}{\intS{r_S, r_S}<0}
$$
by the method described in the proof of Proposition~\ref{prop:LRS}. 
Then we can easily show  that $h_0$ satisfies
$\intL{h_0, r}>0$ for any $r\in \LR(w_0, S)$.
\end{proof}
\begin{remark}
There exist exactly four vectors $\lambda\in U\sperp$
such that $Q(\lambda)=0$.
They correspond to the Leech roots $r\in \LR(w_0)$
such that $r=r_T$.
\end{remark}
From the surjectivity of $\OG(T)\to \OG(q_T)$
and Proposition~\ref{prop:gRRRS},
we obtain the following:
\begin{corollary}\label{cor:action_on_CS0}
The action of $\Aut(X, h_0)$ on $S\tR$ preserves $\chamsz$ and $\wallrs(\chamsz)$.
\end{corollary}
%
%
%
\begin{proposition}\label{prop:CS0}
The maps $r\mapsto r_S$ and $r_S\mapsto (r_S)\sperp_S$
induce bijections
$$
\LR(w_0, S)\;\;\cong\;\; \wallrs(\chamsz)\;\;\cong\;\; \wall(\chamsz).
$$
The action of $\Aut(X,  h_0)$ 
 decomposes $\wallrs(\chamsz)$
into the  three orbits 
$$
\wallrsz{112}:=\wallrs(\chamsz)_{[1, -2]},
\quad
\wallrsz{648}:=\wallrs(\chamsz)_{[2, -4/3]}
\quad\textrm{and}\quad
\wallrsz{5184}:=\wallrs(\chamsz)_{[3, -2/3]}
$$
of cardinalities $112$, $648$ and $5184$, respectively,
where
$$
\wallrs(\chamsz)_{[a, n]}:=\set{r_S\in \wallrs(\chamsz)}{\intS{r_S, h_0}=a, \; \intS{r_S,r_S}=n}.
$$
The set $\wallrsz{112}$ coincides 
with the set of the classes $[\ell_i]$ of lines contained in $X$:
$$
\wallrsz{112}=\{[\ell_1], [\ell_2], \dots, [\ell_{112}]\}.
$$
The sets $\wallrsz{648}$ and  $\wallrsz{5184}$ are the orbits of
\begin{eqnarray*}
b_1&:=&\frac{1}{3}\,[-1, 0, -1, 0, 2, 1, 1, 0, 2, 1, -1, 1, 0, -1, -1, 1, 1, 0, 1, 0, 0, 0]_S
\;\;\in\;\; \wallrsz{648}, \;\;\textrm{and} \\
b_2&:=&\frac{1}{3}\,[0, 1, -1, 0, 2, 0, 2, 1, 1, 0, 0, -1, 2, 1, 0, 1, 1, -1, 0, 0, 0, 0]_S
\;\;\in\;\; \wallrsz{5184},
\end{eqnarray*}
by the action of $\Aut(X,  h_0)$, respectively.
\end{proposition}
\begin{proof}
We have calculated the finite set $\LR(w_0, S)$
in the proof of Proposition~\ref{prop:Weyl0}.
We have also stored the classes $[\ell_i]$ of the $112$ lines
and  the action of $\Aut(X, h_0)$ on $S$ in the computer memory.
Thus the assertions of Proposition~\ref{prop:CS0} are 
verified by  a direct computation,
except for the fact that,
for any $r\in \LR(w_0, S)$,
the hyperplane $(r_S)\sperp_S$ actually bounds $\chamsz$.
This is proved by showing that the point
$$
p:=h_0-\frac{\intS{h_0, r_S}}{\intS{r_S, r_S}} r_S
$$
on $(r_S)\sperp_S$ satisfies $\intL{p, r\sprime}>0$
for any $r\sprime\in \LR(w_0, S)\setminus \{r\}$.
\end{proof}
Since Proposition~\ref{prop:28} implies that the interior point $h_0$ of $\chamsz$ is determined 
by  $\wallrsz{112}$ and since $\OG(T)\to \OG(q_T)$ is surjective,
we obtain the following from Proposition~\ref{prop:gRRRS}:
\begin{corollary}\label{cor:equiv}
For $\gamma\in \OG^+(S)$, the following are equivalent:
{\rm (i)} the interior of $\chamsz^\gamma$ has a common point with $\chamsz$,   
{\rm (ii)} $\chamsz^\gamma=\chamsz$,  
{\rm (iii)} ${\wallrsz{112}}^{\,\gamma}=\wallrsz{112}$, 
{\rm (iv)} $h_0^\gamma=h_0$, and 
{\rm (v)} $h_0^\gamma\in \chamsz$.
\end{corollary}
In particular, we obtain the following:
\begin{corollary}\label{cor:converse}
If $\gamma\in \Aut(X)$ satisfies $h_0^\gamma\in \chamsz$,
then $\gamma$ is in $\Aut(X, h_0)$.
\end{corollary}
\section{The automorphisms  $g_1$ and $g_2$}\label{sec:A}
In order to find  automorphisms $\gamma\in \Aut(X)$ such that 
$h_0^\gamma\notin \chamsz$,
we search for polarizations of degree $2$ that are located on the walls
$(b_1)\sperp_S$  and $(b_2)\sperp_S$.
\par
\medskip 
We fix terminologies and notation.
For a vector $v\in S$, we denote by $\LLL_v\to X$ a  line bundle defined over $\F_9$ whose class is $v$
 (see Remark~\ref{rem:F9}).
We say that a vector $h\in S$ is a \emph{polarization} of degree $d$
if $\intS{h, h}=d$ and the complete linear system $|\LLL_h|$ is nonempty and has no fixed components.
If $h$ is a polarization,
then $|\LLL_h|$ has no base-points by~\cite[Corollary 3.2]{MR0364263} and hence defines a morphism
$$
\Phi_h : X\to \P^N,
$$
where $N=\dim |\LLL_h|$.
\par
A polynomial in $\F_9[w,x,y]$ is said to be \emph{of normal form}
if its degree with respect to $w$ is $\le 3$.
For each polynomial $G\in \F_9[w,x,y]$,
there exists a unique polynomial $\overline{G}$ of normal form such that
$$
G\equiv \overline{G} \mod (w^4+x^4+y^4+1).
$$
We say that $\overline{G}$ is the \emph{normal form} of $G$.
For any $d\in \Z$,
the vector space 
$H^0(X, \LLL_{dh_0})$ over $\F_9$ is naturally identified with
the vector subspace
$$
\Gamma(d):=\set{G\in \F_9[w,x,y]}{\textrm{$G$ is of  normal form with total degree $\le d$}}
$$
of $\F_9[w,x,y]$.
For an ideal $J$ of $\F_9[w,x,y]$,
we put
$$
\Gamma(d, J):=\Gamma(d)\cap J.
$$
A basis of $\Gamma(d, J)$ is easily obtained  by a Gr\"obner basis of $J$.
Let $\ell_i$ be a line contained in $X$.
We denote by $I_i\subset\F_9[w,x,y]$ the affine defining ideal of $\ell_i$ in $\P^3$
(see~Table~\ref{table:lines}), and 
 put
$$
I_i^{(\nu)}:=I_i^\nu +(w^4+x^4+y^4+1)\;\;\subset\;\; \F_9[w,x,y]
$$
for nonnegative integers $\nu$.
Suppose that $v\in S$ is written as
\begin{equation}\label{eq:d}
v=d\, h_0 -\sum_{i=1}^{112} a_i [\ell_i],
\end{equation}
where $a_i$ are nonnegative integers.
Then there exists a natural isomorphism
$$
H^0(X, \LLL_v)\cong \Gamma(d, \bigcap_{i=1}^{112} I_i^{(a_i)})
$$
with the property  that,
for 
$$
v\sprime =d\sprime h_0 -\sum_{i=0}^{112} a\sprime_i [\ell_i]\;\;\in\;\; S
$$
with $a_i\sprime\in \Z_{\ge 0}$,
the  multiplication homomorphism
$$
H^0(X, \LLL_v)\times H^0(X, \LLL_{v\sprime})\to H^0(X, \LLL_{v+v\sprime})
$$
is identified with
$$
\Gamma(d, \bigcap I_i^{(a_i)})\times
\Gamma(d\sprime, \bigcap I_i^{(a_i\sprime)})
\to
\Gamma(d+d\sprime, \bigcap I_i^{(a_i+a_i\sprime)})
$$
given by $(\overline{G}, \overline{G\sprime})\mapsto \overline{GG\sprime}$.
\par
\medskip
Proposition~\ref{prop:phi}  in Introduction is an immediate consequence  of the following:
\begin{proposition}\label{prop:ms}
Consider the vectors
\begin{eqnarray*}
m_1 &:=&[-1, 0, -1, -1, 2, 2, 1, 1, 2, 0, -1, 1, 1, -1, 0, 1, 0, 0, 1, -1, 0, 0]_S \quand \\
m_2 &:=&[2, 2, 1, 2, 1, -1, 1, 1, 1, 1, 0, -1, 0, 0, 0, 0, 0, -1, -1, 1, 0, 1]_S
\end{eqnarray*}
of $S$.
Then each $m_i$ is a polarization of degree $2$.
If we choose a basis of  the vector space $H^0(X, \LLL_{m_i})$ appropriately,
the morphism $\Phi_{m_i}: X\to \P^2$ 
associated with $|\LLL_{m_i}|$ coincides with the morphism
$\phi_i: X\to \P^2$ given in the statement of Proposition~\ref{prop:phi}.
\end{proposition}
\begin{proof}
We have $\intS{m_i, m_i}=2$.
By the method described in~\cite[Section 4.1]{shimadapreprintchar5},
we see that $m_i$ is a polarization;
namely, we verify
that the sets
\begin{eqnarray*}
&&\set{v\in S}{\intS{v, v}=-2, \;  \intS{v, m_i}<0,\;  \intS{v, h_0}>0} \quand \\
&&\set{v\in S}{\intS{v, v}=0, \; \intS{v, m_i}=1}
\end{eqnarray*}
are both empty.
Since
\begin{eqnarray}
m_1&=&3h_0-([\ell_{21}]+[\ell_{22}]+[\ell_{50}]+[\ell_{63}]+[\ell_{65}]+[\ell_{88}]) \quand 
\label{eq:hhm1}\\
m_2&=&5h_0-([\ell_{1}]+[\ell_{3}]+[\ell_{6}]+[\ell_{18}]+ \label{eq:hhm2}\\
&&\hskip 2cm +[\ell_{35}]+[\ell_{74}]+[\ell_{90}]+[\ell_{92}]+[\ell_{110}]+[\ell_{111}]),
 \nonumber
\end{eqnarray}
the vector spaces $H^0(X, \LLL_{m_1})$ and $H^0(X, \LLL_{m_2})$ are identified with
the  subspaces 
\begin{eqnarray*}
\Gamma_1 &:=& \Gamma(3, I_{21}\cap I_{22}\cap I_{50}\cap I_{63}\cap I_{65}\cap I_{88})\quand
\\
\Gamma_2 &:=& 
\Gamma(5, I_{1}\cap I_{3}\cap I_{6}\cap I_{18}\cap I_{35}\cap I_{74}\cap I_{90}\cap I_{92}\cap I_{110}\cap I_{111})
\end{eqnarray*}
of $\F_9[w,x,y]$,
respectively.
We calculate   
a basis of  $\Gamma_i$
by means of Gr\"obner bases of the ideals $I_i$.
The set $\{F_{i0}, F_{i1}, F_{i2}\}$ of polynomials  in Table~\ref{table:Fs} is just 
a basis of $\Gamma_i$ thus calculated.
\end{proof}
\begin{remark}
The polarizations  $m_1$ and $m_2$ in Proposition~\ref{prop:ms}
are located on the hyperplanes $(b_1)\sperp_S$ and $(b_2)\sperp_S$
bounding $\chamsz$, respectively,
where $b_1\in \wallrsz{648}$ and $b_2\in \wallrsz{5184}$
are given in Proposition~\ref{prop:CS0}. 
\end{remark}
%
%
\begin{proof}[Proof of Proposition~\ref{prop:sectics}]
The set  $\Exc(\phi_i)$  of the classes of $(-2)$-curves 
contracted by $\phi_i: X\to \P^2$
is calculated by the method 
described in~\cite[Section 4.2]{shimadapreprintchar5}.
We first calculate
the set
$$
R_i^+:=\set{v\in S}{\intS{v, v}=-2, \; \intS{v, m_i}=0,\;  \intS{v, h_0}>0}.
$$
It turns out that every element of $R_i^+$ is written as a linear combination
with coefficients in $\Z_{\ge 0}$  of elements
$l\in R_i^+$ such that $\intS{l, h_0}=1$.
Hence we have
$$
\Exc(\phi_i)=\set{l\in R_i^+}{\intS{l, h_0}=1}.
$$
The $ADE$-type of the root system $\Exc(\phi_i)$ is equal to
$6A_1+4A_2$ for $i=1$ and $A_1+A_2+2A_3+2A_4$ for $i=2$.
Thus the assertion on the $ADE$-type of the singularities of $Y_i$ is proved.
Moreover 
we have proved that all $(-2)$-curves 
contracted by $\phi_i: X\to \P^2$ are 
lines.
See Tables~\ref{table:contracts1} and~\ref{table:contracts2},
in which the lines $\ell_{k_1}, \dots, \ell_{k_r}$ contracted by $\phi_i$ to a singular point $P$ 
of type $A_r$ are indicated in such an order that 
$\intS{\ell_{k_\nu}, \ell_{k_{\nu+1}}}=1$ holds for $\nu=1, \dots, r-1$.
\begin{table}
$$
\begin{array}{rclc}
\ell_{37}&\mapsto&[1: 1-\sqrmo: 1-\sqrmo]& \textrm{($A_{1}$-point)} \\ 
\ell_{23}&\mapsto&[1: 1+\sqrmo: -(1+\sqrmo)]& \textrm{($A_{1}$-point)} \\ 
\ell_{62}&\mapsto&[1: -(1+\sqrmo): 0]& \textrm{($A_{1}$-point)} \\ 
\ell_{102}&\mapsto&[1: -(1-\sqrmo): 0]& \textrm{($A_{1}$-point)} \\ 
\ell_{68}&\mapsto&[1: 1+\sqrmo: 1+\sqrmo]& \textrm{($A_{1}$-point)} \\ 
\ell_{112}&\mapsto&[1: 1-\sqrmo: -(1-\sqrmo)]& \textrm{($A_{1}$-point)} \\ 
\ell_{49}, \ell_{29}&\mapsto&[1: 1: -\sqrmo]& \textrm{($A_{2}$-point)} \\ 
\ell_{73}, \ell_{60}&\mapsto&[1: 1: \sqrmo]& \textrm{($A_{2}$-point)} \\ 
\ell_{18}, \ell_{10}&\mapsto&[0: 1: -1]& \textrm{($A_{2}$-point)} \\ 
\ell_{16}, \ell_{99}&\mapsto&[0: 1: 1]& \textrm{($A_{2}$-point)} 
\end{array}
$$
\caption{Lines contracted by $\phi_1: X\to \P^2$}\label{table:contracts1}
\end{table}
\begin{table}
$$
\begin{array}{rclc}
\ell_{43}&\mapsto&[0: 1: 0]& \textrm{($A_{1}$-point)} \\ 
\ell_{76}, \ell_{94}&\mapsto&[1: -1: 0]& \textrm{($A_{2}$-point)} \\ 
\ell_{22}, \ell_{49}, \ell_{20}&\mapsto&[1: -1: 1]& \textrm{($A_{3}$-point)} \\ 
\ell_{7}, \ell_{5}, \ell_{103}&\mapsto&[1: -1: -1]& \textrm{($A_{3}$-point)} \\ 
\ell_{10}, \ell_{2}, \ell_{4}, \ell_{91}&\mapsto&[1: 0: 1]& \textrm{($A_{4}$-point)} \\ 
\ell_{33}, \ell_{36}, \ell_{72}, \ell_{83}&\mapsto&[1: 0: -1]& \textrm{($A_{4}$-point)} 
\end{array}
$$
\caption{Lines contracted by $\phi_2: X\to \P^2$}\label{table:contracts2}
\end{table}
\par
The defining equation $f_i=0$ of the branch curve $B_i\subset \P^2$ is calculated by the method 
given in~\cite[Section 5]{shimadapreprintchar5}.
We calculate a basis of the vector space $H^0(X, \LLL_{3m_i})$ of dimension $11$
using~\eqref{eq:hhm1},~\eqref{eq:hhm2}
and Gr\"obner bases of $I_i^{(3)}$.
Note that the ten normal forms $M_{i, 1}, \dots, M_{i, 10}$ 
of the cubic monomials of $F_{i0}, F_{i1}, F_{i2}$ are
contained in $H^0(X, \LLL_{3m_i})$.
We choose a polynomial $G_i \in H^0(X, \LLL_{3m_i})$
that is not contained in the linear span of $M_{i, 1}, \dots, M_{i, 10}$.
In the vector space $H^0(X, \LLL_{6m_i})$ of  dimension $38$,
the $39$ normal forms of the  monomials  of $G_i, F_{i0}, F_{i1}, F_{i2}$
of weighted degree $6$ 
with weight $\deg G_i=3$ and $ \deg F_{ij}=1$  have  a non-trivial linear relation.
Note that this linear relation is quadratic with respect to $G_i$.
Completing the square and 
re-choosing $G_i$ appropriately,
we confirm that
$$
\overline{G_i^2+f_i(F_{i0}, F_{i1}, F_{i2})}=0
$$
holds. 
Hence $Y_i$ is defined by $y^2+f_i(x_0, x_1, x_2)=0$.
\end{proof}
\begin{remark}\label{rem:F3}
In order to obtain a defining equation of $B_i$ with coefficients in $\F_3$,
we have to choose the basis $F_{i0}, F_{i1}, F_{i2}$ of $\Gamma_i=H^0(X, \LLL_{m_i})$ carefully.
See~\cite[Section 6.10]{shimadapreprintchar5} for the method.
\end{remark}
\begin{remark}\label{rem:G}
The polynomial
$$
G_1=G_{1(0)}(x, y)+G_{1(1)}(x, y)\,w+G_{1(2)}(x, y)\,w^2+G_{1(3)}(x, y)\,w^3
$$
is given in Table~\ref{table:G1s}.
The polynomial $G_2$ is too large to be presented in the paper~(see~\cite{shimadawebpage}).
\end{remark}
\begin{table}
\begin{dgroup*}
\begin{dmath*}
G_{1(0)}=
-(1-\sqrmo)\, +(1+\sqrmo)\, x+(1+\sqrmo)\, y+\sqrmo\, x^{2}-(1+\sqrmo)\, xy-(1+\sqrmo)\, y^{2}-\, xy^{2}+(1+\sqrmo)\, y^{3}-(1-\sqrmo)\, x^{4}-(1+\sqrmo)\, x^{3}y-\, xy^{3}-(1-\sqrmo)\, y^{4}-(1-\sqrmo)\, x^{5}-\, x^{3}y^{2}-\sqrmo\, x^{2}y^{3}-(1-\sqrmo)\, xy^{4}+(1+\sqrmo)\, y^{5}-(1-\sqrmo)\, x^{6}+\, x^{5}y+\sqrmo\, x^{4}y^{2}-(1-\sqrmo)\, x^{3}y^{3}+(1+\sqrmo)\, x^{2}y^{4}-\sqrmo\, xy^{5}+(1-\sqrmo)\, y^{6}+(1+\sqrmo)\, x^{7}+\, x^{4}y^{3}+(1+\sqrmo)\, x^{3}y^{4}+\sqrmo\, xy^{6}-\sqrmo\, y^{7}+\sqrmo\, x^{8}+(1+\sqrmo)\, x^{7}y+\sqrmo\, x^{6}y^{2}-\sqrmo\, x^{5}y^{3}+(1-\sqrmo)\, x^{4}y^{4}+\, x^{2}y^{6}+\sqrmo\, xy^{7}-\sqrmo\, y^{8}-(1+\sqrmo)\, x^{9}-\sqrmo\, x^{8}y-(1-\sqrmo)\, x^{7}y^{2}-(1+\sqrmo)\, x^{6}y^{3}+\sqrmo\, x^{5}y^{4}+(1-\sqrmo)\, x^{4}y^{5}-(1+\sqrmo)\, x^{3}y^{6}-(1-\sqrmo)\, x^{2}y^{7}-(1-\sqrmo)\, xy^{8}-(1+\sqrmo)\, y^{9}
\end{dmath*}
\begin{dmath*}
G_{1(1)}=
(1-\sqrmo)\, +(1-\sqrmo)\, x-(1+\sqrmo)\, x^{2}-\, xy+\sqrmo\, y^{2}+\, x^{3}-(1-\sqrmo)\, x^{2}y-\, xy^{2}+(1+\sqrmo)\, x^{4}+(1-\sqrmo)\, xy^{3}-(1-\sqrmo)\, y^{4}-\, x^{5}+\, x^{4}y+\, xy^{4}-(1+\sqrmo)\, y^{5}+(1+\sqrmo)\, x^{5}y-(1+\sqrmo)\, xy^{5}+\, y^{6}-\sqrmo\, x^{7}-\, x^{6}y+\, x^{5}y^{2}-(1-\sqrmo)\, x^{4}y^{3}+(1-\sqrmo)\, x^{2}y^{5}+(1-\sqrmo)\, xy^{6}-(1-\sqrmo)\, y^{7}-\sqrmo\, x^{8}+(1+\sqrmo)\, x^{7}y-\sqrmo\, x^{6}y^{2}-\sqrmo\, x^{5}y^{3}-(1-\sqrmo)\, x^{4}y^{4}+(1+\sqrmo)\, x^{3}y^{5}-(1-\sqrmo)\, x^{2}y^{6}+(1-\sqrmo)\, xy^{7}+(1-\sqrmo)\, y^{8}
\end{dmath*}
\begin{dmath*}
G_{1(2)}=
(1-\sqrmo)\, -(1+\sqrmo)\, x-(1+\sqrmo)\, xy+\, y^{2}-(1-\sqrmo)\, x^{3}-(1+\sqrmo)\, x^{2}y-\sqrmo\, xy^{2}-(1+\sqrmo)\, y^{3}+\, x^{4}-(1+\sqrmo)\, x^{3}y+\, xy^{3}-\, y^{4}-\, x^{5}-\, x^{4}y-\, xy^{4}-\, y^{5}-\sqrmo\, x^{6}+\, x^{4}y^{2}+\sqrmo\, x^{3}y^{3}+(1+\sqrmo)\, xy^{5}-\, y^{6}-(1-\sqrmo)\, x^{7}-\sqrmo\, x^{6}y-\sqrmo\, x^{5}y^{2}-\sqrmo\, x^{4}y^{3}-(1+\sqrmo)\, x^{3}y^{4}-(1+\sqrmo)\, x^{2}y^{5}+(1+\sqrmo)\, xy^{6}-(1+\sqrmo)\, y^{7}
\end{dmath*}
\begin{dmath*}
G_{1(3)}=
(1+\sqrmo)\, x-(1-\sqrmo)\, y-(1-\sqrmo)\, x^{2}-(1-\sqrmo)\, xy+(1+\sqrmo)\, y^{2}-\, x^{3}-\, x^{2}y-\, xy^{2}+\, y^{3}-(1-\sqrmo)\, x^{4}-\sqrmo\, x^{3}y+(1-\sqrmo)\, xy^{3}-\sqrmo\, y^{4}+\sqrmo\, x^{5}+\, x^{4}y+(1+\sqrmo)\, x^{3}y^{2}+(1+\sqrmo)\, x^{2}y^{3}+(1+\sqrmo)\, xy^{4}+(1+\sqrmo)\, y^{5}-\, x^{6}-(1-\sqrmo)\, x^{5}y+(1+\sqrmo)\, x^{4}y^{2}+\sqrmo\, x^{3}y^{3}-(1-\sqrmo)\, x^{2}y^{4}-(1+\sqrmo)\, xy^{5}+(1+\sqrmo)\, y^{6}
\end{dmath*}
\end{dgroup*}
\vskip 5pt
\caption{Polynomial $G_{1}$}\label{table:G1s}
\end{table}
\begin{table}
{\small
$$
\left[
\begin{array}{cccccccccccccccccccccc}
 -1 \shrink & 0 \shrink & 0 \shrink & 0 \shrink & 1 \shrink & 1 \shrink & 1 \shrink & 0 \shrink & 1 \shrink & 0 \shrink & -1 \shrink & 0 \shrink & 0 \shrink & 0 \shrink & 0 \shrink & 1 \shrink & 1 \shrink & 0 \shrink & 1 \shrink & -1 \shrink & 0 \shrink & 0 \\ 
\noalign{\myvskipmat} 0 \shrink & 0 \shrink & 0 \shrink & 0 \shrink & 0 \shrink & 0 \shrink & 0 \shrink & 0 \shrink & 0 \shrink & 0 \shrink & 0 \shrink & 0 \shrink & 0 \shrink & 0 \shrink & 0 \shrink & 0 \shrink & 0 \shrink & 0 \shrink & 0 \shrink & 1 \shrink & 0 \shrink & 0 \\ 
\noalign{\myvskipmat} 1 \shrink & 1 \shrink & 0 \shrink & 1 \shrink & 1 \shrink & 0 \shrink & 0 \shrink & 0 \shrink & 1 \shrink & 0 \shrink & 0 \shrink & 1 \shrink & 0 \shrink & -1 \shrink & -1 \shrink & 0 \shrink & 0 \shrink & 0 \shrink & 0 \shrink & 0 \shrink & 0 \shrink & 0 \\ 
\noalign{\myvskipmat} 0 \shrink & 0 \shrink & 0 \shrink & 0 \shrink & 0 \shrink & 0 \shrink & 0 \shrink & 0 \shrink & 0 \shrink & 1 \shrink & 0 \shrink & 0 \shrink & 0 \shrink & 0 \shrink & 0 \shrink & 0 \shrink & 0 \shrink & 0 \shrink & 0 \shrink & 0 \shrink & 0 \shrink & 0 \\ 
\noalign{\myvskipmat} 1 \shrink & 1 \shrink & 1 \shrink & 1 \shrink & 0 \shrink & 0 \shrink & 0 \shrink & 0 \shrink & 0 \shrink & 0 \shrink & 0 \shrink & 0 \shrink & 0 \shrink & -1 \shrink & -1 \shrink & -1 \shrink & 0 \shrink & 0 \shrink & 0 \shrink & 0 \shrink & 0 \shrink & 0 \\ 
\noalign{\myvskipmat} 0 \shrink & 1 \shrink & 0 \shrink & 1 \shrink & 1 \shrink & 1 \shrink & 0 \shrink & 0 \shrink & 1 \shrink & 0 \shrink & -1 \shrink & 0 \shrink & 0 \shrink & -1 \shrink & 0 \shrink & 1 \shrink & 0 \shrink & 0 \shrink & 1 \shrink & -1 \shrink & 0 \shrink & 0 \\ 
\noalign{\myvskipmat} 0 \shrink & -1 \shrink & 0 \shrink & -1 \shrink & 0 \shrink & -1 \shrink & 1 \shrink & 0 \shrink & -1 \shrink & 1 \shrink & 0 \shrink & -1 \shrink & 0 \shrink & 1 \shrink & 0 \shrink & 1 \shrink & 0 \shrink & -1 \shrink & 0 \shrink & 1 \shrink & 1 \shrink & 1 \\ 
\noalign{\myvskipmat} -2 \shrink & -1 \shrink & -2 \shrink & -2 \shrink & 2 \shrink & 1 \shrink & 1 \shrink & 1 \shrink & 2 \shrink & 1 \shrink & 0 \shrink & 1 \shrink & 1 \shrink & 0 \shrink & 0 \shrink & 1 \shrink & 1 \shrink & 0 \shrink & 0 \shrink & 0 \shrink & 0 \shrink & -1 \\ 
\noalign{\myvskipmat} 0 \shrink & 0 \shrink & 0 \shrink & 0 \shrink & 0 \shrink & 0 \shrink & 0 \shrink & 0 \shrink & 0 \shrink & 0 \shrink & 0 \shrink & 1 \shrink & 0 \shrink & 0 \shrink & 0 \shrink & 0 \shrink & 0 \shrink & 0 \shrink & 0 \shrink & 0 \shrink & 0 \shrink & 0 \\ 
\noalign{\myvskipmat} 0 \shrink & 0 \shrink & 0 \shrink & 1 \shrink & 0 \shrink & 0 \shrink & 0 \shrink & 0 \shrink & 0 \shrink & 0 \shrink & 0 \shrink & 0 \shrink & 0 \shrink & 0 \shrink & 0 \shrink & 0 \shrink & 0 \shrink & 0 \shrink & 0 \shrink & 0 \shrink & 0 \shrink & 0 \\ 
\noalign{\myvskipmat} 0 \shrink & 0 \shrink & 0 \shrink & 0 \shrink & 1 \shrink & 0 \shrink & 1 \shrink & 0 \shrink & 0 \shrink & 1 \shrink & 0 \shrink & 0 \shrink & 0 \shrink & 0 \shrink & -1 \shrink & 1 \shrink & 0 \shrink & -1 \shrink & 0 \shrink & 1 \shrink & 1 \shrink & 0 \\ 
\noalign{\myvskipmat} 0 \shrink & 0 \shrink & 0 \shrink & 0 \shrink & 0 \shrink & 0 \shrink & 0 \shrink & 0 \shrink & 1 \shrink & 0 \shrink & 0 \shrink & 0 \shrink & 0 \shrink & 0 \shrink & 0 \shrink & 0 \shrink & 0 \shrink & 0 \shrink & 0 \shrink & 0 \shrink & 0 \shrink & 0 \\ 
\noalign{\myvskipmat} 0 \shrink & 1 \shrink & 0 \shrink & 1 \shrink & 1 \shrink & 1 \shrink & 0 \shrink & 0 \shrink & 1 \shrink & 0 \shrink & -1 \shrink & 1 \shrink & 0 \shrink & -1 \shrink & 0 \shrink & 0 \shrink & 0 \shrink & 1 \shrink & 1 \shrink & -1 \shrink & -1 \shrink & 0 \\ 
\noalign{\myvskipmat} -1 \shrink & 0 \shrink & -1 \shrink & 0 \shrink & 1 \shrink & 1 \shrink & 0 \shrink & 0 \shrink & 2 \shrink & 0 \shrink & 0 \shrink & 2 \shrink & 1 \shrink & -1 \shrink & 0 \shrink & 0 \shrink & 1 \shrink & 1 \shrink & 1 \shrink & -1 \shrink & -1 \shrink & -1 \\ 
\noalign{\myvskipmat} 1 \shrink & 1 \shrink & 1 \shrink & 1 \shrink & 0 \shrink & 0 \shrink & 1 \shrink & 0 \shrink & 0 \shrink & 1 \shrink & -1 \shrink & -1 \shrink & -1 \shrink & 0 \shrink & -1 \shrink & 0 \shrink & 0 \shrink & -1 \shrink & 0 \shrink & 1 \shrink & 1 \shrink & 1 \\ 
\noalign{\myvskipmat} 0 \shrink & 0 \shrink & 0 \shrink & 0 \shrink & 0 \shrink & 0 \shrink & 0 \shrink & 0 \shrink & 0 \shrink & 0 \shrink & 0 \shrink & 0 \shrink & 0 \shrink & 0 \shrink & 0 \shrink & 1 \shrink & 0 \shrink & 0 \shrink & 0 \shrink & 0 \shrink & 0 \shrink & 0 \\ 
\noalign{\myvskipmat} 1 \shrink & 1 \shrink & 1 \shrink & 1 \shrink & 0 \shrink & 0 \shrink & 0 \shrink & 0 \shrink & 0 \shrink & 0 \shrink & -1 \shrink & -1 \shrink & -1 \shrink & 0 \shrink & 0 \shrink & 0 \shrink & 0 \shrink & 0 \shrink & 0 \shrink & 0 \shrink & 0 \shrink & 0 \\ 
\noalign{\myvskipmat} 2 \shrink & 2 \shrink & 2 \shrink & 3 \shrink & 0 \shrink & -1 \shrink & -1 \shrink & -1 \shrink & 0 \shrink & 0 \shrink & 0 \shrink & 1 \shrink & 0 \shrink & -1 \shrink & -1 \shrink & 0 \shrink & 0 \shrink & 0 \shrink & 0 \shrink & 0 \shrink & -1 \shrink & 0 \\ 
\noalign{\myvskipmat} -2 \shrink & -1 \shrink & -2 \shrink & -2 \shrink & 1 \shrink & 1 \shrink & 1 \shrink & 0 \shrink & 1 \shrink & 0 \shrink & 1 \shrink & 1 \shrink & 1 \shrink & 1 \shrink & 0 \shrink & 0 \shrink & 1 \shrink & 0 \shrink & 0 \shrink & 0 \shrink & 0 \shrink & -1 \\ 
\noalign{\myvskipmat} 0 \shrink & 1 \shrink & 0 \shrink & 0 \shrink & 0 \shrink & 0 \shrink & 0 \shrink & 0 \shrink & 0 \shrink & 0 \shrink & 0 \shrink & 0 \shrink & 0 \shrink & 0 \shrink & 0 \shrink & 0 \shrink & 0 \shrink & 0 \shrink & 0 \shrink & 0 \shrink & 0 \shrink & 0 \\ 
\noalign{\myvskipmat} 2 \shrink & 2 \shrink & 2 \shrink & 3 \shrink & 0 \shrink & 0 \shrink & -1 \shrink & -1 \shrink & 1 \shrink & 0 \shrink & 0 \shrink & 1 \shrink & -1 \shrink & -1 \shrink & -1 \shrink & -1 \shrink & 0 \shrink & 0 \shrink & 0 \shrink & 0 \shrink & -1 \shrink & 0 \\ 
\noalign{\myvskipmat} 1 \shrink & 1 \shrink & 1 \shrink & 1 \shrink & 0 \shrink & 0 \shrink & 0 \shrink & 0 \shrink & 0 \shrink & 0 \shrink & -1 \shrink & 0 \shrink & 0 \shrink & -1 \shrink & 0 \shrink & 0 \shrink & -1 \shrink & 0 \shrink & 0 \shrink & 0 \shrink & 0 \shrink & 0 
\end{array}
\right]
$$
}
\caption{The matrix $A_1$}\label{table:A1}
\end{table}
\begin{table}
{\small
$$
\left[
\begin{array}{cccccccccccccccccccccc}
 1 \shrink & 1 \shrink & -1 \shrink & 0 \shrink & 2 \shrink & 0 \shrink & 2 \shrink & 1 \shrink & 1 \shrink & 0 \shrink & 0 \shrink & -1 \shrink & 2 \shrink & 1 \shrink & 0 \shrink & 1 \shrink & 1 \shrink & -1 \shrink & 0 \shrink & 0 \shrink & 0 \shrink & 0 \\ 
\noalign{\myvskipmat} 0 \shrink & 0 \shrink & 0 \shrink & 1 \shrink & 0 \shrink & 0 \shrink & 0 \shrink & 0 \shrink & 0 \shrink & 0 \shrink & 0 \shrink & 0 \shrink & 0 \shrink & 0 \shrink & 0 \shrink & 0 \shrink & 0 \shrink & 0 \shrink & 0 \shrink & 0 \shrink & 0 \shrink & 0 \\ 
\noalign{\myvskipmat} 0 \shrink & 2 \shrink & -1 \shrink & 0 \shrink & 4 \shrink & 0 \shrink & 4 \shrink & 2 \shrink & 2 \shrink & 0 \shrink & 0 \shrink & -2 \shrink & 4 \shrink & 2 \shrink & 0 \shrink & 2 \shrink & 2 \shrink & -2 \shrink & 0 \shrink & 0 \shrink & 0 \shrink & 0 \\ 
\noalign{\myvskipmat} 0 \shrink & 1 \shrink & 0 \shrink & 0 \shrink & 0 \shrink & 0 \shrink & 0 \shrink & 0 \shrink & 0 \shrink & 0 \shrink & 0 \shrink & 0 \shrink & 0 \shrink & 0 \shrink & 0 \shrink & 0 \shrink & 0 \shrink & 0 \shrink & 0 \shrink & 0 \shrink & 0 \shrink & 0 \\ 
\noalign{\myvskipmat} 0 \shrink & 0 \shrink & 0 \shrink & 0 \shrink & 1 \shrink & 0 \shrink & 0 \shrink & 0 \shrink & 0 \shrink & 0 \shrink & 0 \shrink & 0 \shrink & 0 \shrink & 0 \shrink & 0 \shrink & 0 \shrink & 0 \shrink & 0 \shrink & 0 \shrink & 0 \shrink & 0 \shrink & 0 \\ 
\noalign{\myvskipmat} 0 \shrink & 2 \shrink & -2 \shrink & 0 \shrink & 4 \shrink & 0 \shrink & 4 \shrink & 2 \shrink & 2 \shrink & 0 \shrink & 0 \shrink & -1 \shrink & 4 \shrink & 2 \shrink & 0 \shrink & 2 \shrink & 2 \shrink & -2 \shrink & 0 \shrink & 0 \shrink & 0 \shrink & 0 \\ 
\noalign{\myvskipmat} 0 \shrink & 0 \shrink & 0 \shrink & 0 \shrink & -1 \shrink & -1 \shrink & 0 \shrink & 0 \shrink & 0 \shrink & 1 \shrink & 1 \shrink & 0 \shrink & 0 \shrink & 1 \shrink & 0 \shrink & 0 \shrink & 0 \shrink & 0 \shrink & -1 \shrink & 1 \shrink & 0 \shrink & 0 \\ 
\noalign{\myvskipmat} 0 \shrink & 2 \shrink & 0 \shrink & 1 \shrink & 1 \shrink & 0 \shrink & 2 \shrink & 0 \shrink & 1 \shrink & 0 \shrink & 0 \shrink & -1 \shrink & 2 \shrink & 1 \shrink & 0 \shrink & 1 \shrink & 1 \shrink & -1 \shrink & 0 \shrink & 0 \shrink & 0 \shrink & 0 \\ 
\noalign{\myvskipmat} 4 \shrink & 2 \shrink & 3 \shrink & 3 \shrink & -1 \shrink & -2 \shrink & -1 \shrink & 0 \shrink & -1 \shrink & 0 \shrink & 0 \shrink & -1 \shrink & -1 \shrink & -1 \shrink & -1 \shrink & -1 \shrink & -1 \shrink & -1 \shrink & -1 \shrink & 1 \shrink & 0 \shrink & 1 \\ 
\noalign{\myvskipmat} -3 \shrink & -1 \shrink & -4 \shrink & -3 \shrink & 4 \shrink & 2 \shrink & 3 \shrink & 2 \shrink & 2 \shrink & 0 \shrink & 0 \shrink & 0 \shrink & 3 \shrink & 1 \shrink & 1 \shrink & 2 \shrink & 2 \shrink & 0 \shrink & 1 \shrink & -1 \shrink & 0 \shrink & -1 \\ 
\noalign{\myvskipmat} 2 \shrink & 2 \shrink & 1 \shrink & 2 \shrink & 1 \shrink & -2 \shrink & 2 \shrink & 1 \shrink & 0 \shrink & 1 \shrink & 0 \shrink & -2 \shrink & 1 \shrink & 1 \shrink & -1 \shrink & 1 \shrink & 0 \shrink & -2 \shrink & -1 \shrink & 1 \shrink & 1 \shrink & 1 \\ 
\noalign{\myvskipmat} 0 \shrink & 2 \shrink & -2 \shrink & 0 \shrink & 4 \shrink & 1 \shrink & 4 \shrink & 2 \shrink & 2 \shrink & 0 \shrink & 0 \shrink & -2 \shrink & 4 \shrink & 2 \shrink & 0 \shrink & 2 \shrink & 2 \shrink & -2 \shrink & 0 \shrink & 0 \shrink & 0 \shrink & 0 \\ 
\noalign{\myvskipmat} -1 \shrink & 0 \shrink & -1 \shrink & -1 \shrink & 1 \shrink & 1 \shrink & 0 \shrink & 0 \shrink & 1 \shrink & -1 \shrink & 0 \shrink & 1 \shrink & 1 \shrink & 0 \shrink & 0 \shrink & 0 \shrink & 1 \shrink & 1 \shrink & 1 \shrink & -1 \shrink & -1 \shrink & -1 \\ 
\noalign{\myvskipmat} 1 \shrink & 2 \shrink & 0 \shrink & 1 \shrink & 2 \shrink & 0 \shrink & 2 \shrink & 0 \shrink & 0 \shrink & -1 \shrink & 0 \shrink & -1 \shrink & 2 \shrink & 1 \shrink & 0 \shrink & 1 \shrink & 1 \shrink & -1 \shrink & 0 \shrink & 0 \shrink & 0 \shrink & 0 \\ 
\noalign{\myvskipmat} 1 \shrink & 1 \shrink & 1 \shrink & 1 \shrink & 0 \shrink & 0 \shrink & 0 \shrink & 0 \shrink & 0 \shrink & 0 \shrink & -1 \shrink & -1 \shrink & -1 \shrink & 0 \shrink & 0 \shrink & 0 \shrink & 0 \shrink & 0 \shrink & 0 \shrink & 0 \shrink & 0 \shrink & 0 \\ 
\noalign{\myvskipmat} 0 \shrink & 0 \shrink & -1 \shrink & -1 \shrink & 2 \shrink & 0 \shrink & 3 \shrink & 2 \shrink & 1 \shrink & 1 \shrink & 0 \shrink & -2 \shrink & 2 \shrink & 1 \shrink & 0 \shrink & 1 \shrink & 1 \shrink & -2 \shrink & 0 \shrink & 0 \shrink & 1 \shrink & 1 \\ 
\noalign{\myvskipmat} -2 \shrink & -1 \shrink & -2 \shrink & -2 \shrink & 1 \shrink & 2 \shrink & 1 \shrink & 1 \shrink & 1 \shrink & 0 \shrink & -1 \shrink & 0 \shrink & 1 \shrink & 0 \shrink & 1 \shrink & 1 \shrink & 0 \shrink & 0 \shrink & 1 \shrink & -1 \shrink & 0 \shrink & 0 \\ 
\noalign{\myvskipmat} 1 \shrink & 2 \shrink & 0 \shrink & 1 \shrink & 1 \shrink & -1 \shrink & 1 \shrink & 0 \shrink & 1 \shrink & 0 \shrink & 1 \shrink & 0 \shrink & 2 \shrink & 1 \shrink & 0 \shrink & 0 \shrink & 1 \shrink & 0 \shrink & 0 \shrink & 0 \shrink & -1 \shrink & 0 \\ 
\noalign{\myvskipmat} 0 \shrink & 2 \shrink & 0 \shrink & 1 \shrink & 2 \shrink & 0 \shrink & 2 \shrink & 1 \shrink & 1 \shrink & 0 \shrink & 0 \shrink & -1 \shrink & 2 \shrink & 1 \shrink & 0 \shrink & 1 \shrink & 1 \shrink & -1 \shrink & -1 \shrink & 0 \shrink & 0 \shrink & -1 \\ 
\noalign{\myvskipmat} 0 \shrink & 1 \shrink & 1 \shrink & 1 \shrink & 0 \shrink & 0 \shrink & 0 \shrink & 0 \shrink & 1 \shrink & 0 \shrink & 0 \shrink & 1 \shrink & 0 \shrink & 0 \shrink & 0 \shrink & -1 \shrink & 0 \shrink & 0 \shrink & -1 \shrink & 0 \shrink & -1 \shrink & -1 \\ 
\noalign{\myvskipmat} 4 \shrink & 5 \shrink & 1 \shrink & 4 \shrink & 3 \shrink & -1 \shrink & 2 \shrink & 1 \shrink & 1 \shrink & -1 \shrink & 0 \shrink & -2 \shrink & 3 \shrink & 1 \shrink & -1 \shrink & 1 \shrink & 1 \shrink & -2 \shrink & 0 \shrink & 0 \shrink & -1 \shrink & 0 \\ 
\noalign{\myvskipmat} 0 \shrink & 0 \shrink & 0 \shrink & 0 \shrink & 0 \shrink & 0 \shrink & 0 \shrink & 0 \shrink & 0 \shrink & 0 \shrink & 0 \shrink & 0 \shrink & 0 \shrink & 0 \shrink & 0 \shrink & 0 \shrink & 0 \shrink & 0 \shrink & 0 \shrink & 0 \shrink & 0 \shrink & 1 
\end{array}
\right]
$$
}
\caption{The matrix $A_2$}\label{table:A2}
\end{table}
\begin{proposition}\label{prop:A}
Let $g_1$ and  $g_2$  be the involutions of $X$  defined in Theorem~\ref{thm:main}.
Then the action $g_{i*}$ on $S$ is given by $v\mapsto v A_i$,
where $A_i$ is the matrix given in Tables~\ref{table:A1} and~\ref{table:A2}.
\end{proposition}
\begin{proof} 
Recall that $\Exc(\phi_i)$  is the set of the classes of $(-2)$-curves 
contracted by $\phi_i: X\to \P^2$.
Suppose that $\gamma_1, \dots, \gamma_r \in \Exc(\phi_i)$
are the classes of $(-2)$-curves that are contracted 
to a  singular point $P\in \Sing (B_i)$ of type $A_r$.
We index them in such a way that
$\intS{\gamma_{\nu}, \gamma_{\nu+1}}=1$
holds for $\nu=1, \dots, r-1$.
Then $g_{i*}$ interchanges $\gamma_i$ and $\gamma_{r+1-i}$.
Let $V(P)\subset S\tQ$ denote the linear span of 
the invariant vectors $\gamma_i+\gamma_{r+1-i}$.
Then the eigenspace of  $g_{i*}$ on $S\tQ$
with eigenvalue $1$ is equal to
$$
\gen{m_i}\oplus\bigoplus_{P\in \Sing (B_i)}V(P), 
$$
and the eigenspace with eigenvalue $-1$ is its orthogonal complement.
\end{proof}
Using the matrix representations $A_i$ of $g_{i*}$,
 we verify the following facts:
 \begin{enumerate}
\item The eigenspace of $g_{i*}$ with eigenvalue $1$ is contained in $(b_i)\sperp_S$.
In particular, we have $b_i^{g_i}=-b_i$.
\item
The vector $h_0^{g_i}$ is equal to  
the image of $h_0$ by the reflection into the wall $(b_i)_S\sperp$,   that is
$h_0^{g_1}=h_0+3\,b_1$ and $h_0^{g_2}=h_0+9\,b_2$ hold.
\end{enumerate}
Since $\Aut(X, h_0)$ acts on each of $\wallrsz{648}$ and $\wallrsz{5184}$ transitively,
we obtain the following:
\begin{corollary}\label{cor:imageA}
For any $r_S\in \wallrsz{648} \cup \wallrsz{5184}$,
there exists $\tau \in \Aut(X, h_0)$ such that
$$
h_0^{g_i\tau}=h_0+ c_i\, r_S
$$
holds, where $i=1$ and $c_1=3$ if $r_S\in \wallrsz{648}$ while $i=2$ and $c_2=9$ if $r_S\in  \wallrsz{5184}$.
\end{corollary}
\section{Proof of Theorem~\ref{thm:main}}\label{sec:mainproof}
We denote by
$$
G:=\gen{\Aut(X, h_0), g_1, g_2}
$$
the subgroup of $\Aut(X)$ generated by $\Aut(X, h_0)$, $g_1$ and $g_2$.
Note that the action of $\Aut(X)$ on $S$ preserves the set of nef classes.
\begin{theorem}\label{thm:D}
If $v\in S$ is nef, 
there exists $\gamma\in G$ such that $v^\gamma\in \chamsz$.
\end{theorem}
\begin{proof}
Let $\gamma\in G$ be an element such that $\intS{v^\gamma, h_0}$ attains 
$$
\min\shortset{\intS{v^{\gamma\sprime}, h_0}}{\gamma\sprime\in G}.
$$
We show that $\intS{v^\gamma, r_S}\ge 0$ holds for any $r_S\in \wallrs(\chamsz)$.
If $r_S\in \wallrsz{112}$, then  $r_S=[\ell_i]$ for some line $\ell_i\subset X$,
and hence $\intS{v^\gamma, r_S}\ge 0$ holds because $v^\gamma$ is nef.
Suppose that  $r_S\in \wallrsz{648} \cup \wallrsz{5184}$.
By Corollary~\ref{cor:imageA}, 
there exists $\tau \in \Aut(X, h_0)$  such that
$h_0^{g_i\tau}=h_0+ c_i r_S$
holds, where $i=1$ and $c_1=3$ if $r_S\in \wallrsz{648}$ while $i=2$ and $c_2=9$ if $r_S\in  \wallrsz{5184}$.
Since $\gamma\tau\inv g_i\in G$,
we have
$$
\intS{v^{\gamma}, h_0}\le 
\intS{v^{\gamma \tau\inv g_i}, h_0}
=\intS{v^{\gamma }, h_0^{g_i\tau}}
=\intS{v^{\gamma }, h_0}+ c_i \intS{v^{\gamma }, r_S}.
$$
Therefore $\intS{v^{\gamma }, r_S}\ge 0$ holds.
\end{proof}
The properties (1), (2), (3) of $\chamsz$ stated in Introduction follow from 
Corollaries~\ref{cor:action_on_CS0},~\ref{cor:equiv},~\ref{cor:converse} and 
Theorem~\ref{thm:D}.
\begin{proof}[Proof of Theorem~\ref{thm:main}]
By Corollary~\ref{cor:converse},  
it is enough to show that,
for any $\gamma\in \Aut(X)$,
there exists $\gamma\sprime\in G$ such that $h_0^{\gamma\gamma\sprime}\in \chamsz$ holds.
Since $h_0^{\gamma}$ is nef, this follows from Theorem~\ref{thm:D}.
\end{proof}
\section{The Fermat quartic polarizations for $g_1$ and $g_2$}\label{sec:quartics}
A polarization $h\in S$ of degree $4$ is said to be a \emph{Fermat quartic polarization}
if, by choosing an appropriate basis of $H^0(X, \LLL_h)$,
the morphism $\Phi_h: X\to \P^3$ 
associated with $|\LLL_h|$ induces an automorphism 
of $X\subset \P^3$.
It is obvious that
$h_0^{\gamma}$ is a Fermat quartic polarization for any $\gamma\in \Aut(X)$.
Conversely,
if $h$ is a Fermat quartic polarization, then
the pull-back of $h_0$ by the automorphism $\Phi_h$ of $X$
is $h$.
Therefore the set of Fermat quartic polarizations is the orbit of $h_0$
by  the action of $\Aut(X)$ on $S$.
Consider the Fermat quartic polarizations
\begin{eqnarray*}
h_1 &:=& h_0^{g_1}=h_0 A_1=[0, 1, 0, 1, 2, 1, 1, 0, 2, 1, -1, 1, 0, -1, -1, 1, 1, 0, 1, 0, 0, 0]_S, \\
h_2 &:=& h_0^{g_2}= h_0 A_2=[1, 4, -2, 1, 6, 0, 6, 3, 3, 0, 0, -3, 6, 3, 0, 3, 3, -3, 0, 0, 0, 0]_S.
\end{eqnarray*}
Using the equalities 
\begin{eqnarray}
h_1&=& 6h_0-([\ell_{3}]+[\ell_{6}]+[\ell_{8}]+[\ell_{14}]+[\ell_{15}]+[\ell_{17}]+[\ell_{19}]+ \label{eq:hh1}\\
&& \phantom{a}+ [\ell_{22}]+[\ell_{31}]+[\ell_{34}]+[\ell_{63}]+[\ell_{70}]+[\ell_{79}]+[\ell_{92}]), \quand \nonumber \\
h_2&=& 15h_0- (3\,[\ell_{3}]+4\,[\ell_{6}]+[\ell_{13}]+[\ell_{14}]+ 3\,[\ell_{18}]+[\ell_{22}]+
 \label{eq:hh2}\\
&&  \phantom{a}+  [\ell_{26}]+ [\ell_{27}]+ 2\,[\ell_{35}]+ [\ell_{44}]+ 2\,[\ell_{50}]+ 3\,[\ell_{92}]+  \nonumber \\
&& \phantom{a}+[\ell_{93}]+ [\ell_{106}]+[\ell_{108}]+3\,[\ell_{111}]),
\nonumber 
\end{eqnarray}
we obtain another description of the involutions $g_1$ and $g_2$. 
\begin{theorem}\label{thm:H}
Let $(w,x,y)$ be the affine coordinates of $\P^3$ with $z=1$,
and let  
$$
H_{1j}(w,x,y)=H_{1j0}(x,y)+H_{1j1}(x,y)\, w+H_{1j2}(x,y)\, w^2+H_{1j3}(x,y)\,w^3
$$
be polynomials given in Table~\ref{table:H1s}.
Then the rational map
\begin{equation}\label{isom:H}
(w,x,y)\mapsto [H_{10}: H_{11}: H_{12}: H_{13}]\in \P^3
\end{equation}
gives the involution  $g_1$ of $X$. 
\end{theorem}
\begin{remark}
We have a similar list  of polynomials 
$H_{20}, H_{21}, H_{22}, H_{23}$
that gives the involution $g_2$.
They are, however, too large to be presented 
in the paper~(see~\cite{shimadawebpage}).
\end{remark}
\begin{table}
{\small
\begin{dgroup*}
\begin{dmath*}
H_{100}=
-1\, -(1-\sqrmo)\, x-(1-\sqrmo)\, x^{2}-(1+\sqrmo)\, y^{2}-\sqrmo\, x^{3}-(1-\sqrmo)\, xy^{2}+\, x^{4}-(1-\sqrmo)\, x^{3}y+(1+\sqrmo)\, x^{2}y^{2}+(1+\sqrmo)\, xy^{3}-(1-\sqrmo)\, y^{4}-(1-\sqrmo)\, x^{5}+(1+\sqrmo)\, x^{4}y-(1-\sqrmo)\, x^{3}y^{2}-\sqrmo\, x^{2}y^{3}-\sqrmo\, xy^{4}+\sqrmo\, y^{5}+\sqrmo\, x^{6}-\, x^{5}y-(1+\sqrmo)\, x^{4}y^{2}-(1+\sqrmo)\, x^{3}y^{3}+(1+\sqrmo)\, xy^{5}+(1-\sqrmo)\, y^{6}
\end{dmath*}
\begin{dmath*}
H_{101}=
(1+\sqrmo)\, +\, x+(1-\sqrmo)\, y-\sqrmo\, x^{2}-(1+\sqrmo)\, xy+\sqrmo\, y^{2}+\sqrmo\, x^{3}-\, x^{2}y-\sqrmo\, xy^{2}+\, y^{3}+\, x^{3}y+(1-\sqrmo)\, x^{2}y^{2}-(1-\sqrmo)\, xy^{3}+(1-\sqrmo)\, y^{4}+\, x^{5}+\sqrmo\, x^{4}y+\, x^{3}y^{2}-\sqrmo\, x^{2}y^{3}+(1-\sqrmo)\, xy^{4}-(1-\sqrmo)\, y^{5}
\end{dmath*}
\begin{dmath*}
H_{102}=
\sqrmo\, +\, x+(1+\sqrmo)\, y+\, x^{2}+(1+\sqrmo)\, y^{2}+(1+\sqrmo)\, x^{3}-\, x^{2}y-(1+\sqrmo)\, xy^{2}+\sqrmo\, y^{3}+(1+\sqrmo)\, x^{4}+(1-\sqrmo)\, x^{2}y^{2}+\, xy^{3}+(1+\sqrmo)\, y^{4}
\end{dmath*}
\begin{dmath*}
H_{103}=
(1-\sqrmo)\, -(1-\sqrmo)\, x+(1-\sqrmo)\, y-(1+\sqrmo)\, x^{2}-(1-\sqrmo)\, xy+(1+\sqrmo)\, x^{3}-(1+\sqrmo)\, y^{3}
\end{dmath*}
\begin{dmath*}
\hbox{---------------------------------------------------------------}\mystruthd{10pt}{5pt}
\end{dmath*}
\begin{dmath*}
H_{110}=
-\sqrmo\, +\sqrmo\, x+\, y-(1+\sqrmo)\, x^{2}+\, xy-(1-\sqrmo)\, y^{2}-\, x^{3}-(1-\sqrmo)\, x^{2}y+(1+\sqrmo)\, xy^{2}-\, y^{3}+(1+\sqrmo)\, x^{4}-\sqrmo\, x^{3}y-(1-\sqrmo)\, x^{2}y^{2}+\, xy^{3}+(1-\sqrmo)\, y^{4}-(1+\sqrmo)\, x^{5}+(1-\sqrmo)\, x^{4}y+\sqrmo\, x^{2}y^{3}-(1+\sqrmo)\, xy^{4}-(1+\sqrmo)\, y^{5}-\sqrmo\, x^{6}+(1+\sqrmo)\, x^{4}y^{2}+(1+\sqrmo)\, x^{3}y^{3}+(1-\sqrmo)\, xy^{5}-(1-\sqrmo)\, y^{6}
\end{dmath*}
\begin{dmath*}
H_{111}=
-(1-\sqrmo)\, +\, x+(1+\sqrmo)\, y-(1+\sqrmo)\, x^{2}-\sqrmo\, xy-\sqrmo\, y^{2}+(1-\sqrmo)\, x^{3}-\sqrmo\, x^{2}y-\, y^{3}+(1+\sqrmo)\, x^{4}+(1-\sqrmo)\, x^{3}y-\, x^{2}y^{2}+(1+\sqrmo)\, xy^{3}-\sqrmo\, y^{4}+(1-\sqrmo)\, x^{5}+\sqrmo\, x^{4}y+(1-\sqrmo)\, x^{2}y^{3}+(1-\sqrmo)\, xy^{4}-\, y^{5}
\end{dmath*}
\begin{dmath*}
H_{112}=
-1\, +(1+\sqrmo)\, y+\, x^{2}-(1-\sqrmo)\, xy-(1+\sqrmo)\, y^{2}-\, x^{2}y+(1+\sqrmo)\, xy^{2}-(1+\sqrmo)\, y^{3}+(1-\sqrmo)\, x^{4}+(1+\sqrmo)\, x^{3}y-(1+\sqrmo)\, x^{2}y^{2}-\, xy^{3}-(1+\sqrmo)\, y^{4}
\end{dmath*}
\begin{dmath*}
H_{113}=
(1+\sqrmo)\, -\, x+\, y+\, x^{2}-\sqrmo\, y^{2}-(1-\sqrmo)\, x^{3}+\sqrmo\, x^{2}y-(1-\sqrmo)\, xy^{2}-\sqrmo\, y^{3}
\end{dmath*}
\begin{dmath*}
\hbox{---------------------------------------------------------------}\mystruthd{10pt}{5pt}
\end{dmath*}
\begin{dmath*}
H_{120}=
(1+\sqrmo)\, +(1+\sqrmo)\, x+(1+\sqrmo)\, y+(1-\sqrmo)\, x^{2}+\, y^{2}+(1+\sqrmo)\, x^{3}+(1+\sqrmo)\, x^{2}y-\sqrmo\, xy^{2}-\, y^{3}-(1-\sqrmo)\, x^{3}y+(1-\sqrmo)\, x^{2}y^{2}-(1+\sqrmo)\, xy^{3}+(1-\sqrmo)\, y^{4}+(1-\sqrmo)\, x^{5}-\sqrmo\, x^{4}y+(1-\sqrmo)\, x^{3}y^{2}-(1+\sqrmo)\, x^{2}y^{3}+\sqrmo\, xy^{4}-\, y^{5}+\, x^{6}-(1+\sqrmo)\, x^{5}y-(1-\sqrmo)\, x^{4}y^{2}+\, x^{3}y^{3}-\sqrmo\, x^{2}y^{4}-(1-\sqrmo)\, xy^{5}+(1-\sqrmo)\, y^{6}
\end{dmath*}
\begin{dmath*}
H_{121}=
\sqrmo\, +\, x+\, xy-(1+\sqrmo)\, y^{2}+\, x^{3}-(1+\sqrmo)\, x^{2}y-(1-\sqrmo)\, xy^{2}+(1+\sqrmo)\, y^{3}+\, x^{4}-(1-\sqrmo)\, x^{3}y-(1-\sqrmo)\, x^{2}y^{2}+(1+\sqrmo)\, xy^{3}-(1-\sqrmo)\, y^{4}-(1-\sqrmo)\, x^{5}+(1+\sqrmo)\, x^{3}y^{2}+(1+\sqrmo)\, x^{2}y^{3}+(1-\sqrmo)\, y^{5}
\end{dmath*}
\begin{dmath*}
H_{122}=
(1-\sqrmo)\, -\, x-(1+\sqrmo)\, y+\sqrmo\, x^{2}-(1-\sqrmo)\, xy-(1+\sqrmo)\, y^{2}-\, x^{3}-(1-\sqrmo)\, xy^{2}-\sqrmo\, y^{3}-(1+\sqrmo)\, x^{4}-(1-\sqrmo)\, x^{3}y-(1+\sqrmo)\, x^{2}y^{2}-\, xy^{3}+(1+\sqrmo)\, y^{4}
\end{dmath*}
\begin{dmath*}
H_{123}=
1\, -(1+\sqrmo)\, x+(1-\sqrmo)\, y+\, x^{2}+\sqrmo\, xy+\sqrmo\, y^{2}-(1+\sqrmo)\, x^{3}+(1-\sqrmo)\, x^{2}y-(1+\sqrmo)\, xy^{2}+(1+\sqrmo)\, y^{3}
\end{dmath*}
\begin{dmath*}
\hbox{---------------------------------------------------------------}\mystruthd{10pt}{5pt}
\end{dmath*}
\begin{dmath*}
H_{130}=
-(1-\sqrmo)\, +\sqrmo\, x+(1+\sqrmo)\, y-(1+\sqrmo)\, x^{2}+(1-\sqrmo)\, xy+(1-\sqrmo)\, y^{2}+\, x^{3}-(1+\sqrmo)\, x^{2}y+\sqrmo\, xy^{2}+\sqrmo\, y^{3}-(1+\sqrmo)\, x^{4}+\sqrmo\, x^{3}y+\, x^{2}y^{2}-(1+\sqrmo)\, y^{4}+(1+\sqrmo)\, x^{5}-(1-\sqrmo)\, x^{4}y+(1-\sqrmo)\, x^{3}y^{2}-\, x^{2}y^{3}-(1+\sqrmo)\, y^{5}-(1+\sqrmo)\, x^{6}-(1-\sqrmo)\, x^{5}y-(1+\sqrmo)\, x^{4}y^{2}+\sqrmo\, x^{3}y^{3}+\sqrmo\, x^{2}y^{4}+\sqrmo\, xy^{5}+(1+\sqrmo)\, y^{6}
\end{dmath*}
\begin{dmath*}
H_{131}=
-1\, -\, x+(1+\sqrmo)\, y-(1-\sqrmo)\, x^{2}+(1+\sqrmo)\, xy-\sqrmo\, y^{2}-(1+\sqrmo)\, x^{3}-\sqrmo\, x^{2}y-\, xy^{2}+\sqrmo\, y^{3}-\, x^{4}-\, x^{3}y+\, xy^{3}-(1+\sqrmo)\, y^{4}-(1+\sqrmo)\, x^{5}+\, x^{4}y+(1-\sqrmo)\, x^{3}y^{2}-\sqrmo\, x^{2}y^{3}+(1+\sqrmo)\, xy^{4}
\end{dmath*}
\begin{dmath*}
H_{132}=
(1+\sqrmo)\, +\sqrmo\, x+\, y-\, x^{2}+\, xy+\, y^{2}+\sqrmo\, x^{3}-(1-\sqrmo)\, x^{2}y-(1+\sqrmo)\, xy^{2}-(1-\sqrmo)\, x^{4}-\, x^{2}y^{2}-\sqrmo\, xy^{3}-(1-\sqrmo)\, y^{4}
\end{dmath*}
\begin{dmath*}
H_{133}=
\sqrmo\, -\, y+\, x^{2}+(1+\sqrmo)\, xy-(1-\sqrmo)\, y^{2}
\end{dmath*}
\end{dgroup*}
\vskip 5pt
}
\caption{Polynomials $H_{1j}$}\label{table:H1s}  
\end{table}
\begin{proof}[Proof of Theorem~\ref{thm:H}]
We put
$$
Z:=\{3, 6, 8, 14, 15, 17,
19, 22, 31, 34, 63, 70, 79, 92\},
$$
which is the set of indices of lines on $X$ that appear in the right-hand side 
of~\eqref{eq:hh1}. 
The polynomials 
$H_{10}, H_{11}, H_{12},  H_{13}$ form  a basis of the vector space 
$$
H^0(X, \LLL_{h_1})\cong\Gamma(6, \bigcap_{i\in Z} I_i).
$$
(See Section~\ref{sec:A} for the notation.) 
We can easily verify that
$$
H_{10}^4+ H_{11}^4+ H_{12}^4+ H_{13}^4\equiv 0 \;\;\bmod (w^4+x^4+y^4+1)
$$
holds.
Hence the rational map~\eqref{isom:H} induces an automorphism  $g\sprime$ of $X$.
We prove  $g\sprime=g_1$ by showing that the action $g\sprime_*$ of $g\sprime$
on $S$ is equal to the action $v\mapsto vA_1$ of $g_1$.
We  homogenize the polynomials $H_{1j}$ to $\tilde{H}_{1j}(w,x,y,z)$
so that $g\sprime$ is  given by
$$
[w:x:y:z]\mapsto [\tilde{H}_{10}: \tilde{H}_{11}: \tilde{H}_{12}: \tilde{H}_{13}].
$$
Let $\ell_k$ be a line on $X$ 
whose index $k$  is not in $Z$.
We calculate a parametric representation
$$
[u:v]\mapsto [l_{k0}: l_{k1}: l_{k2}: l_{k3}] 
$$
of $\ell_k$ in $\P^3$,
where $u, v$ are homogeneous coordinates of $\P^1$ and
$l_{k\nu}$ are homogeneous linear polynomials of $u, v$.
We put
$$
L_{1j}^{(k)}:=\tilde{H}_{1j}(l_{k0}, l_{k1}, l_{k2}, l_{k3})
$$
for $j=0, \dots, 3$, 
which are homogeneous polynomials of $u, v$.
Let $M^{(k)}$ be the greatest common divisor of
$L_{10}^{(k)}, L_{11}^{(k)}, L_{12}^{(k)}, L_{13}^{(k)}$ in $\F_9[u, v]$.
Then 
$$
\rho_k\;\; : \;\;[u:v]\mapsto [L_{10}^{(k)}/M^{(k)}: L_{11}^{(k)}/M^{(k)}: L_{12}^{(k)}/M^{(k)}: L_{13}^{(k)}/M^{(k)}]  
$$
is a parametric representation of the image of $\ell_k$ by $g\sprime$.
(If $k\in Z$, then $L_{1j}^{(k)}$ are constantly equal to $0$.)
Pulling back the defining homogeneous ideal of $\ell_{k\sprime}$ by $\rho_k$,
we can calculate the intersection number $\intS{[\ell_{k}]^{g\sprime}, [\ell_{k\sprime}]}$.
Since the classes $[\ell_{k}]$ with $k\notin Z$ span $S\tQ$,
we can calculate the action $g\sprime_*$ of $g\sprime$ on $S$,
which turns out to be equal to $v\mapsto vA_1$.
\end{proof}
\begin{remark}
The polynomials $H_{10}, H_{11}, H_{12},  H_{13}$ are found by the following method.
Let $H\sprime_{0}, H\sprime_{1}, H\sprime_{2},  H\sprime_{3}$ be an arbitrary basis of
$\Gamma(6, \cap_{i\in Z} I_i)\cong H^0(X, \LLL_{h_1})$.
Then the normal forms of the quartic monomials  of $H\sprime_{0}, H\sprime_{1}, H\sprime_{2},  H\sprime_{3}$
are subject to a  linear relation 
of the following form~(see~\cite[n.~3]{MR0213949} or~\cite[Theorem 6.11]{shimadapreprintchar5}):
$$
\sum _{i, j=0}^3 a_{ij}\overline{H\sprime_i H_j^{\prime3}}=0,
$$
where the coefficients $a_{ij}\in \F_9$ satisfy $a_{ji}=a_{ij}^3$ and $\det(a_{ij})\ne 0$;
that is, the matrix $(a_{ij})$ is non-singular Hermitian.
We search for $B\in \GL_{3}(\F_9)$ such that
$$
(a_{ij})=B\;{}^t\hskip -1pt B^{(3)}
$$
holds, where $B^{(3)}$ is obtained from $B$ by applying $x\mapsto x^3$ to the entries,
and put
$$
(H\spprime_{0}, H\spprime_{1}, H\spprime_{2},  H\spprime_{3})=(H\sprime_{0}, H\sprime_{1}, H\sprime_{2},  H\sprime_{3})B.
$$
Then $H\spprime_{0}, H\spprime_{1}, H\spprime_{2},  H\spprime_{3}$ satisfy
$$
H\sp{\pprime4}_{0}+ H\sp{\pprime4}_{1}+ H\sp{\pprime4}_{2}+  H\sp{\pprime4}_{3}\equiv 0
\;\bmod (w^4+x^4+y^4+1).
$$
Therefore $(w,x,y)\mapsto [H\spprime_{0}: H\spprime_{1}: H\spprime_{2}:  H\spprime_{3}]$
induces an automorphism $g\spprime$ of $X$.
Using the method described in the proof of Theorem~\ref{thm:H},  
we calculate the matrix $A\spprime$ such that
the action $g\spprime_*$ of $g\spprime$ on $S$ is given by $v\mapsto vA\spprime$.
Next we search for $\tau\in \PGU_4(\F_9)$ such that
 $A\spprime T_{\tau}$ is equal to $A_1$,
 where $T_\tau\in \OG^+(S)$ is the matrix representation of $\tau$.
Then the polynomials
$$
(H_{10}, H_{11}, H_{12},  H_{13}):=(H\spprime_{0}, H\spprime_{1}, H\spprime_{2},  H\spprime_{3})\tau
$$
have the required property.
\end{remark}
\begin{remark}
We have calculated the images of the $\F_9$-rational points of $X$
by the morphisms $\psi_i: X\to Y_i$ and $g_i: X\to X$,
and confirmed that they are compatible~(see~\cite{shimadawebpage}).
\end{remark}
\section{Generators of $\OG^+(S)$}\label{sec:OS}
%
%
Let $F\in \OG^+(S)$ denote the isometry of $S$ obtained from 
the Frobenius action  $\Fr$ of $\F_9$ over $\F_3$ on $X$.
Calculating  the action of $\Fr$ on the lines 
($\ell_{1}^{\Fr}=\ell_{6}, \ell_{2}^{\Fr}=\ell_{5}, \ell_{3}^{\Fr}=\ell_{8},  \ell_{4}^{\Fr}=\ell_{7}, \dots$), 
we see that 
 $F$ is given $v\mapsto vA_F$,
where $A_F$ is the matrix presented  in Table~\ref{table:Frob}. 
\begin{table}
$$
\left[
\begin{array}{cccccccccccccccccccccc}
 0 \shrink & 0 \shrink & 0 \shrink & 0 \shrink & 0 \shrink & 1 \shrink & 0 \shrink & 0 \shrink & 0 \shrink & 0 \shrink & 0 \shrink & 0 \shrink & 0 \shrink & 0 \shrink & 0 \shrink & 0 \shrink & 0 \shrink & 0 \shrink & 0 \shrink & 0 \shrink & 0 \shrink & 0 \\ 
\noalign{\myvskipmat} 0 \shrink & 0 \shrink & 0 \shrink & 0 \shrink & 1 \shrink & 0 \shrink & 0 \shrink & 0 \shrink & 0 \shrink & 0 \shrink & 0 \shrink & 0 \shrink & 0 \shrink & 0 \shrink & 0 \shrink & 0 \shrink & 0 \shrink & 0 \shrink & 0 \shrink & 0 \shrink & 0 \shrink & 0 \\ 
\noalign{\myvskipmat} 1 \shrink & 1 \shrink & 1 \shrink & 1 \shrink & -1 \shrink & -1 \shrink & -1 \shrink & 0 \shrink & 0 \shrink & 0 \shrink & 0 \shrink & 0 \shrink & 0 \shrink & 0 \shrink & 0 \shrink & 0 \shrink & 0 \shrink & 0 \shrink & 0 \shrink & 0 \shrink & 0 \shrink & 0 \\ 
\noalign{\myvskipmat} 0 \shrink & 0 \shrink & 0 \shrink & 0 \shrink & 0 \shrink & 0 \shrink & 1 \shrink & 0 \shrink & 0 \shrink & 0 \shrink & 0 \shrink & 0 \shrink & 0 \shrink & 0 \shrink & 0 \shrink & 0 \shrink & 0 \shrink & 0 \shrink & 0 \shrink & 0 \shrink & 0 \shrink & 0 \\ 
\noalign{\myvskipmat} 0 \shrink & 1 \shrink & 0 \shrink & 0 \shrink & 0 \shrink & 0 \shrink & 0 \shrink & 0 \shrink & 0 \shrink & 0 \shrink & 0 \shrink & 0 \shrink & 0 \shrink & 0 \shrink & 0 \shrink & 0 \shrink & 0 \shrink & 0 \shrink & 0 \shrink & 0 \shrink & 0 \shrink & 0 \\ 
\noalign{\myvskipmat} 1 \shrink & 0 \shrink & 0 \shrink & 0 \shrink & 0 \shrink & 0 \shrink & 0 \shrink & 0 \shrink & 0 \shrink & 0 \shrink & 0 \shrink & 0 \shrink & 0 \shrink & 0 \shrink & 0 \shrink & 0 \shrink & 0 \shrink & 0 \shrink & 0 \shrink & 0 \shrink & 0 \shrink & 0 \\ 
\noalign{\myvskipmat} 0 \shrink & 0 \shrink & 0 \shrink & 1 \shrink & 0 \shrink & 0 \shrink & 0 \shrink & 0 \shrink & 0 \shrink & 0 \shrink & 0 \shrink & 0 \shrink & 0 \shrink & 0 \shrink & 0 \shrink & 0 \shrink & 0 \shrink & 0 \shrink & 0 \shrink & 0 \shrink & 0 \shrink & 0 \\ 
\noalign{\myvskipmat} 1 \shrink & 0 \shrink & 1 \shrink & 1 \shrink & 0 \shrink & -1 \shrink & 0 \shrink & 0 \shrink & -1 \shrink & 0 \shrink & 0 \shrink & 0 \shrink & 0 \shrink & 0 \shrink & 0 \shrink & 0 \shrink & 0 \shrink & 0 \shrink & 0 \shrink & 0 \shrink & 0 \shrink & 0 \\ 
\noalign{\myvskipmat} 0 \shrink & 1 \shrink & 1 \shrink & 1 \shrink & -1 \shrink & 0 \shrink & 0 \shrink & -1 \shrink & 0 \shrink & 0 \shrink & 0 \shrink & 0 \shrink & 0 \shrink & 0 \shrink & 0 \shrink & 0 \shrink & 0 \shrink & 0 \shrink & 0 \shrink & 0 \shrink & 0 \shrink & 0 \\ 
\noalign{\myvskipmat} -1 \shrink & -1 \shrink & -1 \shrink & -2 \shrink & 1 \shrink & 1 \shrink & 1 \shrink & 1 \shrink & 1 \shrink & 1 \shrink & 0 \shrink & 0 \shrink & 0 \shrink & 0 \shrink & 0 \shrink & 0 \shrink & 0 \shrink & 0 \shrink & 0 \shrink & 0 \shrink & 0 \shrink & 0 \\ 
\noalign{\myvskipmat} 2 \shrink & 1 \shrink & 1 \shrink & 2 \shrink & -1 \shrink & -2 \shrink & -1 \shrink & -1 \shrink & -1 \shrink & 0 \shrink & 1 \shrink & 0 \shrink & 0 \shrink & 0 \shrink & -1 \shrink & 0 \shrink & 0 \shrink & 0 \shrink & 0 \shrink & 1 \shrink & 0 \shrink & 0 \\ 
\noalign{\myvskipmat} 0 \shrink & 0 \shrink & 0 \shrink & 0 \shrink & 1 \shrink & 1 \shrink & 0 \shrink & 1 \shrink & 1 \shrink & 0 \shrink & -1 \shrink & 0 \shrink & 0 \shrink & -1 \shrink & 0 \shrink & 0 \shrink & 0 \shrink & 0 \shrink & 0 \shrink & -1 \shrink & 0 \shrink & 0 \\ 
\noalign{\myvskipmat} 2 \shrink & 2 \shrink & 2 \shrink & 2 \shrink & -1 \shrink & 0 \shrink & -1 \shrink & -1 \shrink & 0 \shrink & -1 \shrink & 0 \shrink & 0 \shrink & -1 \shrink & 0 \shrink & 0 \shrink & -1 \shrink & 0 \shrink & 0 \shrink & 0 \shrink & 0 \shrink & -1 \shrink & 0 \\ 
\noalign{\myvskipmat} 0 \shrink & 1 \shrink & 1 \shrink & 0 \shrink & 0 \shrink & 1 \shrink & 1 \shrink & 0 \shrink & 0 \shrink & 0 \shrink & -1 \shrink & -1 \shrink & 0 \shrink & 0 \shrink & 0 \shrink & 0 \shrink & 0 \shrink & 0 \shrink & 0 \shrink & -1 \shrink & 0 \shrink & 0 \\ 
\noalign{\myvskipmat} 0 \shrink & 0 \shrink & 0 \shrink & 0 \shrink & 0 \shrink & 0 \shrink & 0 \shrink & 0 \shrink & 0 \shrink & 0 \shrink & 0 \shrink & 0 \shrink & 0 \shrink & 0 \shrink & 0 \shrink & 0 \shrink & 0 \shrink & 0 \shrink & 0 \shrink & 1 \shrink & 0 \shrink & 0 \\ 
\noalign{\myvskipmat} 0 \shrink & 0 \shrink & 0 \shrink & 0 \shrink & 0 \shrink & 0 \shrink & 0 \shrink & 0 \shrink & 0 \shrink & 0 \shrink & 0 \shrink & 0 \shrink & 0 \shrink & 0 \shrink & 0 \shrink & 0 \shrink & 0 \shrink & 0 \shrink & 0 \shrink & 0 \shrink & 1 \shrink & 0 \\ 
\noalign{\myvskipmat} 2 \shrink & 1 \shrink & 1 \shrink & 2 \shrink & 0 \shrink & 0 \shrink & -1 \shrink & 0 \shrink & 0 \shrink & -1 \shrink & 0 \shrink & 1 \shrink & 0 \shrink & -1 \shrink & 0 \shrink & -1 \shrink & -1 \shrink & 0 \shrink & 0 \shrink & 0 \shrink & -1 \shrink & 0 \\ 
\noalign{\myvskipmat} 0 \shrink & 0 \shrink & 0 \shrink & 0 \shrink & 0 \shrink & 0 \shrink & 1 \shrink & 0 \shrink & -1 \shrink & 0 \shrink & 0 \shrink & -1 \shrink & 0 \shrink & 1 \shrink & 0 \shrink & 1 \shrink & 0 \shrink & -1 \shrink & 0 \shrink & 0 \shrink & 1 \shrink & 0 \\ 
\noalign{\myvskipmat} -3 \shrink & -2 \shrink & -2 \shrink & -3 \shrink & 1 \shrink & 1 \shrink & 1 \shrink & 1 \shrink & 1 \shrink & 1 \shrink & 0 \shrink & 0 \shrink & 1 \shrink & 0 \shrink & 1 \shrink & 1 \shrink & 1 \shrink & 1 \shrink & 1 \shrink & -1 \shrink & 0 \shrink & 0 \\ 
\noalign{\myvskipmat} 0 \shrink & 0 \shrink & 0 \shrink & 0 \shrink & 0 \shrink & 0 \shrink & 0 \shrink & 0 \shrink & 0 \shrink & 0 \shrink & 0 \shrink & 0 \shrink & 0 \shrink & 0 \shrink & 1 \shrink & 0 \shrink & 0 \shrink & 0 \shrink & 0 \shrink & 0 \shrink & 0 \shrink & 0 \\ 
\noalign{\myvskipmat} 0 \shrink & 0 \shrink & 0 \shrink & 0 \shrink & 0 \shrink & 0 \shrink & 0 \shrink & 0 \shrink & 0 \shrink & 0 \shrink & 0 \shrink & 0 \shrink & 0 \shrink & 0 \shrink & 0 \shrink & 1 \shrink & 0 \shrink & 0 \shrink & 0 \shrink & 0 \shrink & 0 \shrink & 0 \\ 
\noalign{\myvskipmat} 4 \shrink & 3 \shrink & 3 \shrink & 4 \shrink & -1 \shrink & -1 \shrink & -1 \shrink & -1 \shrink & -1 \shrink & -1 \shrink & -1 \shrink & -1 \shrink & -1 \shrink & -1 \shrink & -1 \shrink & -1 \shrink & -1 \shrink & -1 \shrink & 0 \shrink & 0 \shrink & 0 \shrink & 1 
\end{array}
\right]
$$
\vskip 5pt
\caption{Frobenius action on $S$}\label{table:Frob}
\end{table}
Since $h_0^F=h_0$,
 we have $\chamsz^F=\chamsz$ by Corollary~\ref{cor:equiv}.
\begin{proposition}\label{prop:autchamber}
The automorphism group $\Aut(\chamsz)\subset \OG^+(S)$ of 
$\chamsz$ is  the  split extension of $\gen{F}\cong \Z/2\Z$ by  $\Aut(X, h_0)$.
\end{proposition}
\begin{proof}
Since we have calculated the representation~\eqref{eq:Ttau} of $\Aut(X, h_0)$
into $\OG^+(S)$,
we can verify that $F\notin \Aut(X, h_0)$.
Therefore it is enough to show that the order of $\Aut(\chamsz)$ is equal to $2$ times $|\Aut(X, h_0)|$.
Since $|\PGU_4(\F_9)|$ is equal to $4$ times $|\PSU_4(\F_9)|$,
this follows from~\cite[Lemma 2.1]{KatsuraKondochar3} (see also~\cite[p.~52]{MR827219}).
\end{proof}
Since $\intS{[\ell_1], [\ell_1]}=-2$, 
the reflection $s_1: S\tR\to S\tR$ into  
the hyperplane $([\ell_1])_S\sperp$  
is contained in $\OG^+(S)$.
In the same way as the proof of Theorem~\ref{thm:main},
we obtain the following:
\begin{theorem}\label{thm:OplusS} 
The  autochronous orthogonal group $\OG^+(S)$ 
of the N\'eron-Severi lattice $S$ of $X$ is generated by $\Aut(X, h_0)=\PGU_4(\F_9)$,  
$g_1$, $g_2$, $F$ and $s_1$.
\end{theorem}
\bibliographystyle{plain}
\def\cftil#1{\ifmmode\setbox7\hbox{$\accent"5E#1$}\else
  \setbox7\hbox{\accent"5E#1}\penalty 10000\relax\fi\raise 1\ht7
  \hbox{\lower1.15ex\hbox to 1\wd7{\hss\accent"7E\hss}}\penalty 10000
  \hskip-1\wd7\penalty 10000\box7} \def\cprime{$'$} \def\cprime{$'$}
  \def\cprime{$'$} \def\cprime{$'$}

\end{document}